\documentstyle[12pt]{article}
    \title{{\bf  Semi-infinite forms and 
topological vertex operator algebras}}
    \author{Yi-Zhi Huang and Wenhua Zhao}
   \date{}
    \begin{document}
    \bibliographystyle{alpha}
    \maketitle

    \input amssym.def
    \input amssym
    \newtheorem{rema}{Remark}[section]
    \newtheorem{propo}[rema]{Proposition}
    \newtheorem{theo}[rema]{Theorem}
 \newtheorem{conj}[rema]{Conjecture}
   \newtheorem{defi}[rema]{Definition}
    \newtheorem{lemma}[rema]{Lemma}
    \newtheorem{corol}[rema]{Corollary}
     \newtheorem{exam}[rema]{Example}
 \newtheorem{warn}[rema]{Warning}
\newcommand{\binom}[2]{{{#1}\choose {#2}}}
	\newcommand{\nno}{\nonumber}
	\newcommand{\lbar}{\bigg\vert}
\newcommand{\mbar}{\mbox{\large $\vert$}}
	\newcommand{\p}{\partial}
	\newcommand{\dps}{\displaystyle}
	\newcommand{\bra}{\langle}
	\newcommand{\ket}{\rangle}
\newcommand{\kr}{\mbox{\rm Ker}\ }
 \newcommand{\res}{\mbox{\rm Res}}
\renewcommand{\hom}{\mbox{\rm Hom}}
 \newcommand{\pf}{{\it Proof.}\hspace{2ex}}
\newcommand{\epfd}{\hspace{1em}\Box}
 \newcommand{\epf}{\hspace{1em}$\Box$}
 \newcommand{\epfv}{\hspace{1em}$\Box$\vspace{1em}}
\newcommand{\nord}{\mbox{\scriptsize ${\circ\atop\circ}$}}
\newcommand{\wt}{\mbox{\rm wt}\ }
\newcommand{\swt}{{\mbox{\rm \scriptsize wt}}\;}
\newcommand{\clr}{\mbox{\rm clr}\ }
\newcommand{\ideg}{\mbox{\rm Ideg}\ }

\begin{abstract}
Semi-infinite forms on the moduli spaces of genus-zero Riemann
surfaces with punctures and local coordinates are introduced. A
partial operad for semi-infinite forms is constructed.  Using
semi-infinite forms and motivated by a partial suboperad of the
partial operad for semi-infinite forms, topological vertex partial
operads of type $k<0$ and strong topological vertex partial operads of
type $k<0$ are constructed.  It is proved that the category of
(locally-)grading-restricted 
(strong) topological vertex operator algebras of type $k<0$ and the
category of (weakly) meromorphic ${\Bbb Z}\times {\Bbb Z}$-graded 
algebras over
the (strong) topological vertex partial operad of type $k$ are
isomorphic.  As an application of this isomorphism theorem, the
following conjecture of Lian-Zuckerman and Kimura-Voronov-Zuckerman is
proved: A strong topological vertex operator algebra gives a homotopy
Gerstenhaber algebra.  These results hold in particular for  the tensor
product of the moonshine module vertex operator algebra, the vertex
algebra constructed {from} a rank $2$ Lorentz lattice and the ghost
vertex operator algebra, studied in detail first by Lian and Zuckerman.
\end{abstract}

\tableofcontents

\renewcommand{\theequation}{\thesection.\arabic{equation}}
\renewcommand{\therema}{\thesection.\arabic{rema}}
\setcounter{equation}{0}
\setcounter{rema}{0}
\setcounter{section}{0}

\section{Introduction}

In the present paper, we give a geometric and operadic formulation of
the notion of topological vertex 
operator algebra (satisfying certain additional
axioms)
and as a consequence,
we prove a conjecture of Lian and Zuckerman \cite{LZ2} and Kimura,
Voronov and Zuckerman \cite{KVZ}. This geometric and operadic
formulation is an application of the geometric theory of vertex operator
algebras developed in \cite{H-1}, \cite{H0} and \cite{H} and the theory
of semi-infinite forms for the Virasoro algebra introduced by Feigin
\cite{F} and developed by Frenkel, Garland and Zuckerman \cite{FGZ} and
by Lian and Zuckerman \cite{LZ1}. These results will be
useful in the future mathematical study of string backgrounds and $N=2$
superconformal field theories. 

A notion of topological chiral algebra was first introduced by Lian
 and Zuckerman \cite{LZ2}. Following the mathematical terminology,
 these algebras are later called topological vertex operator algebras
 by several authors (though there is a very subtle difference between
 topological chiral algebras in the sense of \cite{LZ2} 
(or topological vertex operator algebras in the sense of 
\cite{KVZ}) and topological vertex operator algebras
 in the sense of \cite{H2}; see below).  
In \cite{LZ2}, Lian and Zuckerman showed that the cohomology of a
(strong) topological chiral algebra has a natural structure of a
Batalin-Vilkovisky (``coboundary Gerstenhaber algebra'') structure.
Based on calculations, they also conjectured that there must be a
certain ``homotopy Gerstenhaber algebra'' structure on a topological
chiral algebra.  Identities satisfied by the corresponding operations in
topological vertex operator algebras were studied by Akman \cite{A}.
After attempts by different authors, there is now a
 precise notion of homotopy Gerstenhaber algebra (see \cite{V}). In
 \cite{KVZ}, Kimura, Voronov and Zuckerman developed a framework for
 the construction of homotopy Gerstenhaber algebras from genus-zero
 topological conformal field theories. Actually the notion of homotopy
 Gerstenhaber algebra used in \cite{KVZ} is not the correct one
 because the construction in \cite{GJ} justifying the definition
 contained an error. But appropriate corrections can be made so that
 with the correct definition of homotopy Gerstenhaber algebra, the
 statement in \cite{KVZ} that a genus-zero topological conformal field
 theory in the sense of Segal \cite{S} and Getzler \cite{G} has a
 structure of a homotopy Gerstenhaber algebra is still correct. See
 \cite{V} for more details. With the precise definition of homotopy
Gerstenhaber algebra, the conjecture of Lian-Zuckerman can be 
formulated  precisely as follows (cf. \cite{KVZ}):

\begin{conj}\label{conj1}
A topological vertex operator algebra gives
a homotopy Gerstenhaber algebra.
\end{conj}

Since it was already proved in \cite{KVZ} and \cite{V} that 
a genus-zero topological conformal field theory gives a 
homotopy Gerstenhaber algebra, it is natural to try to prove Conjecture 
\ref{conj1} by proving
the following conjecture
also by Kimura, Voronov and Zuckerman in \cite{KVZ}:

\begin{conj}\label{conj2}
An appropriate topological completion of a topological vertex operator algebra
has a structure of a genus-zero topological conformal field theory. 
\end{conj}

There is a very subtle issue which we must address here.  In many works
on vertex operator algebras, including \cite{FLM}, \cite{FHL}, \cite{H2}
and \cite{H}, vertex operator algebras or topological vertex operator
algebras are required to satisfy the following grading-restriction
conditions: For a vertex operator algebra $V=\coprod_{n\in {\Bbb
Z}}V_{(n)}$, $\dim V_{(n)}<\infty$ for $n\in {\Bbb Z}$ and $V_{(0)}=0$
for $n$ sufficiently small. But the topological chiral algebras or
topological vertex operator algebras discussed by Lian and Zuckerman in
\cite{LZ2} and by Kimura, Voronov and Zuckerman in \cite{KVZ} in general
do not satisfy the grading-restriction conditions. In fact, the
important examples {from} string backgrounds, including the tensor
product of the moonshine module vertex operator algebra, the vertex
algebra constructed {from} a rank $2$ Lorentz lattice and the ghost
vertex operator algebra (see \cite{LZ3}), in general do not satisfy
these conditions.  

In the present paper, an algebra satisfying the original definition of
Lian and Zuckerman in \cite{LZ2} is called a topological vertex
operator algebra, as is in \cite{KVZ}. (We warn the reader that a
topological vertex operator algebra in this sense is in general not a
vertex operator algebra in the sense of \cite{FLM} and \cite{FHL}
because the grading-restriction conditions are not satisfied.) A
topological vertex operator algebra in the sense of \cite{H2} is
called a grading-restricted topological vertex operator algebra. We
also introduce local grading-restriction conditions.
Grading-restricted topological vertex operator algebras and the
examples {from} string backgrounds are locally-grading-restricted
topological vertex operator algebras. We also introduce the notion of
strong topological vertex operator algebra for which the square of a
certain operator on the algebra is $0$.  We prove Conjecture
\ref{conj2} for a local-grading-restricted strong topological vertex
operator algebra as a consequence of the main theorem of the present
paper stating that a geometric formulation of the notion of
locally-grading-restricted (strong) topological vertex operator
algebra is equivalent to the algebraic formulation.  As a consequence,
Conjecture \ref{conj1} is true for a locally-grading-restricted strong
topological vertex operator algebra.  In particular, Conjecture
\ref{conj1} and Conjecture \ref{conj2} are true for grading-restricted
strong topological vertex operator algebra and for the tensor product
of the moonshine module vertex operator algebra, the vertex algebra
constructed {from} a rank $2$ Lorentz lattice and the ghost vertex
operator algebra (see \cite{LZ3}).  In the case that the topological
vertex operator algebra is not strong, the same construction does not
give a simple geometric structure.

The geometric theory  of vertex operator algebras was initialized by Frenkel
\cite{Fr}.
A geometric and operadic formulation of the notion of 
vertex operator algebra in the sense of \cite{FLM} and \cite{FHL} was
given  by the first author in \cite{H-1}, \cite{H0} and \cite{H} (see
also 
\cite{H1}) and by 
Lepowsky and the first 
author in \cite{HL1} and \cite{HL2}. The 
 corresponding isomorphism theorem was proved by the first author 
in \cite{H-1}, \cite{H0} and \cite{H}. 

In \cite{H2}, the first author gave a
geometric and operadic formulation of the notion of topological vertex
algebra (which may not have a Virasoro element) and an isomorphism
theorem was proved. This isomorphism theorem for topological vertex algebras
combined with a theorem of Cohen \cite{C} or Getzler \cite{G}
gives a geometric construction of the Gerstenhaber or
Batalin-Vilkovisky algebra structure on the cohomology of 
a topological vertex algebra (see \cite{H2}). 

It is natural to look for a geometric formulation of the notion of
grading-restricted topological vertex operator algebra, and more
generally, of the notion of locally-grading-restricted topological
vertex operator algebra.  Because topological vertex operator algebras
and locally-grading-restricted topological vertex operator algebras
have Virasoro elements and elements which give differentials and
fermion gradings, the geometry underlying these algebras is certainly
much more complicated than that underlying topological vertex
algebras.

In the present paper, we find that the correct geometric objects
underlying   locally-grading-restricted  topological 
vertex operator algebras  are 
``universal coverings'' of a partial operad constructed {from} what we call
``semi-infinite forms on the moduli spaces of genus-zero Riemann
surfaces with punctures and local coordinates.''  Semi-infinite forms
and semi-infinite cohomologies for graded Lie algebras were introduced
by Feigin \cite{F} and developed in the context of string theory by
Frenkel, Garland and Zuckerman \cite{FGZ} and further by Lian and Zuckerman
\cite{LZ1}. Using semi-infinite forms for the Virasoro algebra, we
construct semi-infinite forms on the moduli spaces above and a partial
operad for semi-infinite forms.  Using these semi-infinite
forms and motivated by a partial
suboperad of the partial operad for semi-infinite forms, we construct
partial operads  called ``topological vertex partial
operad of type $k$'' and  partial operads  called 
``strong topological vertex partial
operad of type $k$'' for $k<0$. We also introduce the notions of (strong)
topological vertex operator algebra of type $k$ and (strong)
topological 
vertex operator algebras  of type $k$ 
 for $k<0$.
The main theorem of the present paper states that for any $k<0$,
the category of meromorphic ${\Bbb Z}\times {\Bbb Z}$-graded algebras 
over the (strong) topological vertex partial operad of type $k$
is isomorphic to the category of (strong) topological
vertex operator algebras of type $k$, and 
the category of weakly meromorphic ${\Bbb Z}\times {\Bbb Z}$-graded algebras 
over the (strong) topological vertex partial operad of type $k$
is isomorphic to the category of local-grading-restricted (strong) 
topological 
vertex operator algebras  of type $k$.
The proof of  Conjecture \ref{conj2} for 
local-grading-restricted strong  topological 
vertex operator algebras
 follows easily because there is a morphism
{from} the operad for genus-zero topological
conformal field theories to  strong
topological vertex partial operads.

The material in this paper depends heavily on the monograph \cite{H}. 

The present paper is organized as follows: In Section 2, we review 
basic concepts and state the main results, Theorems \ref{theo2},
\ref{theo1} and \ref{main},
 of the present paper. The notions of topological vertex operator 
algebra in the sense of \cite{LZ2} and in the sense of \cite{H2} 
and related notions are recalled in Subsection 2.1. 
In Subsection
2.2, after recalling the notions of Gerstenhaber algebra, homotopy 
Gerstenhaber algebra and genus-zero topological 
conformal field theory and a result of Kimura-Voronov-Zuckerman and Voronov, 
we state our solutions, Theorems \ref{theo2} and \ref{theo1}, 
to  Conjectures \ref{conj2} and 
\ref{conj1}, respectively, for 
local-grading-restricted strong topological 
vertex operator algebras 
and in particular for grading-restricted strong topological 
vertex operator algebras in the sense of \cite{H2}. In Subsection 2.3, the 
main theorem, Theorem \ref{main}, of the present paper is stated.
 In Section 3, we construct and study semi-infinite forms on 
moduli spaces of spheres with tubes. In Subsection 3.1,  we construct
left invariant meromorphic vector fields on $\tilde{K}^{c}(1)$ (see
\cite{H}) satisfying the Virasoro relations. They are needed in the
construction of semi-infinite forms. Semi-infinite forms on $K(0)$
(see \cite{H}) are constructed and studied in Subsection 3.2. In 
Subsections 3.3 and 3.4,
semi-infinite forms on $K(n)$ for $n\ge 0$ and a partial operad ${\frak
G}$ for semi-infinite forms, respectively,
are constructed. In Section 4, we construct the 
(strong) topological vertex partial operad and prove the main results
stated in Subsections 2.2 and 2.3. In Subsection 4.1, some properties of
topological  vertex operator algebras 
 are proved. The topological vertex partial operad and the strong
topological vertex partial operad are constructed in Subsection 4.2.
The main theorem, Theorem \ref{main}, of the present paper is 
proved in Subsection
4.3 and Theorem \ref{theo2} is proved in Subsection 4.4. 
 Section 5 is an appendix, in which 
two types of examples of strong topological vertex
operator algebras and  local-grading-restricted strong topological 
vertex operator algebras
are given. The examples 
obtained by tensoring the ghost vertex operator algebra with 
vertex algebras of central charge $26$ are given in 
Subsection 5.1. The examples obtained by twisting $N=2$ superconformal 
vertex operator superalgebras are given in Subsection 5.2.

\paragraph{Acknowledgment} We are grateful to Sasha Voronov for
explaining to us his correction of the main result in \cite{KVZ} and
for sending us a preliminary version of \cite{V}. We also thank 
Gregg Zuckerman and Bong H. Lian for comments and suggestioins.
Y.Z.H. is supported
in part by by NSF grant DMS-9622961. W.Z. is grateful to the
Department of Mathematics at University of Chicago for the financial
supports.

\renewcommand{\theequation}{\thesection.\arabic{equation}}
\renewcommand{\therema}{\thesection.\arabic{rema}}
\setcounter{equation}{0}
\setcounter{rema}{0}

\section{The main results}

\subsection{Topological vertex operator algebras}

We assume that the reader is familiar with the precise
notions of vertex algebra and
vertex operator algebra (see, for example, \cite{Bor}, 
\cite{FLM}, \cite{FHL}). In particular,  a vertex operator 
algebra $V=\coprod_{n\in {\Bbb Z}}V_{(n)}$ in this paper 
is ${\Bbb Z}$-graded and satisfies the following grading-restriction 
conditions: $\dim V_{(n)}<\infty$ for $n\in {\Bbb Z}$ and $V_{(n)}=0$ for 
$n$ sufficiently small.
In the present paper, we need the following
variants (see, for example, \cite{H2}; cf. \cite{DL}):

\begin{defi}
{\rm A {\it ${\Bbb Z}\times {\Bbb Z}$-graded 
vertex operator algebra without grading restrictions} is a ${\Bbb Z}\times 
{\Bbb Z}$-graded vector space $V$ 
(graded by {\it weights} and by {\it fermion numbers}), that is,
$$V=\coprod_{m, n\in {\Bbb Z}}V_{(n)}^{(m)}=\coprod_{n\in {\Bbb Z}}V_{(n)}
=\coprod_{m\in {\Bbb Z}}V^{(m)}
$$
where 
$$
V_{(n)}=\coprod_{m\in {\Bbb Z}}V_{(n)}^{(m)}, \;\;\;
V^{(m)}=\coprod_{n\in {\Bbb Z}}V_{(n)}^{(m)},
$$
equipped with a vertex operator 
map $Y: V\otimes V\rightarrow V[[x, x^{-1}]]$ which maps
$V^{(m_{1})}\otimes V^{(m_{2})}$ to $V^{(m_{1}+m_{2})}[[x, x^{-1}]]$,
a vacuum ${\bf 1}$ and a Virasoro element $\omega$, satisfying the 
following axioms: 

\begin{enumerate}

\item For $v\in V^{(m)}$, we say that $v$ has fermion number $m$ and
 use $|v|$ to denote $m$. Then for the vacuum
${\bf 1}$ and the Virasoro element $\omega$, $|{\bf 1}|=|\omega|=0$.

\item 
For $u, v\in V$, the (Cauchy-)Jacobi identity
\begin{eqnarray*}
\lefteqn{x_{0}^{-1}\delta 
\left({\displaystyle\frac{x_{1}-x_{2}}{x_{0}}}\right)Y(u, 
x_{1})Y(v, x_{2})}\nno\\
&&\hspace{6em}-(-1)^{|u||v|}x_{0}^{-1} \delta \left({\displaystyle
\frac{x_{2}-x_{1}}{-x_{0}}}\right)Y(v, x_{2})Y(u, x_{1})\nno \\
&&=x_{2}^{-1} \delta \left({\displaystyle\frac{x_{1}-x_{0}}{x_{2}}}\right)
Y(Y(u,
x_{0})v, x_{2})
\end{eqnarray*}
holds.

\item All the other  axioms for 
vertex operator algebras, except the grading-restriction conditions,
hold. 

\end{enumerate}}

\end{defi}

The following notion of  topological vertex operator algebra 
is  the same as the notion of topological chiral algebra 
introduced by Lian and Zuckerman in \cite{LZ2}:

\begin{defi}\label{tvoa}
{\rm A {\it topological vertex operator algebra}
is a  ${\Bbb Z}\times {\Bbb Z}$-graded vertex operator algebra 
without grading restrictions
$(V, Y, {\bf 1}, \omega)$ 
equipped  with three additional distinguished elements $f\in
V_{(1)}$, $q\in V_{(1)}^{(1)}$ and $g\in V_{(2)}^{(1)}$ satisfying the
following axioms:

\begin{enumerate}

\item  Let $Y(f, x)=\sum_{n\in {\Bbb Z}}f_{n}x^{-n-1}$. Then for any $v\in
V^{(m)}$, 
\begin{equation}\label{5.1}
f_{0}v=mv.
\end{equation}

\item  Let $Y(q, x)=\sum_{n\in {\Bbb Z}}q_{n}x^{-n-1}$ and $Q=q_{0}$. Then
\begin{eqnarray*}
L(n)q&=&0, \;\;\;n>0,\\
Q^{2}&=&0.
\end{eqnarray*}

\item  Let $Y(g, x)=\sum_{n\in {\Bbb Z}}g(n)x^{-n-2}$ and $\omega$ the
Virasoro element of $V$. Then
\begin{eqnarray*}
L(n)g&=&0, \;\;\;n>0,\\
Qg&=&\omega.
\end{eqnarray*}
\end{enumerate}

A topological vertex operator algebra is  {\it strong}
if it  satisfies the condition $g(0)^{2}=0$.

A topological vertex operator algebra is {\it 
grading-restricted} if the grading-restriction conditions are satisfied.

A topological vertex operator algebra is {\it locally
grading-restricted} if the following additional axioms are satisfied:

\begin{enumerate}

\item For any element of the topological vertex
operator algebra, the module $W=\coprod_{n\in {\Bbb Z}}W_{(n)}$ for the
Virasoro algebra generated by this element satisfies the
grading-restriction conditions, that is, $\dim W_{(n)}<\infty$ for $n\in
{\Bbb Z}$ and $W_{(n)}=0$ for $n$ sufficiently small. 

\item The vertex subalgebra generated by $\omega, q, f, g$
satisfies the grading-re\-striction conditions. 

\end{enumerate}

Let $k$ be a negative integer. 
A  locally-grading-restricted (strong) topological 
vertex operator algebra is {\it 
of type $k$} if  the weights of the homogeneous nonzero
elements of its  vertex subalgebra generated by 
$\omega$, $q$, $f$ and $g$ are larger than or equal to $k$.

{\it Homomorphisms} and {\it isomorphisms} of topological vertex
operator algebras are defined in the obvious way.}
\end{defi}

We denote the topological vertex operator algebra just defined by 
$$(V, Y, {\bf 1}, \omega, f, q, g)$$ 
or simply by $V$. 
Note that (\ref{5.1}) implies $f\in V_{(1)}^{(0)}$. Also it is clear that 
if a (strong) topological vertex operator algebra is of type $k$, it is also
of type $k+m$ for any $m<0$.

\begin{rema}\label{term}
{\rm We warn the reader that topological vertex operator algebras
defined above in general are not ${\Bbb Z}\times {\Bbb Z}$-graded
vertex operator algebras because the grading-restriction conditions
might not be satisfied. Topological vertex operator algebras defined
in \cite{H2} are grading-restricted topological vertex operator
algebras defined above. The reason why we use the terminology above in
the present paper is to make our terminology agree with those in
\cite{LZ2} and \cite{KVZ}.  We prefer to call an algebra satisfying
all the axioms for vertex operator operator algebras except the
grading-restriction axioms a {\it conformal vertex algbera}. Then a
vertex operator algebra is a {\it grading-restricted conformal vertex
algebra} and a topological vertex operator algebra in the sense of
Definition \ref{tvoa} is called a {\it topological conformal vertex
algebra}. Similarly, other notions can be defined using conformal
vertex algebras. For example, a {\it (${\Bbb Z}\times {\Bbb
Z}$-graded) locally-grading-restricted conformal vertex algebra} is a
(${\Bbb Z}\times {\Bbb Z}$-graded) vertex operator algebra without
grading restrictions but satisfying the local grading-restriction
conditions.  For simplicity, this terminology will be used in
Subsections 4.3 and 5.1.}
\end{rema}

In the appendix (Section 5), we give two types of examples of
(locally-)grading-restricted topological vertex operator algebras: The ones
obtained by tensoring the ghost vertex operator algebra with 
vertex operator
algebras without grading restrictions 
of central charge $26$ satisfying the local grading restriction conditions,
and the ones obtained by
twisting $N=2$ superconformal vertex operator algebras. These examples
are the main interesting 
(locally-)grading-restricted topological 
vertex operator algebras. They
are all strong in the sense that $g^{2}(0)=0$.

\subsection{Homotopy Gerstenhaber algebras and topological
conformal field theories}

\begin{defi}
{\rm A {\it Gerstenhaber algebra} is a graded commutative
algebra $A$ together with a bracket
$[\cdot, \cdot]: A\otimes A\to A$ such that $[A^{(n)}, A^{(m)}]\subset
A^{(n+m-1)}$ (where $A^{(n)}$, $A^{(m)}$, $n, m\in {\Bbb Z}$, are homogeneous
components of $A$), satisfying 
\begin{eqnarray*}
{[a, b]}&=&-(-1)^{(|a|-1)(|b|-1)}[b, a],\\ 
{[a, [b, c]]}&=&[[a,b], c]+(-1)^{(|a|-1)(|b|-1)}[b, [a, c]],\\ 
{[a, bc]}&=&[a, b]c+(-1)^{|a|(|b|-1)}b[a, c]
\end{eqnarray*}
for any homogeneous elements $a,b,c$ of $A$.}
\end{defi}

The definition of homotopy Gerstenhaber algebra is much more
involved. Roughly speaking, a homotopy Gerstenhaber algebra means a
``Gerstenhaber algebra up to homotopy together with homotopies for the
homotopies.'' The precise notion of homotopy Gerstenhaber algebra we
shall use is the one proposed by Voronov in \cite{V}. This definition
needs an operad $E^{1}$ of complexes 
called {\it homotopy Gerstenhaber operad}.  Since the explicit
definition of homotopy Gerstenhaber operad is involved and since
we shall not need the details of this definition in the present paper,
we only mention that $E^{1}$ is 
the first term of the spectral sequence
corresponding to a stratification of a certain moduli space operad
which is homotopy equivalent to the little disks operad. 
See \cite{V} for details.

\begin{defi}
{\rm A {\it homotopy Gerstenhaber algebra} is an algebra over the 
homotopy Gerstenhaber operad $E^{1}$.}
\end{defi}

We need a corrected version by Voronov in \cite{V} of 
a result of Kimura, Voronov and Zuckerman in \cite{KVZ}. 
Recall the suboperad
$K_{{\frak H}_{1}}$ of the partial operad $K$ discussed in Section 6.4
of \cite{H}.  We also need the operad $\wedge TK_{{\frak H}_{1}}$ of complexes
of the exterior
algebras of the tangent bundles of components of $K_{{\frak
H}_{1}}$: For any $n\ge 0$, let $TK_{{\frak H}_{1}}(n)$ be the tangent
bundle of $K_{{\frak H}_{1}}(n)$ and $\wedge TK_{{\frak H}_{1}}(n)$ 
the direct sum of all exterior powers of $TK_{{\frak H}_{1}}(n)$.  Then it
is easy to see that the sequence $\wedge TK_{{\frak H}_{1}}=\{
\wedge TK_{{\frak H}_{1}}(n)\}_{n\ge 0}$ has a natural structure of operad.

Let $(H, Q)$ be a complex. Then the endomorphism operad 
$${\cal E}_{H}=\{\hom(H^{\otimes n}, H)\}_{n\ge 0}$$
has a natural structure of operad of complexes. We shall still use $Q$ 
to denote the differentials on $\hom(H^{\otimes n}, H)$  for $n\ge 0$.
A map $\mu_{n}$ {from} $\wedge TK_{{\frak H}_{1}}(n)$ to
$\hom(H^{\otimes n}, H)$ is said to be {\it fiber-linear} if 
it is linear on the fiber of $\wedge TK_{{\frak H}_{1}}(n)$. 
Then  a fiber-linear smooth map $\mu_{n}$
{from} $\wedge TK_{{\frak H}_{1}}(n)$ to
$\hom(H^{\otimes n}, H)$
can be viewed as a $\hom(H^{\otimes n}, H)$-valued
form on $K_{{\frak H}_{1}}(n)$. In particular, the differential $d$ on 
$\hom(H^{\otimes n}, H)$-valued forms 
on $K_{{\frak H}_{1}}(n)$ acts on $\mu_{n}$.
A morphism $\mu$ {from} $\wedge TK_{{\frak H}_{1}}$ to a smooth 
operad is said to be {\it smooth} if $\mu_{n}$, $n\ge 0$, are smooth and 
is said to be {\it fiber-linear} if $\mu_{n}$, $n\ge 0$, are
fiber-linear.

The following notion is due to Segal \cite{S} and Getzler
\cite{G} (see also \cite{KSV} and \cite{KVZ}):

\begin{defi}
{\rm A {\it genus-zero topological conformal field theory}
is a complex $(H, Q)$
and a fiber-linear smooth morphism $\mu$ of operads
{from} $\wedge TK_{{\frak H}_{1}}$ to 
the endomorphism operad 
$\tilde{\cal E}_{H}$ such that 
\begin{equation}\label{6.1}
d\mu_{n}=Q\mu_{n},
\end{equation}
for $n\ge 0$,
where as discussed above, 
$d$ is the differential on $\hom(H^{\otimes n}, H)$-valued forms 
on $K_{{\frak H}_{1}}(n)$ and $Q$ is the differential on 
$\hom(H^{\otimes n}, H)$
induced {from} $Q: H\to H$.}
\end{defi}

The following corrected version of a result in \cite{KVZ} is 
proved by Voronov in \cite{V}:

\begin{theo}\label{KVZ}
A genus-zero topological conformal field theory has a structure 
of a homotopy Gerstenhaber algebra.
\end{theo}

In Subsection 4.4, the following result is proved:

\begin{theo}\label{theo2}
For a  locally-grading-restricted strong topological   vertex operator
 algebra $V$,
there is a locally convex completion $H$ of $V$ and 
an extension (still denoted $Q$) to $H$ of the operator $Q$ on $V$
such that 
$(H, Q)$ has a structure of  a
genus-zero topological conformal field theory. In particular, 
for a grading-restricted 
strong topological vertex operator algebra, the same conclusion
holds.
\end{theo}

This theorem proves Conjecture \ref{conj2} for 
 locally-grading-restricted strong topological   vertex operator
 algebras and in particular, for grading-restricted strong topological
vertex operator algebras. It is easy to see that for 
locally-grading-restricted topological  vertex operator
 algebras which are not strong, the same construction
does not give genus-zero topological conformal field theories. 

Combining Theorem \ref{KVZ} 
and Theorem \ref{theo2}, we obtain:

\begin{theo}\label{theo1}
For a  locally-grading-restricted strong topological vertex operator
 algebra $V$, there is a locally convex completion $H$ of $V$ such that 
 the Gerstenhaber algebra structure
on the cohomology of $V$ is extended to a homotopy Gerstenhaber algebra
structure on $H$. In particular, 
for a grading-restricted 
strong topological vertex operator algebra, the same conclusion
holds. \epf
\end{theo}

This theorem proves Conjecture \ref{conj1} for 
locally-grading-restricted strong topological   vertex operator
 algebras and in particular, for grading-restricted strong topological
vertex operator algebras.

\subsection{The main theorem}

Theorem \ref{theo2} and consequently Theorem 
\ref{theo1} are corollaries of the main theorem, Theorem \ref{main} below,
of the present paper.

In Subsection 4.2, for $k<0$,
a topological vertex partial operad ${\cal T}^{k}$
of type $k$
and a strong topological vertex partial operad 
$\bar{\cal T}^{k}$ of type $k$
are 
constructed. There is a natural 
morphism of partial operads from the operad $\wedge TK_{{\frak H}_{1}}$ 
to the partial operad ${\cal T}^{k}$ for any $k<0$. The notions of (weakly)
meromorphic algebra over ${\cal T}^{k}$ and
$\bar{\cal T}^{k}$ are introduced in the same subsection.
In Subsection 4.3, we prove
the following main theorem of the present paper:

\begin{theo}\label{main}
Let $k$ be an integer less than $0$. 
The category of  meromorphic ${\Bbb Z}\times {\Bbb Z}$-graded algebras
over ${\cal T}^{k}$ ($\bar{\cal T}^{k}$) 
is isomorphic to the category of grading-restricted (strong) topological
vertex operator algebras of type $k$. The category of  
weakly meromorphic ${\Bbb Z}\times {\Bbb Z}$-graded algebras
over ${\cal T}^{k}$ ($\bar{\cal T}^{k}$) 
is isomorphic to the category of  locally-grading-restricted
(strong) topological
vertex operator algebras of type $k$.
\end{theo}

\renewcommand{\theequation}{\thesection.\arabic{equation}}
\renewcommand{\therema}{\thesection.\arabic{rema}}
\setcounter{equation}{0}
\setcounter{rema}{0}

\section{Semi-infinite forms}

\subsection{Left invariant meromorphic 
vector fields satisfying the Virasoro relations}

For $c\in {\Bbb C}$, recall the $c/2$-th power $\tilde{K}^{c}(1)$
of the determinant line bundle over the moduli space $K(1)$
of spheres with tubes of type $(1, 1)$ constructed in \cite{H}.
In this subsection, we construct left invariant meromorphic 
vector fields on $\tilde{K}^{c}(1)$
 satisfying the Virasoro relations. These vector fields are needed 
in the construction of semi-infinite forms on the moduli spaces 
of spheres with tubes.

In \cite{H}, the holomorphic line bundle $\tilde{K}^{c}(1)$ over $K(1)$
was identified  with the space 
$$\{(A^{(0)}, (a_{0}, A^{(1)}); C)\;|\; A^{(0)}, A^{(1)}\in H\subset
{\Bbb C}^{\infty},
a_{0}\in {\Bbb C}^{\times}, C\in {\Bbb C}\},$$
where $H$ is the subset of ${\Bbb C}^{\infty}$ consisting of sequences
$A=\{A_{j}\}_{j>0}$ such that 
$\exp\left(\sum_{j>0}A_{j}x^{j+1}\frac{d}{dx}\right)x$ is convergent 
in a neighborhood  (depending on $A$) of $0$.  
There is an associative partial operation
$_{^{1}}\widetilde{\infty}^{c}_{^{0}}$ on $\tilde{K}^{c}(1)$ with identity
$\tilde I=({\bf 0}, (1, {\bf 0}); 1)$, here ${\bf 0}$ 
denote the infinite sequence with all 
components to be 0.
 Also
in \cite{H}, the notions of meromorphic functions on $\tilde{K}^{c}(1)$ and
the meromorphic tangent space of $\tilde{K}^{c}(1)$ at $\tilde{Q}
\in \tilde{K}^{c}(1)$ were introduced.

The following result is obvious:

\begin{propo}
The 
meromorphic tangent space $T_{\tilde{Q}}\tilde{K}^{c}(1)$ of 
$\tilde{K}^{c}(1)$ 
at any $\tilde{Q}\in \tilde{K}^{c}(1)$
is linearly isomorphic to the 
vector space of  infinite linear combinations of 
the meromorphic tangent vectors
$\frac{\p}{\p A_{j}^{(0)}}\mbar_{\tilde{Q}}$, 
$\frac{\p}{\p A_{j}^{(1)}}\mbar_{\tilde{Q}}$, 
$j>0$, $\frac{\p}{\p a_{0}}\mbar_{\tilde{Q}}$, and 
$\frac{\p}{\p C}\mbar_{\tilde{Q}}$. \epf
\end{propo}

A {\it meromorphic tangent field on $\tilde{K}^{c}(1)$} is a tangent 
vector fields on $\tilde{K}^{c}(1)$ mapping meromorphic functions to 
meromorphic functions. Then we have:

\begin{propo}
A tangent vector field on $\tilde{K}^{c}(1)$ 
is meromorphic if and only if it
is a an infinite linear combination of
the tangent fields
$\frac{\p}{\p A_{j}^{(0)}}$, 
$\frac{\p}{\p A_{j}^{(1)}}$, 
$j>0$, $\frac{\p}{\p a_{0}}$ and $\frac{\p}{\p C}$
with meromorphic functions on $\tilde{K}^{c}(1)$ as coefficients. \epf
\end{propo}

We denote 
the space of meromorphic tangent fields on $\tilde{K}^{c}(1)$ by
$\Gamma(T\tilde{K}^{c}(1))$. This space has a subspace consists of 
meromorphic tangent fields which are finite linear combinations of 
the tangent fields
$\frac{\p}{\p A_{j}^{(0)}}$, 
$\frac{\p}{\p A_{j}^{(1)}}$, 
$j>0$, $\frac{\p}{\p a_{0}}$ and $\frac{\p}{\p C}$
with meromorphic functions on $\tilde{K}^{c}(1)$ as coefficients. 
We denote this
subspace by  $\Gamma(\hat{T}\tilde{K}^{c}(1))$.

In \cite{H}, the first author found formal series $\Psi_{j}= \Psi_{j}(A,
B, a_{0})$, $j\in {\Bbb Z}$, in $A_{j}$, $B_{j}$ and $a_{0}$ such that
\begin{equation} e_{{\cal A}}(x) (\alpha_{0})^{\displaystyle
x\frac{d}{dx}}e^{-1}_{{\cal B}}(x^{-1})=e_{\Psi^{-}}(x^{-1})
e^{-1}_{\Psi^{+}}(x) (\alpha_{0})^{\displaystyle
x\frac{d}{dx}}e^{\displaystyle -\Psi_{0}x\frac{d}{dx}}, 
\end{equation}
where $\Psi^{+}=\{\Psi_{j}(A, B, a_{0})\}_{j>0}$,
$\Psi^{-}=\{\Psi_{-j}(A, B, a_{0})\}_{j>0}$,
$$e_{A}(x)=\exp\left(\sum_{j>0}A_{j}x^{j+1}\frac{d}{dx}\right)$$ 
for any
sequence $A$ and
$e_{A}^{-1}(x)$ is the compositional inverse of $e_{A}(x)$. 

For any $\tilde{P} \in \tilde{K}^{c}(1)$, we define the partial map 
$\ell_{\tilde{P}} : \tilde{K}^{c}(1) \to \tilde{K}^{c}(1)$ by 
$\ell_{\tilde{P}} 
(\tilde{Q})=\tilde{P} {}_{^{1}}\widetilde {\infty}^c_{^0} \tilde{Q}$ for
any $\tilde{Q} \in \tilde{K}^{c}(1)$ whenever the sewing makes sense. 
Recall from \cite{H} that the
tangent space $\hat{T}_{\tilde{I}} \tilde{K}^{c}(1)$ has a basis
$${\cal L}(j) =-\frac {\partial} {\partial A_{-j}^{(0)}}\lbar_{\tilde{I}}$$
for $j<0$,
$${\cal L}(j) =-\frac {\partial} {\partial A_{j}^{(1)}}\lbar_{\tilde{I}}$$
for $j>0$,
$${\cal L}(0) =-\frac {\partial}{\partial a_0}\lbar_{\tilde{I}},$$
$$ {\cal K}=C\frac {\p}{\p C}\lbar_{\tilde I}=\frac {\p}{\p C}
\lbar_{\tilde I}$$
and a bracket operation $[\cdot, \cdot]$ coming {from} the sewing operation
${}_1\widetilde \infty^c_0$ and satisfying the Virasoro relations 
with central charge $c$: 
$$
[{\cal L}(m), {\cal L}(n)]= (m-n) {\cal L}(m+n)+\frac{c}{12}(m^{3}-m)
\delta_{m+n, 0}{\cal K}.
$$

We define tangent vector fields ${\Bbb K}$ and ${\Bbb L}(j)$ for
$j \in {\Bbb Z}$ by 
$${\Bbb K}\mbar_{\tilde{P}}= (\ell_{\tilde{P}})_\ast ({\cal K})$$
and
$$
{\Bbb L}(j)\mbar_{\tilde{P}}= (\ell_{\tilde{P}})_\ast ({\cal L}(j))
$$
for $j \in {\Bbb Z}$ and $\tilde{P} \in \tilde{K}^{c}(1)$, where 
$(\ell_{\tilde{P}})_\ast$ is
the map from $T_{\tilde{I}} \tilde{K}^{c}(1)$ to $T_{\tilde{P}} 
\tilde{K}^{c}(1)$ 
induced from  
$\ell_{\tilde{P}}$. Observe that the identity element $\tilde I$ 
can be sewn to any
element $\tilde{P} \in \tilde{K}^{c}(1)$, therefore the 
the vector fields ${\Bbb K}$ and 
${\Bbb L}(j)$, $j \in {\Bbb Z}$, are well defined everywhere 
and are left invariant
with respect to the sewing operation.
Next we want to write down explicitly the vector fields defined
above. 
To do this we need more notations and lemmas.

Let 
$$\tilde{P}=( A^{(0)}, (a_0, A^{(1)}); C_{\tilde{P}})\in \tilde{K}^{c}(1)$$ 
and
$$\tilde{Q}=(B^{(0)}, (b_0, B^{(1)}); C_{\tilde{Q}})\in \tilde{K}^{c}(1)$$
such that 
$\tilde{P} {}_{^{1}}\widetilde \infty^c_{^{0}} \tilde{Q}$ exists.
In \cite{H}, it was shown that 
$$
\tilde{P} {}_{^{1}}\widetilde \infty^c_{^{0}} \tilde{Q}=(C^{(0)}, 
(c_0, C^{(1)}); 
C_{\tilde{P} {}_{^{1}}\widetilde \infty^c_{^{0}} \tilde{Q}} ) 
$$
where
$C^{(0)}=C^{(0)}(\tilde{P}, \tilde{Q})\in H$, $C^{(1)}=C^{(1)}(\tilde{P}, 
\tilde{Q})\in H$ and $c_0=
c_0(\tilde{P}, \tilde{Q})\in {\Bbb
C}^{\times}$ are determined by
\begin{eqnarray*}
\exp \left(\sum_{j>0} C^{(0)}_j x^{j+1} \frac d{dx}\right)x
&=&\exp \left(-\sum_{j>0} \Psi_{-j} x^{j+1} \frac d{dx}\right)\cdot\nno\\
&&\hspace{1em}\cdot\exp \left(\sum_{j>0} A^{(0)}_j x^{j+1} \frac d{dx}\right)x,\\
c_0^{(1)} \exp \left(\sum_{j>0} C^{(1)}_j x^{j+1} \frac d{dx}\right)x
&=&\exp \left(-\sum_{j>0} \Psi_j x^{j+1} \frac d{dx}\right)
(a_0 e^{-\Psi_0})^{\dps x\frac{d}{dx}}\cdot\nno\\
&&\hspace{1em}\cdot\exp \left(\sum_{j>0} ( B^{(1)}_j x^{j+1} 
\frac d{dx}\right)
(b_{0}^{(1)})^{\dps x\frac{d}{dx}}x
\end{eqnarray*}
and 
$$C_{\tilde{P} {}_{^{1}}\widetilde \infty^c_{^{0}} \tilde{Q}}=
e^{c\Gamma(A^{(1)}, B^{(0)}, a_{0}^{(1)})}C_{\tilde{P}}C_{\tilde{Q}},$$
where $\Gamma=\Gamma(A^{(1)}, B^{(0)}, a_{0}^{(1)})$ is a convergent 
series in $a_{0}^{(1)}$, $A_{j}^{(1)}$, $B^{(0)}_{j}$, $j>0$.
Recall the compositions $\circ$ for ${\Bbb C}^{\infty}$ and
for ${\Bbb C}^{\times}\times {\Bbb C}^{\infty}$
defined in \cite{H}, pages 46 and 47. Using these compositions, we have
\begin{eqnarray}
C^{(0)}(\tilde{P}, \tilde{Q}) 
&=&(-\Psi^{-})\circ A^{(0)},\label{c.1}\\
(c_{0}^{(1)}(\tilde{P}, \tilde{Q}), C^{(1)} (\tilde{P}, \tilde{Q})) &=&
(a_{0}^{(1)}e^{-\Psi_{0}}, -\Psi^{+}) \circ (b_{0}^{(1)}, 
B^{(1)}),\label{c.2}\\
c_{0}^{(1)}(\tilde{P}, \tilde{Q})
&=&a_{0}^{(1)}b_{0}^{(1)}e^{-\Psi_{0}},\label{c.3}\\
C^{(1)} (\tilde{P}, \tilde{Q})&=&(-\Psi^{+}) \circ 
B^{(1)}(a_{0}^{(1)}e^{-\Psi_{0}}).\label{c.4}
\end{eqnarray}
The following lemma is needed later:

\begin{lemma}\label{ideg}
Define the {\it weights} in the space of formal series in 
$a^{(1)}_0, b^{(1)}_0$, $A_j^{(0)}, A_j^{(1)}, B_j^{(0)}, B_j^{(1)}$, 
$j>0$, 
by setting
\begin{eqnarray*}
\wt A_j^{(0)}&=&-j,\\
\wt A_j^{(1)}&=&j,\\
\wt B_j^{(0)}&=&-j,\\
\wt B_j^{(1)}&=&j
\end{eqnarray*}
for $j \in \Bbb Z$, 
$$\wt a^{(1)}_0=0$$
and 
$$\wt b^{(1)}_0=0,$$
Then we have 
\begin{eqnarray*}
\wt \Gamma&=&0,\\
\wt \Psi_j&=&j,\\
\wt C^{(0)}_j&=&-j,\\ 
\wt C^{(1)}_j&=&j
\end{eqnarray*}
for $j>0$.
\end{lemma}
\pf 
By keeping track of the weights of all the series and 
operators in the relevant
proofs in \cite{H}, it is easy to see that the conclusion is true.
\epfv

\begin{propo}\label{tgt-fld}
For any $\tilde{P}=(A^{(0)}, (a_0, A^{(1)}); C_{\tilde{P}})
\in \tilde{K}^{c}(1)$, 
we have
\begin{eqnarray}
{\Bbb K}\mbar_{\tilde P}&=&C \frac {\p}{\p C}
\label{k}\lbar_{\tilde{P}},\\
{\Bbb L}(0)\mbar_{\tilde{P}} &=&-a_0 \frac{\p}{\p a_0}\lbar_{\tilde{P}}, 
\label{l0}\\
{\Bbb L}(j)\mbar_{\tilde{P}} &=&-\sum_{k>0} \frac{\p ((-\Psi^{+}) \circ 
B^{(1)}(a_{0}^{(1)}e^{-\Psi_{0}}))_{k}}{\p B^{(1)}_{j}}
\lbar_{B^{(0)}=B^{(1)}={\bf 0}}
 \frac {\partial} {\partial A^{(1)}_k}\lbar_{\tilde{P}} \nno\\
&=& -a_0^j  \frac {\partial} {\partial A^{(1)}_j}\lbar_{\tilde{P}}
\nno\\
&&-\sum_{k>j} \frac{\p ((-\Psi^{+}) \circ 
B^{(1)}(a_{0}^{(1)}e^{-\Psi_{0}}))_{k}}{\p B^{(1)}_{j}}
\lbar_{B^{(0)}=B^{(1)}={\bf 0}}
 \frac {\partial} {\partial A^{(1)}_k}\lbar_{\tilde{P}},\label{l+} \\
 {\Bbb L}(-j)\mbar_{\tilde{P}}&=&\sum_{k>0}\frac {\p \Psi_k}{\p B_j^{(0)}}
\lbar_{B^{(0)}={\bf 0}} \frac {\p}{\p A_k^{(1)}} \lbar_{\tilde{P}}+ 
\frac {\p \Psi_0}{\p B_j^{(0)}} \lbar_{ B^{(0)}={\bf 0}}
a_0 \frac {\p}{\p a_0}\lbar_{\tilde{P}}\nno\\
&& -\sum_{k>0} 
\frac {\partial C^{(0)}_{k}}
{\partial B_j^{(0)}}
\lbar_{B^{(0)}={\bf 0}} \frac {\partial} {\partial A_k^{(0)}}
\lbar_{\tilde{P}} + c\frac {\p \Gamma}{\p B_j^{(0)}}
\lbar_{B^{(0)}={\bf 0}}
 C \frac {\p}{\p C}\lbar_{\tilde{P}}\label{l-}
\end{eqnarray}
for $j>0$.
\end{propo}
\pf
By the chain rule,
\begin{eqnarray*}
{\Bbb L}(0)&=&
-\frac {\partial c_{0}^{(1)}}
{\partial b_0}\lbar_{B^{(0)}=B^{(1)}={\bf 0}, b_{0}=1} 
 \frac{\partial}{\partial a_0}\lbar_{\tilde{P}}\nno\\
&& - 
 \sum_{k>0} 
\frac {\partial C_k^{(0)}(\tilde{P}, \tilde{Q}) }
{\partial b_0}
\lbar_{B^{(0)}=B^{(1)}={\bf 0}, b_{0}=1} 
\frac{\partial}{\partial A_k^{(0)}}\lbar_{\tilde{P}} \\
&& -\sum_{k>0} \frac {\partial C_k^{(1)}(\tilde{P}, \tilde{Q}) }
{\partial b_0}
\lbar_{B^{(0)}=B^{(1)}={\bf 0}, b_{0}=1} 
\frac{\partial}{ \partial A_k^{(1)}}\lbar_{\tilde{P}}-
\frac {\p C_{\tilde P {}_1\widetilde \infty^c_0 \tilde Q}}{\p b_0}
\frac {\p}{\p C}\lbar_{\tilde P}
\end{eqnarray*}
\begin{eqnarray*}
{\Bbb L}(j)\mbar_{\tilde{P}}&=& 
-\frac {\partial c_{0}^{(1)}}
{\partial B^{(1)}_j}\lbar_{B^{(0)}=B^{(1)}={\bf 0}, b_{0}=1} 
\frac{ \partial}{\partial a_0}
\lbar_{\tilde{P}}\nno\\
&& -
 \sum_{k>0} 
\frac {\partial C_k^{(0)}(A, B) }
{\partial B_j^{(1)}}\lbar_{B^{(0)}=B^{(1)}={\bf 0}, b_{0}=1} 
\frac{\partial}{\partial A_k^{(0)}}\lbar_{\tilde{P}} \\
&& - \sum_{k>0} \frac {\partial C_k^{(1)}(A, B) }
{\partial B_j^{(1)}}\lbar_{B^{(0)}=B^{(1)}={\bf 0}, b_{0}=1} 
\frac{\partial}{ \partial A_k^{(1)}}\lbar_{\tilde{P}}-
\frac {\p C_{\tilde P {}_1\widetilde \infty^c_0 \tilde Q}}{\p B_j^{(1)}}
\frac {\p}{\p C}\lbar_{\tilde P}
\end{eqnarray*}
$j>0$, and
\begin{eqnarray*}
{\Bbb L}(-j)&=& 
-\frac {\partial c_{0}^{(1)}}
{\partial B^{(0)}_j}\lbar_{B^{(0)}=B^{(1)}={\bf 0}, b_{0}=1} 
\frac{ \partial }{\partial a_0}\lbar_{\tilde{P}}\nno\\
&& -
 \sum_{k>0}  
\frac {\partial C_k^{(0)}(A, B ) }
{\partial B_j^{(0)}}_{B^{(0)}=B^{(1)}={\bf 0}, b_{0}=1} 
\frac{ \partial}{ \partial A_k^{(0)}} \lbar_{\tilde{P}}\\
&& - \sum_{k>0}  \frac {\partial C_k^{(1)}(A, B) }
{\partial B_j^{(0)}}\lbar_{B^{(0)}=B^{(1)}={\bf 0}, b_{0}=1} 
\frac{\partial}{ \partial A_k^{(1)}}\lbar_{\tilde{P}}\nno\\
&&-\frac {\p C_{\tilde P {}_1\widetilde \infty^c_0 \tilde Q}}{\p B_j^{(0)}}
\lbar_{B^{(0)}=B^{(1)}={\bf 0}, b_{0}=1} 
\frac {\p}{\p C}\lbar_{\tilde P}
\end{eqnarray*}
$j>0$.
Using (\ref{c.1})--(\ref{c.4}), Lemma \ref{ideg} above
 and Proposition 2.2.5 in \cite{H}, 
we obtain the 
formulas (\ref{l0})--(\ref{l-}).
\epfv 

\begin{propo}
The tangent vector fields ${\Bbb L}(j)$, $j\in {\Bbb Z}$, are 
meromorphic and their restrictions ${\Bbb L}(j)\mbar_{\tilde{Q}}$, 
$j\in {\Bbb Z}$,
at any point
$\tilde Q\in T_{\tilde Q} {\tilde K}^c (1)$
together with ${\frac {\p}{\p C}}\mbar_{\tilde{Q}}$, form a basis at 
$\tilde Q$.
\end{propo}
\pf
The first conclusion follows {from} 
arguments using weights for all relevant
formal series. In fact,  all coefficients,
except 
$\frac {\p C^{(0)}_k}{\p B_j^{(0)} }$, in the explicit expressions
(\ref{l0})--(\ref{l-})
 of ${\Bbb L}(j)$, $j\in {\Bbb Z}$,
are weight-homogeneous 
polynomials in $a^{(1)}_0$, $A^{(1)}_{j}$, $j>0$. 
The coefficient $\frac {\p C^{(0)}_k}{\p B_j^{(0)} }$ is 
a weight-homogeneous 
polynomial in $a^{(1)}_0$, $A^{(1)}_{j}$, $A^{(0)}_{j}$, $j>0$.

To prove the second conclusion, note that 
we have the decomposition 
$\tilde Q=\tilde Q_1 {}_{^{1}}\widetilde \infty_{^{0}}^{c} \tilde Q_2$, 
where the local
coordinates of $\tilde Q_1$ and $\tilde Q_2$  
at $0$ and $\infty$, respectively, are trivial.  
In the case that $C_{\tilde{Q}}\ne 0$, the maps
 $$(\ell_{\tilde Q_2})_\ast : T_{\tilde I} {\tilde K}^c (1) \to
T_{\tilde Q_2} {\tilde K}^c (1)$$ 
and
$$(\ell_{\tilde Q_1})_\ast : T_{\tilde Q_2} {\tilde K}^c (1) \to
T_{\tilde Q} {\tilde K}^c (1)$$ are linear isomorphisms,
since $\tilde Q_1$ and
$\tilde Q_2$ are invertible in $\tilde K^c(1)$. So
$(\ell_{\tilde Q})_\ast=(\ell_{\tilde Q_1})_\ast 
\circ (\ell_{\tilde Q_2})_\ast$
is also a linear isomorphism, proving the second conclusion in this case. 
In the case of 
$C_{\tilde Q}=0$, since $\tilde{K}^{c}(1)=K(1)\times {\Bbb C}$, 
$\tilde Q \in K(1)\times \{0\}$. Let $\tilde{Q}=(Q; 0)$ where $Q\in K(1)$. 
Then in the decomposition 
$\tilde Q=\tilde Q_1 {}_{^{1}}\widetilde \infty_{^{0}}^{c} \tilde Q_2$,
we can also choose $\tilde Q_1$ and $\tilde Q_2$ to be in $K(1)\times \{0\}$
such that $\tilde Q_1=(Q_{1}; 0)$ and $\tilde{Q}_{2}=(Q_{2}, 0)$.
Since $Q_{1}$ and $Q_{2}$ are invertible, 
$(\ell_{\tilde Q_2})_\ast$ and $(\ell_{\tilde Q_1})_\ast$ induce linear 
isomorphisms {from} $T_{\tilde I} {\tilde K}^c (1)/{\Bbb C}
{\frac {\p}{\p C}}\mbar_{\tilde{I}}$ to
$T_{\tilde{Q}_{2}} {\tilde K}^c (1)/{\Bbb C}
{\frac {\p}{\p C}}\mbar_{\tilde{Q}_{2}}$ and {from} 
$T_{\tilde{Q}_{2}} {\tilde K}^c (1)/{\Bbb C}
{\frac {\p}{\p C}}\mbar_{\tilde{Q}_{2}}$ to
$T_{\tilde{Q}} {\tilde K}^c (1)/{\Bbb C}
{\frac {\p}{\p C}}\mbar_{\tilde{Q}}$, respectively. Thus the second 
conclusion is true in this case.
\epfv

Let ${\frak L}$ be the Virasoro algebra and 
let $L_{j}$, $j\in {\Bbb Z}$, and $d$, be the usual basis of ${\frak L}$. 
Recall that in \cite{H}, the space of meromorphic functions on 
$\tilde{K}^{c}(1)$
is denoted by $\tilde{D}^{c}_{1}$. The eigenspaces of 
${\Bbb L}(0)$ as an operator
on $\tilde{D}^{c}_{1}$ give  a ${\Bbb Z}$-grading to 
$\tilde{D}^{c}_{1}$. Let $\tilde{D}^{c; 1}_{1}$ be the space of 
meromorphic functions on $\tilde{K}^{c}(1)$ consisting of elements 
which are linear in the last coordinate $C$ of 
$\tilde{Q}=(A^{(0)}, (a_{0}^{(1)}, A^{(1)}); C)\in \tilde{K}^{c}(1)$.

\begin{propo}\label{vir-rel}
The restriction ${\Bbb L}(j)\mbar_{\tilde{D}^{c; 1}_{1}}$, 
$j\in {\Bbb Z}$, and
${\Bbb K}\mbar_{\tilde{D}^{c; 1}_{1}}$,  of the 
meromorphic vector fields ${\Bbb L}(j)$, $j\in {\Bbb Z}$, and
${\Bbb K}$, to $\tilde{D}^{c; 1}_{1}$ satisfy the 
Virasoro relations with central charge $c$:
\begin{eqnarray}
[{\Bbb L}(m)\mbar_{\tilde{D}^{c; 1}_{1}}, 
{\Bbb L}(n)\mbar_{\tilde{D}^{c; 1}_{1}}]
&=&(m-n){\Bbb L}(m+n)\mbar_{\tilde{D}^{c; 1}_{1}}+\frac{c}{12}(m^{3}-m)
\delta_{m+n, 0} {\Bbb K}\mbar_{\tilde{D}^{c; 1}_{1}}, \nno\\
&&\label{vir-rel1}\\
{[{\Bbb L}(m)\mbar_{\tilde{D}^{c; 1}_{1}}, {\Bbb
K}\mbar_{\tilde{D}^{c; 1}_{1}}]}
&=&0. \label{vir-rel2}
\end{eqnarray}
In particular, 
the ${\Bbb Z}$-graded vector space $\tilde{D}^{c;1}_{1}$ equipped with
the map $\pi: {\frak L}\to \mbox{\rm End}\; \tilde{D}^{c;1}_{1}$
defined by 
$\pi(L_{j})={\Bbb L}(j)\mbar_{\tilde{D}^{c; 1}_{1}}$, $j\in {\Bbb Z}$ and 
$\pi(d)=c{\Bbb K}\mbar_{\tilde{D}^{c; 1}_{1}}=cI_{\tilde{D}^{c;
1}_{1}}$ ($I_{\tilde{D}^{c;
1}_{1}}$ being the identity map on $\tilde{D}^{c;
1}_{1}$), is a module for the 
Virasoro algebra
with central charge $c$.
\end{propo}
\pf
The second equation (\ref{vir-rel2}) is obvious {from} the Proposition 
\ref{tgt-fld}. 
We prove the first equation (\ref{vir-rel1}). 
Fix $\tilde Q_0 \in \tilde K^c (1)$.
Let $F$ be any meromorphic function on $\tilde K^c (1)$. Then for any
$m>0$, $n<0$ and 
$\tilde Q \in \tilde K^c (1)$, 
by definition, we have 
\begin{eqnarray*}
(\Bbb L(i) F)(\tilde Q) &=& 
-\frac {\p}{\p A^{(1)}_i}\lbar_{\tilde P_1=\tilde I}F(\tilde Q {}_1 
{\widetilde \infty}^c_0 \tilde P_1 ), \\
(\Bbb L(j) F)(\tilde Q) &=& -\frac {\p}{\p B^{(0)}_j}
\lbar _{\tilde P_2=\tilde I}
 F(\tilde Q {}_1 {\widetilde \infty}^c_0 \tilde P_2 ), 
\end{eqnarray*}
where 
$\tilde P_1=( A^{(0)}, (\alpha_0, A^{(1)}), C_1)$ 
$\tilde P_2=( B^{(0)}, (\beta_0, B^{(1)}), C_2)$ 
are elements near $\tilde I$. 

We have
\begin{eqnarray*}
(\Bbb L(m) \Bbb L(n) )\mbar_{\tilde Q_0} F &=&
\Bbb L(m)\mbar_{\tilde Q_0} (\Bbb L(n)  F) \\
&=&-\frac {\p}{\p A^{(1)}_m}\lbar_{\tilde P_1=\tilde I}
(\Bbb L(n)  F)(\tilde Q_0 {}_{^{1}} {\widetilde \infty} ^c_{^{0}} 
\tilde P_1 )\\
&=&\frac {\p}{\p A^{(1)}_m}\lbar_{\tilde P_1=\tilde I}
\frac {\p}{\p B^{(0)}_n}\lbar _{\tilde P_2=\tilde I}
 F(\tilde Q_0 {}_{^{1}} {\widetilde \infty}^c_{^{0}} \tilde P_1 
{}_{^{1}} {\widetilde \infty}^c_{^{0}} \tilde P_2  ) \\
&=& \frac {\p}{\p A^{(1)}_m}\lbar_{\tilde P_1=\tilde I}
\frac {\p}{\p B^{(0)}_n}\lbar _{\tilde P_2=\tilde I}
 ((\ell_{Q_0})^* F)( \tilde P_1 
{}_{^{1}} {\widetilde \infty}^c_{^{0}} \tilde P_2 ).
\end{eqnarray*}
Similarly, 
\begin{eqnarray*}
(\Bbb L(n) \Bbb L(m) )\mbar_{\tilde Q_0} F &=&
\Bbb L(n)\mbar_{\tilde Q_0} (\Bbb L(m)  F) \\
&=&-\frac {\p}{\p A^{(0)}_n}\lbar_{\tilde P_1=\tilde I}
(\Bbb L(m)  F)(\tilde Q_0 {}_{^{1}} {\widetilde \infty} ^c_{^{0}} 
\tilde P_1 )\\
&=&\frac {\p}{\p A^{(0)}_n}\lbar_{\tilde P_1=\tilde I}
\frac {\p}{\p B^{(1)}_m}\lbar _{\tilde P_2=\tilde I}
 F(\tilde Q_0 {}_{^{1}} {\widetilde \infty}^c_{^{0}} \tilde P_1 
{}_{^{1}} {\widetilde \infty}^c_{^{0}} \tilde P_2  ) \\
&=&\frac {\p}{\p A^{(0)}_n}\lbar_{\tilde P_1=\tilde I}
\frac {\p}{\p B^{(1)}_m}\lbar _{\tilde P_2=\tilde I} 
((\ell_{\tilde Q_0 })^* F)
 ( \tilde P_1 
{}_{^{1}} {\widetilde \infty}^c_{^{0}} \tilde P_2  ). 
\end{eqnarray*}
Thus, by Proposition 3.5.2 in \cite{H},
\begin{eqnarray*}
\lefteqn{[\Bbb L(m), \Bbb L(n)]\mbar_{\tilde Q_0} F }\nno\\
&&=
\left(\frac {\p}{\p A^{(1)}_m}\lbar_{\tilde P_1=\tilde I}
\frac {\p}{\p B^{(0)}_n}\lbar _{\tilde P_2=\tilde I} -
\frac {\p}{\p A^{(0)}_n}\lbar_{\tilde P_1=\tilde I}
\frac {\p}{\p B^{(1)}_m}\lbar _{\tilde P_2=\tilde I}\right)\cdot\nno\\ 
&&\quad\quad\quad\quad\cdot ((\ell_{\tilde Q_0 })^* F)
 ( \tilde P_1 
{}_1 {\widetilde \infty}^c_0 \tilde P_2  ) \\
&&= ( (m-n) {\cal L} (m+n) +\frac {c}{12} (m^3-m)
 \delta_{m+n, 0} C\frac {\p}{\p C} )
 (\ell_{\tilde Q_0 })^* (F) \\
&&=(\ell_{\tilde Q_0 })_\ast
 ( (m-n) {\cal L} (m+n) +\frac {c}{12} (m^3-m) \delta_{m+n, 0} 
C\frac {\p}{\p C}) F\\
&&= ( (m-n) {\Bbb L} (m+n)+\frac {c}{12} (m^3-m) \delta_{m+n, 0} {\Bbb K} )
\mbar_{\tilde Q_0} F,
\end{eqnarray*}
proving the Virasoro relation (\ref{vir-rel1}) for $m>0, n<0$. 
We can prove (\ref{vir-rel1}) for other $m$ and $n$ similarly.
\epfv

\subsection{Semi-infinite forms on $K(0)$}

In this subsection we introduce and study semi-infinite forms on $K(0)$.
We shall use structures on  $K(1)$ to define 
semi-infinite forms on $K(0)$. 

We need the embedding map ${\frak E}$ {from} $K(0)$ to
$K(1)$ defined by 
$${\frak E}(A)= (A, (1, {\bf 0})).$$ 
It has a left inverse 
${\frak
D}$ defined by 
$${\frak
D}(A^{(0)}, (a_{0}^{(1)}, A^{(1)}))=
(A^{(0)}, (a_{0}^{(1)}, A^{(1)}))_{^{1}}\infty_{^{0}}{\bf 0}.$$ 
These two maps induce $\tilde{\frak E}: \tilde{K}^{c}(0)\to
\tilde{K}^{c}(1)$ and $\tilde{\frak D}: \tilde{K}^{c}(1)\to
\tilde{K}^{c}(0)$. 

For $\tilde{Q}\in \tilde{K}^{c}(1)$, we denote the subspace of 
$T_{\tilde{Q}}\tilde{K}^{c}(1)$ 
 consisting of finite linear combinations of 
${\Bbb L}(j)\mbar_{\tilde{Q}}$,
$j\in {\Bbb Z}$, $\frac{\p}{\p C}\mbar_{\tilde{P}}$, 
by $\hat{T}_{\tilde{Q}}\tilde{K}^{c}(1)$. 
This space has a ${\Bbb Z}$-grading 
called {\it weight} defined by $\wt {\Bbb L}(j)\mbar_{\tilde{P}}=-j$, 
$j\in {\Bbb Z}$,
$\wt {\Bbb K}\mbar_{\tilde{P}}=0$. 
For $\tilde{P}\in \tilde{K}^{c}(0)$, we denote the subspace 
of $T_{\tilde{P}}\tilde{K}^{c}(0)$ 
consisting of finite linear combinations of 
$\tilde{\frak D}_{*}({\Bbb L}(j))\mbar_{\tilde{P}}$,
$j<-1$, $\tilde{\frak D}_{*}(\frac{\p}{\p C})\mbar_{\tilde{P}}$,
by $\hat{T}_{\tilde{P}}\tilde{K}^{c}(0)$. Clearly, 
$$\hat{T}_{\tilde{P}}\tilde{K}^{c}(0)=(\tilde{\frak
D}_{*})\mbar_{\tilde{\frak E}(\tilde{P})}(\hat{T}_{\tilde{\frak E}(\tilde{P})}
\tilde{K}^{c}(1)).$$
The union of $\hat{T}_{\tilde{P}}\tilde{K}^{c}(0)$ for $\tilde{P}\in
\tilde{K}^{c}(0)$ is a holomorphic vector bundle 
$\hat{T}\tilde{K}^{c}(0)$ over $\tilde{K}^{c}(0)$. Using 
the flat section $\psi_{0}$ of $\tilde{K}^{c}(0)$ constructed in 
\cite{H}, we pull 
$\hat{T}\tilde{K}^{c}(0)$  back to 
 a holomorphic vector bundle $\psi_{0}^{*}(\hat{T}\tilde{K}^{c}(0))$ 
over $K(0)$.

For $\tilde{Q}\in \tilde{K}^{c}(1)$, we denote the graded dual space of
$\hat{T}_{\tilde{Q}}\tilde{K}^{c}(1)$ by
$\hat{T}'_{\tilde{Q}}\tilde{K}^{c}(1)$. The union of 
$\hat{T}'_{\tilde{Q}}\tilde{K}^{c}(0)$ for $\tilde{Q}\in
\tilde{K}^{c}(1)$ is a holomorphic vector bundle 
$\hat{T}'\tilde{K}^{c}(1)$ over $\tilde{K}^{c}(1)$.
Let 
${\Bbb L}'(j)$,  
$j\in {\Bbb Z}$, and $(\frac{\p}{\p C})'$ 
be the sections of $\hat{T}'\tilde{K}^{c}(1)$ such that
for $\tilde{Q}\in \tilde{K}^{c}(1)$,
${\Bbb L}'(j)\mbar_{\tilde{Q}}$, $j\in {\Bbb Z}$, and 
$(\frac{\p}{\p C})'\mbar_{\tilde{Q}}$ form the dual basis
of the basis ${\Bbb L}(j)\mbar_{\tilde{Q}}$, $j\in {\Bbb Z}$, and 
$\frac{\p}{\p C}\mbar_{\tilde{Q}}$. 
For $\tilde{P}\in
\tilde{K}^{c}(0)$, we  denote the subspace of the dual space
of $T_{\tilde{P}}\tilde{K}^{c}(0)$
consisting of finite linear combinations of $\tilde{\frak E}^{*}({\Bbb
L}'(j))\mbar_{\tilde{P}}$, $j<-1$, and 
$\tilde{\frak E}^{*}((\frac{\p}{\p C})')\mbar_{\tilde{P}}$ by 
$\hat{T}'_{\tilde{P}}\tilde{K}^{c}(0)$. Clearly, 
$$\hat{T}'_{\tilde{P}}\tilde{K}^{c}(0)=\tilde{\frak
E}^{*}\mbar_{\tilde{\frak E}(\tilde{P})}(\hat{T}'_{\tilde{\frak E}(\tilde{P})}
\tilde{K}^{c}(1)).$$ 
Note that the kernel of $\tilde{\frak
E}^{*}\mbar_{\tilde{\frak E}(\tilde{P})}$ 
is spanned by ${\Bbb L}'(j)$, 
$j\ge -1$. 
The union of kernels of $\tilde{\frak
E}^{*}\mbar_{\tilde{\frak E}(\tilde{P})}$  for $\tilde{P}\in
\tilde{K}^{c}(0)$ forms a holomorphic vector bundle 
$\kr \tilde{\frak
E}^{*}$ over $\tilde{K}^{c}(0)$. Using 
the flat section $\psi_{0}$ of $\tilde{K}^{c}(0)$, we pull 
$\kr \tilde{\frak
E}^{*}$  back to 
 a holomorphic vector bundle $\psi_{0}^{*}(\kr \tilde{\frak
E}^{*})$ 
over $K(0)$. 

Motivated by the semi-infinite cohomology of 
graded Lie algebras introduced by Feigin in \cite{F} and developed 
by Frenkel-Garland-Zuckerman in \cite{FGZ} and by Lian-Zuckerman 
\cite{LZ1}, we now consider 
the holomorphic vector bundle
$(\wedge \psi_{0}^{*}(\kr \tilde{\frak
E}^{*}))\wedge (\wedge \psi_{0}^{*}(\hat{T}\tilde{K}^{c}(0)))$
where $\wedge \psi_{0}^{*}(\kr \tilde{\frak
E}^{*})$ and $\wedge \psi_{0}^{*}(\hat{T}\tilde{K}^{c}(0))$ 
are the wedge product
bundles of
$\psi_{0}^{*}(\kr \tilde{\frak
E}^{*})$ and $\psi_{0}^{*}(\hat{T}\tilde{K}^{c}(0))$, respectively. 
For any $P\in K(0)$, the fiber of 
$(\wedge \psi_{0}^{*}(\kr \tilde{\frak
E}^{*}))\wedge (\wedge \psi_{0}^{*}(\hat{T}\tilde{K}^{c}(0)))$ 
over $P$
is a
${\Bbb Z}\times 
{\Bbb Z}$-graded vector spaces spanned by the canonical basis
$${\Bbb L}'(i_{1})\mbar_{\tilde{\frak E}(\psi_{0}(P))}\wedge \dots \wedge
{\Bbb L}'(i_{m})\mbar_{\tilde{\frak E}(\psi_{0}(P))}\wedge \tilde{\frak
D}_{*}({\Bbb L}(j_{1}))\mbar_{\psi_{0}(P)}
\wedge \dots \wedge \tilde{\frak
D}_{*}({\Bbb L}(j_{n}))\mbar_{\psi_{0}(P)},$$
$m, n>0$, $i_{1}, \dots, i_{m}\ge -1$, $j_{1}, \dots, j_{n}<-1$. 
The integer $m+n$ is called the {\it fermion number} or 
{\it ghost number} of these
basis elements. If 
$m+n\in 2{\Bbb Z}$ ($\in 2{\Bbb Z}+1$), 
these basis
elements are said to be {\it even} ({\it odd}). 
The integers $-\sum_{k=1}^{m}i_{k}+\sum_{k=1}^{n}j_{k}$
are called 
{\it weights}  or {\it degrees} of the canonical basis 
elements.

This bundle has holomorphic sections 
$$P\mapsto 
{\Bbb L}'(i_{1})\mbar_{\tilde{\frak E}(\psi_{0}(P))}\wedge \dots \wedge
{\Bbb L}'(i_{m})\mbar_{\tilde{\frak E}(\psi_{0}(P))}\wedge 
\tilde{\frak D}_{*}({\Bbb L}(j_{1}))\mbar_{\psi_{0}(P)}
\wedge \dots \wedge \tilde{\frak D}_{*}({\Bbb L}(j_{n}))
\mbar_{\psi_{0}(P)},$$
$m, n\ge 0$, $i_{1}, \dots, i_{m}\ge -1$, $j_{1}, \dots, j_{n}<-1$.
We shall write these sections simply as 
$$\psi_{0}^{*}({\Bbb L}'(i_{1}))\wedge \dots \wedge
\psi_{0}^{*}({\Bbb L}'(i_{m}))\wedge 
\psi_{0}^{*}(\tilde{\frak D}_{*}({\Bbb L}(j_{1})))
\wedge \dots \wedge \psi_{0}^{*}(\tilde{\frak D}_{*}({\Bbb L}(j_{n}))).$$
We denote the space of these sections by $\wedge_{\infty}
\psi_{0}^{*}(\tilde{\frak D}_{*}({\Bbb L}))$.

We also consider 
the graded dual bundle
$$((\wedge \psi_{0}^{*}(\kr \tilde{\frak
E}^{*}))\wedge (\wedge \psi_{0}^{*}(\hat{T}\tilde{K}^{c}(0))))'
=(\wedge \psi_{0}^{*}(\kr \tilde{\frak
E}^{*}))'\wedge (\wedge \psi_{0}^{*}(\hat{T}\tilde{K}^{c}(0)))'$$
of $(\wedge \psi_{0}^{*}(\kr \tilde{\frak
E}^{*}))\wedge (\wedge \psi_{0}^{*}(\hat{T}\tilde{K}^{c}(0)))$.
For any $P\in K(0)$, the fiber of 
$$(\wedge \psi_{0}^{*}(\kr \tilde{\frak
E}^{*}))'\wedge (\wedge \psi_{0}^{*}(\hat{T}\tilde{K}^{c}(0)))'$$
is a
${\Bbb Z}\times 
{\Bbb Z}$-graded vector spaces spanned by the canonical basis
$${\Bbb L}(i_{1})\mbar_{\tilde{\frak E}(\psi_{0}(P))}\wedge \dots \wedge
{\Bbb L}(i_{m})\mbar_{\tilde{\frak E}(\psi_{0}(P))}\wedge \tilde{\frak
E}^{*}({\Bbb L}'(j_{1}))\mbar_{\psi_{0}^{*}(P)}
\wedge \dots \wedge \tilde{\frak
E}^{*}({\Bbb L}'(j_{n}))\mbar_{\psi_{0}^{*}(P)},$$
$m, n\ge 0$, $m+n>0$,  $i_{1}, \dots, i_{m}\ge -1$, $j_{1}, \dots, j_{n}<-1$. 
The integer $m+n$ is called the {\it fermion number} or 
{\it ghost number} of these
basis elements. If 
$m+n\in 2{\Bbb Z}$($\in 2{\Bbb Z}+1$), 
these basis
elements are said to be {\it even} ({\it odd}). 
The integers $-\sum_{k=1}^{m}i_{k}+\sum_{k=1}^{n}j_{k}$
are called 
{\it weights}  or {\it degrees} of the canonical basis 
elements.

This bundle has holomorphic sections 
$$P\mapsto {\Bbb L}(i_{1})\mbar_{\tilde{\frak E}(\psi_{0}(P))}
\wedge \dots \wedge
{\Bbb L}(i_{m})\mbar_{\tilde{\frak E}(\psi_{0}(P))}\wedge \tilde{\frak
E}^{*}({\Bbb L}'(j_{1}))\mbar_{\psi_{0}(P)}
\wedge \dots \wedge \tilde{\frak
E}^{*}({\Bbb L}'(j_{n}))\mbar_{\psi_{0}(P)},$$
$m, n>0$, $i_{1}, \dots, i_{m}\ge -1$, $j_{1}, \dots, j_{n}<-1$.
We shall write these sections simply as 
$$\psi_{0}^{*}({\Bbb L}(i_{1}))\wedge \dots \wedge
\psi_{0}^{*}({\Bbb L}(i_{m}))\wedge \psi_{0}^{*}(\tilde{\frak
E}^{*}({\Bbb L}'(j_{1})))
\wedge \dots \wedge \psi_{0}^{*}(\tilde{\frak
E}^{*}({\Bbb L}'(j_{n}))).$$
We denote the space of these sections by 
$\wedge_{\infty}\psi_{0}^{*}(\tilde{\frak
E}^{*}({\Bbb L}'))$. It is clear that $\wedge_{\infty}\psi_{0}^{*}(\tilde{\frak
E}^{*}({\Bbb L}'))$ is naturally isomorphic to the graded dual 
of $\wedge_{\infty}
\psi_{0}^{*}(\tilde{\frak D}_{*}({\Bbb L}))$. We shall 
identify $\wedge_{\infty}\psi_{0}^{*}(\tilde{\frak
E}^{*}({\Bbb L}'))$ with the graded dual of $\wedge_{\infty}
\psi_{0}^{*}(\tilde{\frak D}_{*}({\Bbb L}))$.

Note that both $\wedge_{\infty}\psi_{0}^{*}(\tilde{\frak D}_{*}({\Bbb L}))$
and $\wedge_{\infty}\psi_{0}^{*}(\tilde{\frak
E}^{*}({\Bbb L}'))$
are linearly isomorphic to 
the space $\wedge_{\infty}{\frak W}$ of semi-infinite forms
on the Witt algebra ${\frak W}$.

A holomorphic section of 
$$(\wedge \psi_{0}^{*}(\kr \tilde{\frak
E}^{*}))\wedge (\wedge \psi_{0}^{*}(\hat{T}\tilde{K}^{c}(0)))$$
or $$(\wedge \psi_{0}^{*}(\kr \tilde{\frak
E}^{*}))'\wedge (\wedge \psi_{0}^{*}(\hat{T}\tilde{K}^{c}(0)))'$$
is said to be {\it meromorphic} if it is a linear combination of 
sections of the form
$$P\mapsto f(P)\psi_{0}^{*}({\Bbb L}'(i_{1}))\wedge \dots \wedge
\psi_{0}^{*}({\Bbb L}'(i_{m}))\wedge 
\psi_{0}^{*}(\tilde{\frak D}_{*}({\Bbb L}(j_{1})))
\wedge \dots \wedge \psi_{0}^{*}(\tilde{\frak D}_{*}({\Bbb L}(j_{n})))$$
or 
$$P\mapsto f(P)\psi_{0}^{*}({\Bbb L}(i_{1}))\wedge \dots \wedge
\psi_{0}^{*}({\Bbb L}(i_{m}))\wedge 
\psi_{0}^{*}(\tilde{\frak E}^{*}({\Bbb L}'(j_{1})))
\wedge \dots \wedge \psi_{0}^{*}(\tilde{\frak E}^{*}({\Bbb L}'(j_{n}))),$$
respectively,
where $f$ is a meromorphic function on $K(0)$
and $m, n>0$, $i_{1}, \dots, i_{m}\ge -1$, $j_{1}, \dots, j_{n}<-1$.
In other words, the spaces  of meromorphic sections of 
$$(\wedge \psi_{0}^{*}(\kr \tilde{\frak
E}^{*}))\wedge (\wedge \psi_{0}^{*}(\hat{T}\tilde{K}^{c}(0)))$$
and 
$$(\wedge \psi_{0}^{*}(\kr \tilde{\frak
E}^{*}))'\wedge (\wedge \psi_{0}^{*}(\hat{T}\tilde{K}^{c}(0)))'$$
are linearly isomorphic to 
$D_{0}\otimes \wedge_{\infty}\psi_{0}^{*}(\tilde{\frak D}_{*}({\Bbb L}))$
and 
$D_{0}\otimes \wedge_{\infty}\psi_{0}^{*}(\tilde{\frak
E}^{*}({\Bbb L}'))$, respectively.
 Note that 
$D_{0}\otimes \wedge_{\infty}\psi_{0}^{*}(\tilde{\frak D}_{*}({\Bbb L}))$ is 
the semi-infinite analogue of the space of skew-symmetric holomorphic 
poly-vector fields
on a complex manifold while $D_{0} \otimes 
\wedge_{\infty}\psi_{0}^{*}(\tilde{\frak
E}^{*}({\Bbb L}'))$ is the semi-infinite analogue of 
the space of holomorphic forms on a complex manifold. 

The semi-infinite analogue $D_{0} \otimes 
\wedge_{\infty}\psi_{0}^{*}(\tilde{\frak
E}^{*}({\Bbb L}'))$ of 
the space of holomorphic forms on a complex 
manifold can also be constructed 
as the semi-infinite analogue on another submanifold of $K(1)$
of the space of skew-symmetric holomorphic poly-vector fields
on a complex manifold. Since in the next subsection we shall need 
not only $\wedge_{\infty}\psi_{0}^{*}(\tilde{\frak D}_{*}({\Bbb L}))$
but also $\wedge_{\infty}\psi_{0}^{*}(\tilde{\frak
E}^{*}({\Bbb L}'))$ and since
this construction gives a conceptual explanation of why we need 
$\wedge_{\infty}\psi_{0}^{*}(\tilde{\frak
E}^{*}({\Bbb L}'))$, we give this construction in detail here.

Let $K_{\ge -1}(1)$ be the submanifold of 
$K(1)$ consisting of elements of the form
$(A(1; z), (a_{0}^{(1)}, A^{(1)}))$, 
$z\in {\Bbb C}$, $a_{0}\in {\Bbb C}^{\times}$
and $A^{(1)}\in H$ (recalling {from} \cite{H}
that $A(1; z)$ is the sequence with the 
first component $z$ and the other components $0$). Then it is clear that
$K_{\ge -1}(1)$ is closed under the sewing operation. 

We need the embedding map ${\frak I}$ {from} $K_{\ge -1}(1)$ to
$K(1)$. It has a left inverse ${\frak
N}$ defined by 
$${\frak
N}(A^{(0)}, (a_{0}^{(1)}, A^{(1)}))=
(A(1; A_{1}^{(0)}), (a_{0}^{(1)}, A^{(1)})).$$
Let $\tilde{K}^{c}_{\ge -1}(1)$ be the restriction of 
the line bundle $\tilde{K}^{c}(1)$ to $K_{\ge -1}(1)$. 
Then the two maps above induce $\tilde{\frak I}: \tilde{K}^{c}_{\ge -1}(1)\to
\tilde{K}^{c}(1)$ and $\tilde{\frak N}: \tilde{K}^{c}(1)\to
\tilde{K}^{c}_{\ge -1}(1)$.

For $\tilde{P}\in \tilde{K}^{c}_{\ge -1}(1)$, we 
denote
the subspace 
of $T_{\tilde{P}}\tilde{K}^{c}_{\ge -1}(1)$ 
consisting of finite linear combinations of 
$\tilde{\frak N}_{*}({\Bbb L}(j))\mbar_{\tilde{P}}$,
$j\ge -1$, $\tilde{\frak N}_{*}(\frac{\p}{\p C})\mbar_{\tilde{P}}$,
by $\hat{T}_{\tilde{P}}\tilde{K}^{c}_{\ge -1}(1)$. Then
$$\hat{T}_{\tilde{P}}\tilde{K}^{c}_{\ge -1}(1)=(\tilde{\frak
N}_{*})\mbar_{\tilde{\frak I}(\tilde{P})}(\hat{T}_{\tilde{\frak I}(\tilde{P})}
\tilde{K}^{c}(1)).$$
The union of $\hat{T}_{\tilde{P}}\tilde{K}^{c}_{\ge -1}(1)$ for $\tilde{P}\in
\tilde{K}^{c}_{\ge -1}(1)$ is a holomorphic vector bundle 
$\hat{T}\tilde{K}^{c}_{\ge -1}(1)$ over $\tilde{K}^{c}_{\ge -1}(1)$. 
For simplicity, we also use the same notation 
$\psi_{1}$ to denote the
flat section $\psi_{1}$ of $\tilde{K}^{c}(1)$ constructed in \cite{H} and
its restriction to $K_{\ge -1}(1)$ of . Using 
$\psi_{1}$, we pull 
$\hat{T}\tilde{K}^{c}_{\ge -1}(1)$  back to 
 a holomorphic vector bundle $\psi_{1}^{*}(\hat{T}\tilde{K}^{c}_{\ge -1}(1))$ 
over $K_{\ge -1}(1)$.

For $\tilde{P}\in
\tilde{K}^{c}_{\ge -1}(1)$, we  denote the subspace
of the dual space of $T_{\tilde{P}}\tilde{K}^{c}_{\ge -1}(1)$
consisting of finite linear combinations of  $\tilde{\frak I}^{*}({\Bbb
L}'(j))\mbar_{\tilde{P}}$, $j\ge -1$, and 
$\tilde{\frak I}^{*}((\frac{\p}{\p C})')\mbar_{\tilde{P}}$ by 
$\hat{T}'_{\tilde{P}}\tilde{K}^{c}_{\ge -1}(1)$. Then
$$\hat{T}'_{\tilde{P}}\tilde{K}^{c}_{\ge -1}(1)=\tilde{\frak
I}^{*}\mbar_{\tilde{\frak I}(\tilde{P})}(\hat{T}'_{\tilde{\frak I}(\tilde{P})}
\tilde{K}^{c}(1)).$$ 
Note that the kernel of $\tilde{\frak
I}^{*}\mbar_{\tilde{\frak I}(\tilde{P})}$ 
is spanned by ${\Bbb L}'(j)$, 
$j< -1$. 
The union of kernels of $\tilde{\frak
I}^{*}\mbar_{\tilde{\frak I}(\tilde{P})}$  for $\tilde{P}\in
\tilde{K}^{c}_{\ge -1}(1)$ forms a holomorphic vector bundle 
$\kr \tilde{\frak
I}^{*}$ over $\tilde{K}^{c}_{\ge -1}(1)$. Using 
the flat section $\psi_{1}$ of $\tilde{K}^{c}_{\ge -1}(1)$, we pull 
$\kr \tilde{\frak
I}^{*}$  back to 
 a holomorphic vector bundle $\psi_{1}^{*}(\kr \tilde{\frak
I}^{*})$ 
over $K_{\ge -1}(1)$.

We consider 
the holomorphic vector bundle
$$(\wedge \psi_{1}^{*}(\kr \tilde{\frak
I}^{*}))\wedge (\wedge \psi_{1}^{*}(\hat{T}\tilde{K}^{c}_{\ge -1}(1))$$
where $\wedge \psi_{1}^{*}(\kr \tilde{\frak
I}^{*})$ and $\wedge \psi_{1}^{*}(\hat{T}\tilde{K}^{c}_{\ge -1}(1))$ 
are the wedge product
bundles of
$\psi_{1}^{*}(\kr \tilde{\frak
I}^{*})$ and $\psi_{1}^{*}(\hat{T}\tilde{K}^{c}_{\ge -1}(1))$. 
For any $P\in K_{\ge -1}(1)$, the fiber of 
$$(\wedge \psi_{1}^{*}(\kr \tilde{\frak
I}^{*}))\wedge (\wedge \psi_{1}^{*}(\hat{T}\tilde{K}^{c}_{\ge -1}(1))$$
over $P$
is a
${\Bbb Z}\times 
{\Bbb Z}$-graded vector spaces spanned by the canonical basis
$${\Bbb L}'(i_{1})\mbar_{\tilde{\frak I}(\psi_{1}(P))}\wedge \dots \wedge
{\Bbb L}'(i_{m})\mbar_{\tilde{\frak I}(\psi_{1}(P))}\wedge \tilde{\frak
N}_{*}({\Bbb L}(j_{1}))\mbar_{\psi_{1}(P)}
\wedge \dots \wedge \tilde{\frak
N}_{*}({\Bbb L}(j_{n}))\mbar_{\psi_{1}(P)},$$
$m, n>0$, $i_{1}, \dots, i_{m}< -1$, $j_{1}, \dots, j_{n}\ge -1$.

This bundle has holomorphic sections 
$$P\mapsto 
{\Bbb L}'(i_{1})\mbar_{\tilde{\frak I}(\psi_{1}(P))}\wedge \dots \wedge
{\Bbb L}'(i_{m})\mbar_{\tilde{\frak I}(\psi_{1}(P))}\wedge 
\tilde{\frak N}_{*}({\Bbb L}(j_{1}))\mbar_{\psi_{1}(P)}
\wedge \dots \wedge \tilde{\frak N}_{*}({\Bbb L}(j_{n}))
\mbar_{\psi_{1}(P)},$$
$m, n>0$, $i_{1}, \dots, i_{m}< -1$, $j_{1}, \dots, j_{n}\ge -1$.
We shall denote the space of these holomorphic sections 
by $\wedge_{\infty}
\psi_{1}^{*}(\tilde{\frak N}_{*}({\Bbb L}))$.
It is clear that $\wedge_{\infty}
\psi_{1}^{*}(\tilde{\frak N}_{*}({\Bbb L}))$
 is naturally 
isomorphic to $\wedge_{\infty}\psi_{0}^{*}(\tilde{\frak
E}^{*}({\Bbb L}'))$. Note that these sections are the semi-infinite 
analogue of skew-symmetric holomorphic poly-vector fields. So we obtain the
other construction of $\wedge_{\infty}\psi_{0}^{*}(\tilde{\frak
E}^{*}({\Bbb L}'))$ we need.

Since the space $D_{0}$ of meromorphic functions on $K(0)$ is 
linearly isomorphic to 
$\tilde{D}_{0}^{c;1}$,  the spaces  
$D_{0}\otimes \wedge_{\infty}\psi_{0}^{*}(\tilde{\frak D}_{*}({\Bbb L}))$
and 
$D_{0}\otimes \wedge_{\infty}\psi_{0}^{*}(\tilde{\frak
E}^{*}({\Bbb L}'))$
of meromorphic sections of 
$$(\wedge \psi_{0}^{*}(\kr \tilde{\frak
E}^{*}))\wedge (\wedge \psi_{0}^{*}(\hat{T}\tilde{K}^{c}(0)))$$
and 
$$(\wedge \psi_{0}^{*}(\kr \tilde{\frak
E}^{*}))'\wedge (\wedge \psi_{0}^{*}(\hat{T}\tilde{K}^{c}(0)))',$$ 
respectively,
are linearly isomorphic to 
the space 
$\tilde{D}_{0}^{c;1}\otimes \wedge_{\infty}{\frak W}$. 
By Proposition \ref{vir-rel}, we know that $\tilde{D}_{0}^{c;1}$ 
is a module for the Virasoro algebra with central charge $c$. 
A space of the form 
$V\otimes \wedge_{\infty}{\frak W}$ where $V$ is a module for the 
Virasoro algebra has been 
studied in detail in 
the theory of the semi-infinite cohomology of the Virasoro algebra
\cite{F} \cite{FGZ}  \cite{LZ1} \cite{LZ2}. 
The space
$\wedge_{\infty}{\frak W}$ is in fact the space of 
 semi-infinite forms for the Virasoro algebra relative to 
the center. Thus we can apply their theory in our case.   
In the remaining part
of this subsection,  we recall 
some basic facts we need in the case $V=\tilde{D}_{1}^{c;1}$.  
We first discuss $\tilde{D}_{1}^{c;1}\otimes \wedge_{\infty} 
\psi_{0}^{*}(\tilde{\frak D}_{*}({\Bbb L}))$. 

We define a module structure on
$\wedge_{\infty}
\psi_{0}^{*}(\tilde{\frak D}_{*}({\Bbb L}))$ for the Virasoro algebra. 
We define linear operators
$\varepsilon(\psi_{0}^{*}({\Bbb L}'(j)))$, 
$\iota(\psi_{0}^{*}(\tilde{\frak D}_{*}({\Bbb L}(j))))$, $j\in {\Bbb Z}$, on 
$\wedge_{\infty}\psi_{0}^{*}(\tilde{\frak D}_{*}({\Bbb L}))$
 as follows: For any $m, n\ge 0$, $i_{1}, 
\dots, i_{m}\ge -1$, $j_{1}, \dots, j_{n}<-1$,
\begin{eqnarray*}
\lefteqn{\varepsilon(\psi_{0}^{*}({\Bbb L}'(j)))
(\psi_{0}^{*}({\Bbb L}'(i_{1}))\wedge \dots \wedge
\psi_{0}^{*}({\Bbb L}'(i_{m}))}\nno\\
&&\hspace{6em}\wedge\psi_{0}^{*}(\tilde{\frak D}_{*}({\Bbb L}(j_{1}))
\wedge \dots \wedge \psi_{0}^{*}(\tilde{\frak D}_{*}({\Bbb L}(j_{n})))\nno\\
&&=
\psi_{0}^{*}({\Bbb L}'(j))\wedge \psi_{0}^{*}({\Bbb L}'(i_{1}))
\wedge \dots \wedge
\psi_{0}^{*}({\Bbb L}'(i_{m}))\nno\\
&&\hspace{6em}\wedge 
\psi_{0}^{*}(\tilde{\frak D}_{*}({\Bbb L}(j_{1}))
\wedge
 \dots \wedge \psi_{0}^{*}(\tilde{\frak D}_{*}({\Bbb L}(j_{n}))
\end{eqnarray*}
for $j\ge -1$,
\begin{eqnarray*}
\lefteqn{\varepsilon(\psi_{0}^{*}({\Bbb L}'(j)))
(\psi_{0}^{*}({\Bbb L}'(i_{1}))\wedge \dots \wedge
\psi_{0}^{*}({\Bbb L}'(i_{m}))}\nno\\
&&\hspace{6em}\wedge 
\psi_{0}^{*}(\tilde{\frak D}_{*}({\Bbb L}(j_{1}))
\wedge \dots \wedge \psi_{0}^{*}(\tilde{\frak D}_{*}({\Bbb L}(j_{n})))\nno\\
&&=\sum_{k=1}^{n}(-1)^{m+k-1}\delta_{jj_{k}}
\psi_{0}^{*}({\Bbb L}'(i_{1}))\wedge \dots \wedge
\psi_{0}^{*}({\Bbb L}'(i_{m}))\wedge 
\psi_{0}^{*}(\tilde{\frak D}_{*}({\Bbb L}(j_{1}))\wedge \dots \nno\\
&&\hspace{6em} \wedge 
\widehat{\psi_{0}^{*}(\tilde{\frak D}_{*}({\Bbb L}(j_{k}))}
\wedge \dots \wedge \psi_{0}^{*}(\tilde{\frak D}_{*}({\Bbb L}(j_{n})))
\end{eqnarray*}
for $j<-1$,
\begin{eqnarray*}
\lefteqn{\iota(\psi_{0}^{*}(\tilde{\frak D}_{*}({\Bbb L}(j))))
(\psi_{0}^{*}({\Bbb L}'(i_{1}))
\wedge \dots \wedge
\psi_{0}^{*}({\Bbb L}'(i_{m}))}\nno\\
&&\hspace{6em}\wedge 
\psi_{0}^{*}(\tilde{\frak D}_{*}({\Bbb L}(j_{1})))
\wedge \dots \wedge \psi_{0}^{*}(\tilde{\frak D}_{*}({\Bbb L}(j_{n}))))\nno\\
&&=\sum_{k=1}^{m}(-1)^{k-1}\delta_{ji_{k}}
\psi_{0}^{*}({\Bbb L}'(i_{1}))\wedge \dots \wedge 
\widehat{\psi_{0}^{*}({\Bbb L}'(i_{k}))}\wedge 
\dots \wedge
\psi_{0}^{*}({\Bbb L}'(i_{m}))\nno\\
&&\quad\quad\quad\quad \wedge
 \psi_{0}^{*}(\tilde{\frak D}_{*}({\Bbb L}(j_{1})))
\wedge \dots \wedge \psi_{0}^{*}(\tilde{\frak D}_{*}({\Bbb L}(j_{n})))
\end{eqnarray*}
for $j\ge -1$, and 
\begin{eqnarray*}
\lefteqn{\iota(\psi_{0}^{*}(\tilde{\frak D}_{*}({\Bbb L}(j))))
(\psi_{0}^{*}({\Bbb L}'(i_{1}))\wedge \dots \wedge
\psi_{0}^{*}({\Bbb L}'(i_{m}))}\nno\\
&&\hspace{6em}\wedge 
\psi_{0}^{*}(\tilde{\frak D}_{*}({\Bbb L}(j_{1}))
\wedge \dots \wedge \psi_{0}^{*}(\tilde{\frak D}_{*}({\Bbb L}(j_{n})))\nno\\
&&=(-1)^{m}\psi_{0}^{*}({\Bbb L}'(i_{1}))\wedge \dots \wedge
\psi_{0}^{*}({\Bbb L}'(i_{m}))\wedge 
\psi_{0}^{*}(\tilde{\frak D}_{*}({\Bbb L}(j)))\nno\\
&&\quad\quad\quad\quad \wedge 
\psi_{0}^{*}(\tilde{\frak D}_{*}({\Bbb L}(j_{1})))
\wedge \dots \wedge \psi_{0}^{*}(\tilde{\frak D}_{*}({\Bbb L}(j_{n})))
\end{eqnarray*}
for $j<-1$.

Recall that {from} the definition, $\wedge_{\infty}
\psi_{0}^{*}(\tilde{\frak D}_{*}({\Bbb L}))$
is
${\Bbb Z}\times {\Bbb Z}$-graded. This ${\Bbb Z}\times {\Bbb Z}$-grading 
(fermion numbers and weights)
gives a ${\Bbb Z}\times {\Bbb Z}$-grading (fermion numbers and 
weights) on the space 
of linear operators on $\wedge_{\infty}
\psi_{0}^{*}(\tilde{\frak D}_{*}({\Bbb L}))$.
Given any two homogeneous operators $O_{1}$ and $O_{2}$ on 
$\wedge_{\infty}\psi_{0}^{*}(\tilde{\frak D}_{*}({\Bbb L}))$, we denote their 
fermion numbers by $|O_{1}|$ and $|O_{2}|$, respectively,
 and we define
$$[O_{1}, O_{2}]=O_{1}O_{2}-(-1)^{|O_{1}||O_{2}|}O_{2}O_{1}.$$
  Then we have:

\begin{propo}
The operators $\varepsilon(\psi_{0}^{*}({\Bbb L}'(j)))$, 
$\iota(\psi_{0}^{*}(\tilde{\frak D}_{*}({\Bbb L}(j))))$, 
$j\in {\Bbb Z}$,
are all odd operators. For any $j\in {\Bbb Z}$,
the weight of $\varepsilon(\psi_{0}^{*}({\Bbb L}'(j)))$ 
is $-j$ and the weight of 
$\iota(\psi_{0}^{*}(\tilde{\frak D}_{*}({\Bbb L}(j))))$ is $j$. 
These operators 
satisfy the following bracket formulas:
\begin{eqnarray*}
[\varepsilon(\psi_{0}^{*}({\Bbb L}'(i))), 
\varepsilon(\psi_{0}^{*}({\Bbb L}'(j)))]&=&0\\
{[\iota(\psi_{0}^{*}(\tilde{\frak D}_{*}({\Bbb L}(i)))), 
\iota(\psi_{0}^{*}(\tilde{\frak D}_{*}({\Bbb L}(j))))]}&=&0\\
{[\varepsilon(\psi_{0}^{*}({\Bbb L}'(i))), 
\iota(\psi_{0}^{*}(\tilde{\frak D}_{*}({\Bbb L}(j))))]}&=&\delta_{ij}
\end{eqnarray*}
for $i, j\in {\Bbb Z}$. \epf
\end{propo}

For $j\in {\Bbb Z}$, we define the operators
$$L_{\wedge}(j)=\sum_{i\in {\Bbb Z}}:
\varepsilon(\psi_{0}^{*}({\Bbb L}'(i)))
\iota(\psi_{0}^{*}(\tilde{\frak D}_{*}([{\Bbb L}(i), 
{\Bbb L}(j)]))):$$
on $\wedge_{\infty}\psi_{0}^{*}(\tilde{\frak D}_{*}({\Bbb L}))$,
 where the normal ordering $:\cdot:$ is defined by 
$$:\varepsilon(\psi_{0}^{*}({\Bbb L}'(i)))
\iota(\psi_{0}^{*}(\tilde{\frak D}_{*}({\Bbb L}(j)))):=
\left\{\begin{array}{ll}
-\iota(\psi_{0}^{*}(\tilde{\frak D}_{*}({\Bbb L}(j))))
\varepsilon(\psi_{0}^{*}({\Bbb L}'(i)))&j< -1\\
\varepsilon(\psi_{0}^{*}({\Bbb L}'(i)))
\iota(\psi_{0}^{*}(\tilde{\frak D}_{*}({\Bbb L}(j)))&j\ge -1.
\end{array}\right.$$

\begin{propo}
The operators $L_{\wedge}(j)$, $j\in {\Bbb Z}$, satisfy the bracket 
relations
$$[L_{\wedge}(i), L_{\wedge}(j)]=(i-j)
L_{\wedge}(i+j)+\frac{-26}{12}(i^{3}-i)
\delta_{i+j, 0}$$
for $i, j\in {\Bbb Z}$. \epf
\end{propo}

By this proposition, we see that $\wedge_{\infty}
\psi_{0}^{*}(\tilde{\frak D}_{*}({\Bbb L}))$
is a 
module of central charge $-26$ for the Virasoro algebra.
Thus the space $\tilde{D}_{0}^{c;1}\otimes \wedge_{\infty}
\psi_{0}^{*}(\tilde{\frak D}_{*}({\Bbb L}))$,
or equivalently, $D_{0}\otimes \wedge_{\infty}{\frak W}$
of meromorphic forms on $K(1)$ is a module  of
central charge $c-26$ for the Virasoro algebra.

As in the case of skew-symmetric poly-vector fields
 on a finite-dimensional manifold, we 
have a differential $\delta$ on the space 
$\tilde{D}_{0}^{c;1}\otimes \wedge_{\infty}
\psi_{0}^{*}(\tilde{\frak D}_{*}({\Bbb L}))$
of meromorphic forms on $K(1)$. In fact such a differential has been defined 
for $V\otimes \wedge_{\infty}{\frak W}$, or
equivalently, $V\otimes \wedge_{\infty}
\psi_{0}^{*}(\tilde{\frak D}_{*}({\Bbb L}))$
 for any module $V$ in the category 
${\cal O}$ for the Virasoro 
algebra. If the Virasoro operators on $V$ are denoted by
$L(j)$, $j\in {\Bbb Z}$,
then 
\begin{eqnarray*}
\delta&=&\sum_{j\in {\Bbb Z}}L(j)\otimes 
\varepsilon(\psi_{0}^{*}({\Bbb L}'(j)))\nno\\
&&-\frac{1}{2}
I_{V}\otimes \sum_{i, j\in {\Bbb Z}}:
\iota(\psi_{0}^{*}(\tilde{\frak D}_{*}([{\Bbb L}(i), {\Bbb L}(j)])))
\varepsilon(\psi_{0}^{*}({\Bbb L}'(i)))
\varepsilon(\psi_{0}^{*}({\Bbb L}'(j))):,
\end{eqnarray*}
where $I_{V}$ is the identity operator on $V$.
In the case that $V=\tilde{D}_{0}^{c;1}$, we have
\begin{eqnarray*}
\delta&=&\sum_{j\in {\Bbb Z}}{\Bbb L}(j)\otimes 
\varepsilon(\psi_{0}^{*}({\Bbb L}'(j)))\nno\\
&&-\frac{1}{2}
I_{\tilde{D}_{0}^{c;1}}\otimes \sum_{i, j\in {\Bbb Z}}:
\iota(\psi_{0}^{*}(\tilde{\frak D}_{*}([{\Bbb L}(i), {\Bbb L}(j)])))
\varepsilon(\psi_{0}^{*}({\Bbb L}'(i)))
\varepsilon(\psi_{0}^{*}({\Bbb L}'(j))):.
\end{eqnarray*}

\begin{propo}\label{delta^2=0}
When $c=26$, $\delta^{2}=0$. \epf
\end{propo}

We also have a  {\it fermion number operator} or
{\it ghost number operator}
$$
U=\sum_{j\in {\Bbb Z}}:\varepsilon(\psi_{0}^{*}({\Bbb L}'(j)))
\iota(\psi_{0}^{*}(\tilde{\frak D}_{*}({\Bbb L}(j)))):.
$$
This operator gives a ${\Bbb Z}$-grading to $\tilde{D}_{0}^{c;1}\otimes 
\psi_{0}^{*}(\tilde{\frak D}_{*}({\Bbb L}))$ and this grading is 
the same as the first ${\Bbb Z}$-grading or the grading given by 
fermion numbers.

Since $\wedge_{\infty}\psi_{0}^{*}(\tilde{\frak
E}^{*}({\Bbb L}'))$ is the graded dual of $\wedge_{\infty}
\psi_{0}^{*}(\tilde{\frak D}_{*}({\Bbb L}))$, 
the adjoint operator $d$ of $\delta$ defines a differential on 
$\wedge_{\infty}\psi_{0}^{*}(\tilde{\frak
E}^{*}({\Bbb L}'))$ and we have:

\begin{propo}\label{d^2=0}
When $c=26$, $d^{2}=0$.\epf
\end{propo}

The adjoint operator $U'$ of $U$ is called the {\it fermion number operator}
or {\it ghost number operator} on 
$\wedge_{\infty}\psi_{0}^{*}(\tilde{\frak
E}^{*}({\Bbb L}'))$, and it gives the first  ${\Bbb Z}$-grading 
or the grading given by 
fermion numbers.

Let $b(j)=
\iota(\psi_{0}^{*}(\tilde{\frak D}_{*}({\Bbb L}(j))))$ and 
$c(j)=
\varepsilon(\psi_{0}^{*}({\Bbb L}'(-j)))$ for $j\in {\Bbb Z}$
and 
let 
$b(x)=\sum_{j\in {\Bbb Z}}b(j)x^{-j-2}$ and $c(x)
=\sum_{j\in {\Bbb Z}}c(j)x^{-j+1}$.
We define a vertex operator map 
$$Y: 
\wedge_{\infty}\psi_{0}^{*}(\tilde{\frak D}_{*}({\Bbb L}))\otimes  
\wedge_{\infty}\psi_{0}^{*}(\tilde{\frak D}_{*}({\Bbb L}))
\to  \wedge_{\infty}\psi_{0}^{*}(\tilde{\frak D}_{*}({\Bbb L}))((x))$$
as follows:  We use recurrence to define
$Y$. First we define 
$Y(1, x)$ to be the identity map on 
$\wedge_{\infty}\psi_{0}^{*}(\tilde{\frak D}_{*}({\Bbb L}))$.
Assume that $u$ is of the form 
$c(i_{1})\cdots c(i_{m})1$, $m\ge 0$, $i_{1}<\cdots< i_{m}\le -1$, and
$Y(u, x)$ is already defined. 
For $v=c(i)u$,  we define 
\begin{eqnarray*}
Y(v, x)
&=&\res_{x_{1}}(x_{1}-x)^{i-2}
c(x_{1})Y(u, x)\nno\\
&& -(-1)^{|u|}\res_{x_{1}}
(-x+x_{1})^{i-2}Y(u, x)c(x_{1}).
\end{eqnarray*}
Assume that $u$ is of the form 
$b(j_{1})\cdots b(j_{n})c(i_{1})\cdots c(i_{m})1$, 
$m, n\ge 0$, $i_{1}<\cdots <i_{m}\le  -1$, $j_{1}<\cdots< j_{n}<-1$,
 and
$Y(u, x)$ is already defined. For $v=b(i)u$,  we define
\begin{eqnarray*}
Y(v, x)
&=&\res_{x_{1}}(x_{1}-x)^{i+1}
b(x_{1})Y(u, x)\nno\\
&& -(-1)^{|u|}\res_{x_{1}}
(-x+x_{1})^{i+1}Y(u, x)b(x_{1}).
\end{eqnarray*}
By definition,  elements of the form
$b(j_{1})\cdots b_{j_{n}}c(i_{1})\cdots c(i_{m})1$,
$m, n\ge 0$, $i_{1}<\cdots <i_{m}\le  -1$, $j_{1}<\cdots< j_{n}<-1$,
 form a basis of 
$\wedge_{\infty}\psi_{0}^{*}(\tilde{\frak D}_{*}({\Bbb L}))$. 
Thus the above procedure
indeed defines a linear map
$$Y: 
\wedge_{\infty}\psi_{0}^{*}(\tilde{\frak D}_{*}({\Bbb L}))\otimes  
\wedge_{\infty}\psi_{0}^{*}(\tilde{\frak D}_{*}({\Bbb L}))
\to  \wedge_{\infty}\psi_{0}^{*}(\tilde{\frak D}_{*}({\Bbb L}))((x)).$$

\begin{propo}\label{3-1}
The ${\Bbb Z}\times {\Bbb Z}$-graded vector space 
$\wedge_{\infty} \psi_{0}^{*}(\tilde{\frak D}_{*}({\Bbb L}))$, 
equipped with the vertex operator map $Y$ defined above, the vacuum
$1$ and the Virasoro element
\begin{eqnarray*}
\omega_{\wedge}&=&2c(0)b(-2)1+c(1)b(-3)1
\\
&=&2\psi_{0}^{*}({\Bbb L}'(0))\wedge 
\psi_{0}^{*}(\tilde{\frak D}_{*}({\Bbb L}(-2)))
+\psi_{0}^{*}({\Bbb L}'(-1))
\wedge \psi_{0}^{*}(\tilde{\frak D}_{*}({\Bbb L}(-3)))
\end{eqnarray*}
 is a ${\Bbb Z}\times {\Bbb Z}$-graded vertex operator
algebra (satisfying the grading-restriction conditions). Moreover,
$$Y(\omega_{\wedge}, x)=:c(x)\frac{d}{dx}b(x):
+2:\left(\frac{d}{dx}c(x)\right)b(x):.$$
\end{propo}
\pf
The proof of the first conclusion
is similar to the proof in \cite{H1} that the
quotient $M_{c, 0}/\langle L(-1)1\rangle$ 
of the Verma module $M_{c, 0}$ of lowest
weight $0$
for the Virasoro algebra  is a vertex operator algebra. We omit it here. 
The second conclusion is a direct calculation.
\epfv

The vertex operator algebra 
$\wedge_{\infty} \psi_{0}^{*}(\tilde{\frak D}_{*}({\Bbb L}))$ is sometimes 
called the {\it ghost vertex operator algebra}.

In Section 6.6 and 7.3 of \cite{H}, it was shown that the 
vertex operator algebra
$M_{c, 0}/\langle L(-1)1\rangle$  is isomorphic to a 
 geometrically constructed vertex operator algebra $M(c)$. We have:

\begin{propo}\label{3-2}
The ${\Bbb Z}$-graded space $\tilde{D}_{0}^{c;1}$ has a structure of 
 a  module for the
vertex operator algebra
$M(c)$.
\end{propo}
\pf
{From} Chapter 6 of \cite{H}, we know that 
$M(c)$ is isomorphic to the vector space spanned by the
restrictions to $\psi_{0}({\bf 0})$ of 
differential operators (of all orders) 
on the space $\tilde{D}^{c}_{0}$ of meromorphic functions 
on $\tilde{K}^{c}(0)$. We define a vertex operator map 
$$M(c)\otimes \tilde{D}_{0}^{c;1}\to \tilde{D}_{0}^{c;1}((x))$$
as follows: 
Let $v\in M(c)$ and $f\in \tilde{D}_{0}^{c;1}$. We think of $v$ as the
restriction to $\psi_{0}({\bf 0})$ of 
a derivative on $\tilde{D}^{c}_{0}$. For any $z\in {\Bbb C}^{\times}$,
let $\psi_{2}$ be the flat section on $\tilde{K}^{c}(2)$ and $P(z)$  the
element of $K(2)$ introduced in \cite{H}. Then  for fixed $\tilde{Q}_{2}\in
\tilde{K}^{c}(0)$,
$$f((\psi_{2}(P(z))_{^{1}}\widetilde{\infty}^{c/2}_{^{0}}
\tilde{Q}_{1})_{^{1}}\widetilde{\infty}^{c/2}_{^{0}}\tilde{Q}_{2})$$
as a function of $\tilde{Q}_{1}\in \tilde{K}^{c}(0)$
is  analytic. It is clear that $v$ can always
be extended so that it acts on analytic functions, not only
on meromorphic functions. Thus $v$ acts
on this function of $\tilde{Q}_{1}$ and 
$$v(f((\psi_{2}(P(z))_{^{1}}\widetilde{\infty}^{c/2}_{^{0}}
\tilde{Q}_{1})_{^{1}}\widetilde{\infty}^{c/2}_{^{0}}\tilde{Q}_{2}))$$
as a function of $\tilde{Q}_{2}\in
\tilde{K}^{c}(0)$ is analytic. In addition, this function is  also analytic
 in $z$. Thus we can expand it as a Laurent series in $z$.  It
is easy to see {from} the definition of the sewing operation that this
series has a pole at $z=0$ and the coefficients of this expansion are
meromorphic functions of  $\tilde{Q}_{2}\in 
\tilde{K}^{c}(0)$. In addition, since $f$ is
proportional to $C$, so are all the coefficients of this
expansion. Thus these coefficients give an element of
$\tilde{D}_{0}^{c;1}((x))$. We define the image $Y(v, x)f$ 
of $v\otimes f$ under
$Y$ to be
this element. 

Since our vertex operator map is defined geometrically, it is now easy
to use the geometry and the geometric proof that 
$M(c)$ is a vertex operator algebra in \cite{H} to prove that
$\tilde{D}^{c; 1}_{0}$ with the vertex operator map just defined is
indeed a module for $M(c)$.
\epfv

By Propositions \ref{3-1} and \ref{3-2}, we have:
\begin{corol}
The ${\Bbb Z}\times {\Bbb Z}$-graded vector space
$\tilde{D}^{c; 1}_{0}\otimes \wedge_{\infty}
\psi_{0}^{*}(\tilde{\frak D}_{*}({\Bbb L}))$ 
is a 
module for the 
tensor product vertex operator 
algebra $M(c)\otimes \wedge_{\infty}
\psi_{0}^{*}(\tilde{\frak D}_{*}({\Bbb L}))$. \epf
\end{corol}

We shall denote the vertex operator algebra 
$M(c)\otimes \wedge_{\infty}
\psi_{0}^{*}(\tilde{\frak D}_{*}({\Bbb L}))$ by $G$. 
We shall use the same notation to denote the vertex operator maps
for $\wedge_{\infty}
\psi_{0}^{*}(\tilde{\frak D}_{*}({\Bbb L}))$,  for
$G$ and for
the module $\tilde{D}^{c; 1}_{0}\otimes \wedge_{\infty}
\psi_{0}^{*}(\tilde{\frak D}_{*}({\Bbb L}))$.

\begin{corol}
Let 
$$b=b(-2)1=\psi_{0}^{*}(\tilde{\frak D}_{*}({\Bbb L}(-2)))$$
and 
$$c=c(1)1=\psi_{0}^{*}({\Bbb L}'(-1)).$$ 
Then
\begin{eqnarray*}
b(x)&=&Y(b, x),\\
c(x)&=&Y(c, x).\epfd
\end{eqnarray*}
\end{corol}

\begin{propo}
Let 
\begin{eqnarray*}
q_{G}&=&{\Bbb L}(-2)1\otimes c
+1\otimes (b(-2)c(1)c(0)1)\\
&=&{\Bbb L}(-2)1\otimes c+1\otimes 
\psi_{0}^{*}({\Bbb L}'(0))\wedge \psi_{0}^{*}({\Bbb L}'(-1))\wedge 
\psi_{0}^{*}(\tilde{\frak D}_{*}({\Bbb L}(-2))),\\
f_{G}&=&1\otimes (c(1)b(-2)1)\\
&=&1\otimes (\psi_{0}^{*}({\Bbb L}'(-1))\wedge 
\psi_{0}^{*}(\tilde{\frak D}_{*}({\Bbb L}(-2)))),\\
\omega_{G}&=&(L(-2)1)\otimes 1+1\otimes \omega_{\wedge}
\end{eqnarray*}
in $G$. Then 
$\omega_{G}$ is the Virasoro element of $G$ and in 
$\tilde{D}^{c; 1}_{0}\otimes \wedge_{\infty}
\psi_{0}^{*}(\tilde{\frak D}_{*}({\Bbb L}))$,
\begin{eqnarray*}
Y(q_{G}, x)&=&:L(x)c(x):+:b(x)c(x)\frac{d}{dx}c(x):,\\
\delta&=&\res_{x}Y(q_{G}, x),\\
Y(f_{G}, x)&=&:c(x)b(x):,\\
U&=&\res_{x}Y(f_{G}, x).\epfd
\end{eqnarray*} 
\end{propo}

By definition, we know that the ${\Bbb Z}\times 
{\Bbb Z}$-graded vertex operator algebra 
$\psi_{0}^{*}(\tilde{\frak D}_{*}({\Bbb L}))$ is generated by $b$ and $c$. 

For $j\in {\Bbb Z}$, let $b'(-j)$ and $c'(-j)$ be the adjoint operators
of $b(j)$ and $c(j)$, respectively.
Since $\wedge_{\infty}\psi_{0}^{*}(\tilde{\frak
E}^{*}({\Bbb L}'))$ is the graded dual of $\wedge_{\infty}
\psi_{0}^{*}(\tilde{\frak D}_{*}({\Bbb L}))$, $b'(-j)$ and $c'(-j)$, 
$j\in {\Bbb Z}$, act on $\wedge_{\infty}\psi_{0}^{*}(\tilde{\frak
E}^{*}({\Bbb L}'))$, and we have:

\begin{propo}\label{3-3}
The ${\Bbb Z}\times {\Bbb Z}$-graded vector space
$\wedge_{\infty}\psi_{0}^{*}(\tilde{\frak
E}^{*}({\Bbb L}'))$ is a module for the 
${\Bbb Z}\times {\Bbb Z}$-graded vertex operator 
algebra $\wedge_{\infty}
\psi_{0}^{*}(\tilde{\frak D}_{*}({\Bbb L}))$.
The ${\Bbb Z}\times {\Bbb Z}$-graded vector space
$\tilde{D}^{c; 1}_{0}\otimes \wedge_{\infty}\psi_{0}^{*}(\tilde{\frak
E}^{*}({\Bbb L}'))$ 
is a 
module for $G$. If we still use $Y$ to denote
the vertex operator map for the second module, then
\begin{eqnarray*}
Y(b, x)&=&\sum_{j\in {\Bbb Z}}b'(-j)x^{-x-2},\\
Y(c, x)&=&\sum_{j\in {\Bbb Z}}c'(-j)x^{-x+1},\\
d&=&\res_{x}Y(q, x),\\
U'&=&\res_{x}Y(f, x). \;\;\;\;\;\;\;\;\;\;\;\;\;\;\Box
\end{eqnarray*}
\end{propo}

\subsection{Semi-infinite forms on $K(n)$, $n\ge 0$}

In this subsection, we introduce semi-infinite forms on $K(n)$ for 
$n\ge 0$.

We fix $n\ge 0$. {From} the definition of $K(n)$, it is easy to see 
that there is a canonical injective 
map ${\frak K}_{n}: K(n)\to K(0)\times (K_{\ge -1}(1))^{n}$ and 
${\frak K}_{n}(K(n))$ is in fact an open subset of 
$K(0)\times (K_{\ge -1}(1))^{n}$. 
Over $K(0)$ and $K_{\ge -1}(1)$, we have the holomorphic bundles
$$(\wedge \psi_{0}^{*}(\kr \tilde{\frak
E}^{*}))\wedge (\wedge \psi_{0}^{*}(\hat{T}\tilde{K}^{c}(0)))$$
and 
$$(\wedge \psi_{1}^{*}(\kr \tilde{\frak
I}^{*}))\wedge (\wedge \psi_{1}^{*}(\hat{T}\tilde{K}^{c}_{\ge -1}(1)),$$
respectively. Thus we have 
the exterior product bundle (not the wedge product bundle)
\begin{equation}\label{4.0.1}
(\wedge \psi_{0}^{*}(\kr \tilde{\frak
E}^{*}))\wedge (\wedge \psi_{0}^{*}(\hat{T}\tilde{K}^{c}(0)))
\boxtimes ((\wedge \psi_{1}^{*}(\kr \tilde{\frak
I}^{*}))\wedge 
(\wedge \psi_{1}^{*}(\hat{T}\tilde{K}^{c}_{\ge -1}(1)))^{\boxtimes n}
\end{equation}
over $K(0)\times (K_{\ge -1}(1))^{n}$. 
The map ${\frak K}_{n}$ pulls this  holomorphic bundle 
back to a holomorphic bundle over $K(n)$. We denote this 
holomorphic bundle over $K(n)$ by ${\cal G}(n)$. 

By the definition of exterior product bundle, we have 
linear injective maps {from} the space of sections of 
$(\wedge \psi_{0}^{*}(\kr \tilde{\frak
E}^{*}))\wedge (\wedge \psi_{0}^{*}(\hat{T}\tilde{K}^{c}(0)))$
and the spaces of sections of copies of 
$(\wedge \psi_{1}^{*}(\kr \tilde{\frak
I}^{*}))\wedge 
(\wedge \psi_{1}^{*}(\hat{T}\tilde{K}^{c}_{\ge -1}(1))$ 
to the space of 
sections of the bundle (\ref{4.0.1}). Moreover, the images of 
these spaces of sections in the space of 
sections of (\ref{4.0.1}) intersect with each other at $0$. 
Let $\wedge_{\infty}
\psi_{0}^{*}(\tilde{\frak D}_{*}({\Bbb L}^{(0)}))$ 
be a copy of $\wedge_{\infty}
\psi_{0}^{*}(\tilde{\frak D}_{*}({\Bbb L}))$ and
for $i=1, \dots, n$, let $\wedge_{\infty}
\psi_{1}^{*}(\tilde{\frak N}_{*}({\Bbb L}^{(i)}))$
be a copy of $\wedge_{\infty}
\psi_{1}^{*}(\tilde{\frak N}_{*}({\Bbb L}))$.
Then we have a linear injective map {from} 
\begin{equation}\label{4.0.2}
\wedge_{\infty}
\psi_{0}^{*}(\tilde{\frak D}_{*}({\Bbb L}^{(0)}))\wedge  
(\wedge_{i=1}^{n} \wedge_{\infty}
\psi_{1}^{*}(\tilde{\frak N}_{*}({\Bbb L}^{(i)})))
\end{equation}
to the space of 
sections of (\ref{4.0.1}). 
Since ${\frak K}_{n}$ is injective and 
${\frak K}_{n}(K(n))$ is in fact an open subset of 
$K(0)\times (K_{\ge -1}(1))^{n}$, 
${\frak K}_{n}^{*}$ is an isomorphism {from} the 
space of holomorphic sections of (\ref{4.0.1}) 
to the space of holomorphic sections of   ${\cal G}(n)$. Thus we obtain a 
linear injective map {from} (\ref{4.0.2})
to the space of holomorphic sections of ${\cal G}(n)$. 
We denote the image of this map by
$\hat{\Gamma}({\cal G}(n))$. Then $\hat{\Gamma}({\cal G}(n))$ is 
isomorphic to (\ref{4.0.2}). The tensor product space $D_{n}\otimes 
\hat{\Gamma}({\cal G}(n))$
is the semi-infinite analogue on $K(n)$ of skew-symmetric holomorphic 
poly-vector fields
on a complex manifold.

In the preceding subsection, we have
shown that $\wedge_{\infty}
\psi_{1}^{*}(\tilde{\frak N}_{*}({\Bbb L}))$ is naturally
isomorphic to $\wedge_{\infty}\psi_{0}^{*}(\tilde{\frak
E}^{*}({\Bbb L}'))$. Thus if for $i=1, \dots, n$, we use 
$\wedge_{\infty}\psi_{0}^{*}(\tilde{\frak
E}^{*}(({\Bbb L}^{(i)})'))$ to denote a copy 
$\wedge_{\infty}\psi_{0}^{*}(\tilde{\frak
E}^{*}({\Bbb L}'))$, then $\hat{\Gamma}({\cal G}(n))$ is 
isomorphic to 
\begin{equation}\label{4.0.3}
\wedge_{\infty}
\psi_{0}^{*}(\tilde{\frak D}_{*}({\Bbb L}^{(0)}))\wedge 
(\wedge_{i=1}^{n} \wedge_{\infty}\psi_{0}^{*}(\tilde{\frak
E}^{*}(({\Bbb L}^{(i)})'))).
\end{equation}
We shall  
 identify $\hat{\Gamma}({\cal G}(n))$ with this space.

The space $\hat{\Gamma}({\cal G}(n))$ is ${\Bbb Z}\times {\Bbb Z}$-graded.
We denote the graded dual of $\hat{\Gamma}({\cal G}(n))$ 
by $\hat{\Omega}(K(n))$. The tensor product space 
$D_{n}\otimes \hat{\Omega}(K(n))$ is the semi-infinite analogue on $K(n)$
of the space of holomorphic forms on a complex manifold.

We now define differentials on the spaces $D_{n}\otimes 
\hat{\Gamma}({\cal G}(n))$ and $D_{n}\otimes \hat{\Omega}(K(n))$. We 
discuss $D_{n}\otimes 
\hat{\Gamma}({\cal G}(n))$ first.
Since $\hat{\Gamma}({\cal G}(n))$
is isomorphic to (\ref{4.0.3}) and we have a differential $\delta$
for $\wedge_{\infty}
\psi_{0}^{*}(\tilde{\frak D}_{*}({\Bbb L}))$ and a differential 
$d$ for $\wedge_{\infty}\psi_{0}^{*}(\tilde{\frak
E}^{*}(({\Bbb L}^{(i)})'))$, 
we have a differential $\delta_{n}$ for $\hat{\Gamma}({\cal G}(n))$.
Since ${\frak K}_{n}(K(n))$ is in fact an open subset of 
$K(0)\times (K_{\ge -1}(1))^{n}$, tangent fields on 
$K(0)\times (K_{\ge -1}(1))^{n}$ can be restricted to 
${\frak K}_{n}(K(n))$ to obtain tangent fields on ${\frak K}_{n}(K(n))$
which can be further pulled back to tangent fields on $K(n)$. 
On $K(0)\times (K_{\ge -1}(1))^{n}$, we have the tangent fields 
${\Bbb L}^{(0)}(j)$, $j<-1$, which are the push-forwards of 
the tangent  fields ${\Bbb L}(j)$ on $K(0)$ and for $i=1, \dots, n$, we
have the tangent fields 
${\Bbb L}^{(i)}(j)$, $j\ge -1$, 
which are the push-forwards of 
the tangent fields ${\Bbb L}(j)$ on the $i$-th copy of $K_{\ge -1}(1)$.
We shall use the same notations ${\Bbb L}^{(0)}(j)$, $j<-1$,
and ${\Bbb L}^{(i)}(j)$, $j\ge -1$, $i=1, \dots, n$,
 to denote the pullbacks to 
$K(n)$ of their
restrictions to ${\frak K}_{n}(K(n))$. Recall {from} the preceding 
subsection that on $\wedge_{\infty}
\psi_{0}^{*}(\tilde{\frak D}_{*}({\Bbb L}))$ and 
$\wedge_{\infty}\psi_{0}^{*}(\tilde{\frak
E}^{*}(({\Bbb L})'))$, we have the operators $c(j)$ and $c'(-j)$,
respectively, $j\in {\Bbb Z}$. We denote $c(j)$, $j\in {\Bbb Z}$, 
on $\wedge_{\infty}
\psi_{0}^{*}(\tilde{\frak D}_{*}({\Bbb L}^{(0)}))$ by $c^{(0)}(j)$ and 
$c'(-j)$, $j\in {\Bbb Z}$, on $\wedge_{\infty}\psi_{0}^{*}(\tilde{\frak
E}^{*}(({\Bbb L}^{(i)})'))$ for $i=1, \dots, n$ by $(c')^{(i)}(-j)$.
The operators $c^{(0)}(j)$, $j<-1$, 
and for $i=1, \dots, n$, the operators $(c')^{(i)}(j)$, $j\ge -1$, 
act on $\hat{\Gamma}({\cal G}(n))$ naturally, and 
we shall use the same notation to 
denote their actions on $\hat{\Gamma}({\cal G}(n))$.

We define a differential on $D_{n}\otimes 
\hat{\Gamma}({\cal G}(n))$, still denoted $\delta_{n}$, 
by 
$$\delta_{n}=\sum_{j<-1}{\Bbb L}^{(0)}(j)\otimes c^{(0)}(-j)
+\sum_{i=1}^{n}\sum_{j\ge -1}{\Bbb L}^{(i)}(j)\otimes (c')^{(i)}(j)
+I_{D_{n}}\otimes \delta_{n}.$$

{From} this definition and Propositions \ref{delta^2=0} and \ref{d^2=0},
we conclude:

\begin{propo}
When $c=26$, $\delta_{n}^{2}=0$.\epf
\end{propo}

We define a differential $d_{n}$ on $D_{n}\otimes 
\hat{\Omega}(K(n))$ using 
$\delta_{n}$ in the obvious way. Then we have:

\begin{propo}
When $c=26$, $d_{n}^{2}=0$.\epf
\end{propo}

\subsection{A partial operad for semi-infinite forms}

In this subsection, 
we construct a partial operad ${\frak G}$ 
{from} semi-infinite forms introduced in 
the preceding two subsections. We also construct a partial suboperad ${\cal
T}_{G}$ of ${\frak G}$ which is the motivation of the construction of 
the (strong) topological vertex partial operad in Subsection 4.2. 

In \cite{H}, $M(c)$ is constructed 
as a space of 
differential operators (of all orders) on $\tilde{K}^{c}(0)$ restricted 
to $\psi_{0}({\bf 0})$. Using the tangent fields 
$\tilde{\frak D}_{*}({\Bbb L}(j))$, $j<-1$,  we see that
$M(c)$ is isomorphic to a space of differential operators (without 
restricting to $\psi_{0}({\bf 0})$). The union of the vector spaces of 
restrictions of elements of this space to 
$Q$ for all $Q\in \tilde{K}^{c}(0)$ is a
holomorphic bundle over $\tilde{K}^{c}(0)$ and $M(c)$ is a space 
of holomorphic sections of this holomorphic bundle. This holomorphic bundle
is pulled back to a holomorphic bundle ${\cal O}(0)$ over $K(0)$ and 
$M(c)$ is in fact isomorphic to a  space of 
holomorphic sections
of this pullback bundle ${\cal O}(0)$. 
{From} now on we shall identify $M(c)$ with this 
space of holomorphic
sections.

We now consider the space of differential operators (of all orders) on
$\tilde{K}^{c}_{\ge -1}(1)$ generated by $\tilde{\frak N}_{*}({\Bbb
L}(j))$, $j\ge -1$. The union of the vector spaces of 
restrictions of elements of this space to 
$Q$ for all $Q\in \tilde{K}^{c}_{\ge -1}(1)$ is a
holomorphic bundle over $\tilde{K}^{c}_{\ge -1}(1)$ and this 
space is a  space of holomorphic 
sections of this bundle. This bundle and
this space of holomorphic sections are pulled back by $\psi_{1}^{*}$
to a holomorphic bundle ${\cal O}_{\ge -1}(1)$ over $K_{\ge -1}(1)$
and a space $M_{\ge -1}(c)$ of holomorphic sections of ${\cal O}_{\ge
-1}(1)$, respectively.  Note that $M(c)$ and $M_{\ge -1}(c)$ are in fact
isomorphic to the universal enveloping algebras of the Lie subalgebras
of the Virasoro algebra spanned by $L(j)$ for $j<-1$ and by $L(j)$ for
$j\ge -1$, respectively.
Since $M(c)$ is a module for the Virasoro algebra, both
$M(c)$ and $M_{\ge -1}(c)$ act on $M(c)$ naturally. These actions can
also be interpreted geometrically as follows: The sewing operation
gives the map $_{^{1}}\infty_{^{0}}: K_{\ge -1}(1)\times K(0)\to
K(0)$. Note that elements of $M_{\ge -1}(c)\otimes M(c)$ can be viewed
as operators on functions on $K_{\ge -1}(1)\times K(0)$.  The map
$_{^{1}}\infty_{^{0}}$ pushes the restriction to $(I, {\bf 0})\in
K_{\ge -1}(1)\times K(0)$ of an operator in $M_{\ge -1}(c)\otimes
M(c)$ to an operator on the space of functions holomorphic at ${\bf
0}\in K(0)$. Using the left multiplication in $K(0)$ defined by the 
embedding {from} $K(0)$ to $K(1)$ and the sewing operation, 
we see that this push-forward operator is the restriction to ${\bf
0}$ of an operator in $M(c)$. Thus we obtain a map {from} $M_{\ge
-1}(c)\otimes M(c)$ to $M(c)$. This map is the action of $M_{\ge
-1}(c)$ on $M(c)$. Similarly, the action of $M(c)$ on $M(c)$ has a
geometric meaning. For the details of the geometric construction of
the action of the Virasoro algebra on $M(c)$,
see \cite{H}.

Both $M(c)$ and $M_{\ge -1}(c)$ are ${\Bbb Z}$-graded (by weights).
Note that the homogeneous subspaces of $M(c)$ of weights less than $0$
are all $0$ and homogeneous subspaces of $M_{\ge -1}(c)$ of weights
larger than $1$ are all $0$. 
Also note that the homogeneous
subspaces of $\wedge_{\infty} \psi_{0}^{*}(\tilde{\frak D}_{*}({\Bbb
L}))$ of weights less than $0$ are all $0$ and homogeneous subspaces of
$\wedge_{\infty} \psi_{1}^{*}(\tilde{\frak N}_{*}({\Bbb L}))$ of weights
larger than $1$ are all $0$. 

{From} the definition, we see that the space 
$\wedge_{\infty}
\psi_{1}^{*}(\tilde{\frak N}_{*}({\Bbb L}))$ is isomorphic to 
the space  spanned by operators on $\wedge_{\infty}
\psi_{0}^{*}(\tilde{\frak D}_{*}({\Bbb L}))$ of the form
$$b(i_{1})\cdots b(i_{m})c(-j_{1})\cdots c(-j_{n}),$$
$m, n>0$, $i_{1}, \dots, i_{m}\ge -1$, $j_{1}, \dots, j_{n}<-1$.
We shall identify $\wedge_{\infty}
\psi_{1}^{*}(\tilde{\frak N}_{*}({\Bbb L}))$ with this space of operators.
The space $\wedge_{\infty}
\psi_{0}^{*}(\tilde{\frak D}_{*}({\Bbb L}))$ is isomorphic to
the space spanned by operators on itself of the form 
$$c(-i_{1})\cdots c(-i_{m})b(j_{1})\cdots b(j_{n}),$$
$m, n>0$, $i_{1}, \dots, i_{m}\ge -1$, $j_{1}, \dots, j_{n}<-1$.
We shall identify $\wedge_{\infty}
\psi_{0}^{*}(\tilde{\frak D}_{*}({\Bbb L}))$ with this space of operators.

Consider the tensor products 
$$T_{\wedge}=M(c)\otimes
\wedge_{\infty} \psi_{0}^{*}(\tilde{\frak D}_{*}({\Bbb
L}))$$ and 
$$T_{\wedge}^{\ge -1}=M_{\ge -1}(c)\otimes \wedge_{\infty} 
\psi_{1}^{*}(\tilde{\frak N}_{*}({\Bbb L})).$$
Then $T_{\wedge}$ is a ${\Bbb Z}\times {\Bbb Z}$-graded 
vertex operator algebra and both 
$T_{\wedge}$ and $T_{\wedge}^{\ge -1}$ can be identified with 
spaces of operators on $T_{\wedge}$ spanned by homogeneous vectors. The
homogeneous subspaces of $T_{\wedge}$ of weights less than $0$
are all $0$ and homogeneous subspaces of $T_{\wedge}^{\ge -1}$ of weights
larger than $1$ are all $0$. We have the algebraic completions
$\overline{T}_{\wedge}$ and $\overline{T}_{\wedge}^{\ge -1}$ of 
$T_{\wedge}$ and $T_{\wedge}^{\ge -1}$, respectively. When 
$c=26$, the central charge of $T_{\wedge}$ is $0$. From now on, we
shall take $c=26$.

Since $M(c)$ and $\wedge_{\infty} 
\psi_{0}^{*}(\tilde{\frak D}_{*}({\Bbb
L}))$ ($M_{\ge -1}(c)$ and $\wedge_{\infty} 
\psi_{1}^{*}(\tilde{\frak N}_{*}({\Bbb L}))$) are
spaces of holomorphic sections of holomorphic vector bundles over 
$K(0)$ ($K_{\ge -1}(1)$), the unions of the vector spaces of 
restrictions of elements of $T_{\wedge}$ and 
$T_{\wedge}^{\ge -1}$ to 
$Q$ for all $Q\in \tilde{K}^{c}_{\ge -1}(1)$ form
holomorphic bundles ${\frak G}(0)$ and ${\frak G}_{\ge -1}(1)$
over 
$K(0)$ and $K_{\ge -1}(1)$, respectively, and $T_{\wedge}$ and 
$T_{\wedge}^{\ge -1}$ are spaces of holomorphic 
sections of ${\frak G}(0)$ and ${\frak G}_{\ge -1}(1)$, respectively.

In \cite{H6} and \cite{H7}, a locally convex topological completion 
$H^{V}\subset \overline{V}$
of a vertex operator algebra $V$ is constructed such that 
the maps in $\hom(V^{\otimes n}, \overline{V})$
corresponding to elements of $K_{{\frak H}_{1}}(n)$  
 can in fact be extended to
maps in $\hom((H^{V})^{\otimes n}, H^{V})$. In particular, 
when the central charge of $V$ is $0$, $H^{V}$
has a structure of an algebra over the (non-partial) 
suboperad $K_{{\frak H}_{1}}$ 
of $K$. These results generalize to ${\Bbb Z}\times {\Bbb Z}$-graded 
vertex operator algebras
in an obvious way. 

Applying this result to
the vertex operator algebra $T_{\wedge}$ of central charge $0$, 
we obtain a locally 
convex completion $H^{T_{\wedge}}\subset \overline{T}_{\wedge}$ 
such that $H^{T_{\wedge}}$
has a structure of an algebra over  $K_{{\frak H}_{1}}$. 
Elements of $\hom((H^{T_{\wedge}})^{\otimes n}, H^{T_{\wedge}})$
give elements of $\hom(T_{\wedge}^{\otimes n}, 
\overline{T}_{\wedge})$. Thus we have a map {from} 
$\hom((H^{T_{\wedge}})^{\otimes n}, H^{T_{\wedge}})$ to 
$\hom(T_{\wedge}^{\otimes n}, 
\overline{T}_{\wedge})$. In fact this map is linear and injective.
{From} now on, we shall view $\hom((H^{T_{\wedge}})^{\otimes n}, 
H^{T_{\wedge}})$ as a subspace of $\hom(T_{\wedge}^{\otimes n}, 
\overline{T}_{\wedge})$.

Note that since elements of $T_{\wedge}$ and $T_{\wedge}^{\ge -1}$
have been identified with spaces of operators on $T_{\wedge}$,
elements of the algebraic completions
$\overline{T}_{\wedge}$ and $\overline{T}_{\wedge}^{\ge -1}$
can be identified with elements of $\hom(T_{\wedge}, 
\overline{T}_{\wedge})$. 
Let $\widetilde{T}_{\wedge}$ ($\widetilde{T}_{\wedge}^{\ge -1}$)
be the subspaces of   $\overline{T}_{\wedge}$ 
($\overline{T}_{\wedge}^{\ge -1}$)
consisting of elements $u\in \overline{T}_{\wedge}$ ($u\in 
\widetilde{T}_{\wedge}^{\ge -1}$) satisfying the 
following property: There 
exists a positive  number $a$ such that $a^{L(0)}ua^{-L(0)}$ 
as an element 
of $\hom(T_{\wedge}, 
\overline{T}_{\wedge})$
is in fact in the subspace $\hom(H^{T_{\wedge}}, H^{T_{\wedge}})$. 

Using the same construction of 
 ${\frak G}(0)$ and ${\frak G}_{\ge -1}(1)$ in terms of 
$T_{\wedge}$ and 
$T_{\wedge}^{\ge -1}$, respectively, we construct 
holomorphic vector bundles 
$\widetilde{\frak G}(0)$ and $\widetilde{\frak G}_{\ge -1}(1)$
over $K(0)$ and $K_{\ge -1}(1)$, respectively, such that 
$\widetilde{T}_{\wedge}$ and $\widetilde{T}_{\wedge}^{\ge -1}$ 
are spaces of holomorphic sections of $\widetilde{\frak G}(0)$ and
$\widetilde{\frak G}_{\ge -1}(1)$, respectively.

We have the exterior product bundle 
$$\widetilde{\frak G}(0)\boxtimes 
(\widetilde{\frak G}_{\ge -1}(1))^{\boxtimes n}$$
over $K(0)\times (K(1))^{n}$.
The space $\widetilde{T}_{\wedge}$ can be embedded 
into the space of holomorphic sections
of this bundle and 
there are $n$ embeddings of 
$\widetilde{T}_{\wedge}^{\ge -1}$ into the space 
of holomorphic sections of this bundle.
The images of these embeddings 
can be pulled back to subspaces of the space of holomorphic 
sections of the
holomorphic  bundle
$${\frak K}_{n}^{*}(\widetilde{\frak G}(0)\boxtimes 
(\widetilde{\frak G}_{\ge -1}(1))^{\boxtimes n}).$$
We shall denote these pullbacks by $\widetilde{T}^{(0)}_{\wedge}$
and 
by $\widetilde{T}^{(i)}_{\wedge}$, $i=1, \dots, n$.
Let 
$$\widetilde{\cal T}_{\wedge}(n)
={\frak K}_{n}^{*}(\widetilde{T}^{(0)}_{\wedge}
\otimes \otimes_{i=1}^{n}\widetilde{T}^{(i)}_{\wedge}).$$
Then this subspace  of the space of 
holomorphic sections of 
$${\frak K}_{n}^{*}(\widetilde{\frak G}(0)\boxtimes 
(\widetilde{\frak G}_{\ge -1}(1))^{\boxtimes n})$$
gives a holomorphic 
bundle $\widetilde{\frak G}(n)$ over $K(n)$ for each $n\ge 0$.

We now construct a partial operad structure on 
$\widetilde{\frak G}=\{\widetilde{\frak G}(n)\}_{n\ge 0}$.
We first define the composition maps for $\widetilde{\frak G}$.

We need some notations. Let 
$$Q=(z_{1}, \dots, z_{n-1}; A^{(0)}, (a_{0}^{(1)}, A^{(1)}),
\dots, (a_{0}^{(n)}, A^{(n)}))\in K(n).$$
We shall use $z_{1}(Q), \dots, z_{n}(Q)$,
$A^{(0)}(Q), \dots, A^{(n)}(Q)$ and $a_{0}^{(1)}(Q), \dots,
a_{0}^{(n)}(Q)$ to denote  $z_{1}, \dots, z_{n}$,
$A^{(0)}, \dots, A^{(n)}$ and $a_{0}^{(1)}, \dots,
a_{0}^{(n)}$, respectively. Also recall the notations 
\begin{eqnarray*}
e^{L^{+}(A)}&=&\exp\left(\sum_{j>0}A_{j}L(j)\right),\\
e^{L^{-}(A)}&=&\exp\left(\sum_{j>0}A_{j}L(-j)\right)
\end{eqnarray*}
in \cite{H} for any sequence $A\in {\Bbb C}^{\infty}$.

Let $n_{1}>0$, $n_{2}\ge 0$ and $1\le i\le n_{1}$. Let 
$(Q_{1}, \alpha_{1})\in \widetilde{\frak G}(n_{1})$ and 
$(Q_{2}, \alpha_{2})\in \widetilde{\frak G}(n_{2})$ where $Q_{1}\in K(n_{1})$,
$Q_{2}\in K(n_{2})$ and $\alpha_{1}$ and $\alpha_{2}$ 
are elements of the fibers at $Q_{1}$ and $Q_{2}$, respectively,
of $\widetilde{\frak G}(n_{1})$ and $\widetilde{\frak G}(n_{2})$.
We would like to define 
$(Q_{1}, \alpha_{1})\circ_{i}(Q_{2}, \alpha_{2})$ when 
$Q_{1}{}_{^{i}}\infty_{^{0}}Q_{2}$ exists and an 
additional condition (see below) on $\alpha_{1}$ and $\alpha_{2}$ holds.

By the definition of $\widetilde{\frak G}(n_{1})$ and 
$\widetilde{\frak G}(n_{2})$,
we need only consider those $\alpha_{1}$ and $\alpha_{2}$ which 
are the restrictions to $Q_{1}$ and $Q_{2}$, respectively,
of the sections 
of the forms 
$${\frak K}_{n_{1}}^{*}(u^{(0)}\otimes
u^{(1)}\otimes \cdots\otimes u^{(n_{1})})$$
and 
$${\frak K}_{n_{2}}^{*}(v^{(0)}\otimes v^{(1)}
\otimes \cdots\otimes v^{(n_{2})}),$$
respectively,
where 
$u^{(0)}, v^{(0)}\in \widetilde{T}_{\wedge}$, and
$u^{(k)}, v^{(l)}\in \widetilde{T}_{\wedge}^{\ge -1}$
for $k=1, \dots, n_{1}$ and $l=1, \dots, n_{2}$. 

{From} the definition of sewing operation,
we know that it is possible to find positive numbers 
$r_{0}, \dots, r_{n_{1}}, 
s_{0}, \dots, s_{n_{2}}$ such that
\begin{eqnarray*}                      
\lefteqn{(\cdots (
(({\bf 0}, (r_{0}, {\bf 0}))_{^{1}}\infty_{^{0}}Q_{1})_{^{1}}\infty_{^{0}}
({\bf 0}, (r_{1}, {\bf 0})))_{^{1}}\infty_{^{0}}\cdots)_{^{1}}\infty_{^{0}}
({\bf 0}, (r_{n_{1}}, {\bf 0})))}\nno\\
&&=(\cdots (
((A^{(0)}, (a_{0}^{(0)}, {\bf 0}))_{^{1}}
\infty_{^{0}}\hat{Q}_{1})_{^{1}}\infty_{^{0}}
({\bf 0}, (a^{(1)}_{0}, A^{(1)})))_{^{1}}
\infty_{^{0}}\cdots)\nno\\
&&\hspace{8em} _{^{1}}\infty_{^{0}}
({\bf 0}, (a_{0}^{(n_{1})}, A^{(n_{1})})))\in K_{{\frak H}_{1}}(n_{1})
\end{eqnarray*}
and 
\begin{eqnarray*}                      
\lefteqn{(\cdots (
(({\bf 0}, (s_{0}, {\bf 0}))_{^{1}}\infty_{^{0}}Q_{2})_{^{1}}\infty_{^{0}}
({\bf 0}, (s_{1}, {\bf 0})))_{^{1}}\infty_{^{0}}\cdots)_{^{1}}\infty_{^{0}}
({\bf 0}, (s_{n_{2}}, {\bf 0})))}\nno\\
&&=(\cdots (
((B^{(0)}, (b_{0}^{(0)}, {\bf 0}))_{^{1}}
\infty_{^{0}}\hat{Q}_{2})_{^{1}}\infty_{^{0}}
({\bf 0}, (b^{(1)}_{0}, B^{(1)})))_{^{1}}
\infty_{^{0}}\cdots)_{^{1}}\nno\\
&&\hspace{8em} \infty_{^{0}}
({\bf 0}, (b_{0}^{(n_{2})}, B^{(n_{2})})))\in K_{{\frak H}_{1}}(n_{2}),
\end{eqnarray*}
where 
$$(A^{(0)}, (a_{0}^{(0)}, {\bf 0})),
({\bf 0}, (a^{(k)}_{0}, A^{(k)})), 
(B^{(0)}, (b_{0}^{(0)}, {\bf 0})),\;\;\;({\bf 0}, (b^{(l)}_{0}, B^{(l)}))\in 
K_{{\frak H}_{1}}(1),$$
for $k=1, \dots, n_{1}$,
$l=1, \dots, n_{2}$, and
\begin{eqnarray*}
\hat{Q}_{1}&=&(\xi_{1}, \dots, \xi_{n_{1}}; {\bf 0}, ({\bf 0},
(c_{0}^{(1)}, {\bf 0})), \dots, ({\bf 0},
(c_{0}^{(n_{1})}, {\bf 0})))\in K_{{\frak H}_{1}}(n_{1}),\\
\hat{Q}_{2}&=&(\eta_{1}, \dots, \eta_{n_{1}}; {\bf 0}, ({\bf 0},
(d_{0}^{(1)}, {\bf 0})), \dots, ({\bf 0},
(d_{0}^{(n_{2})}, {\bf 0})))\in K_{{\frak H}_{1}}(n_{2}).
\end{eqnarray*}

If we can find such positive numbers such that in addition to the
property above,
$s_{0}=r_{i}^{-1}$ and 
\begin{eqnarray*}
\lefteqn{r_{i}^{L(0)}u^{(i)}s_{0}^{L(0)}e^{-L^{-}(B^{(0)})}
(b_{0}^{(0)})^{-L(0)}
v^{(0)}}\nno\\
&&=r_{i}^{L(0)}u^{(i)}r_{i}^{-L(0)}e^{-L^{-}(B^{(0)})}
(b_{0}^{(0)})^{-L(0)}
v^{(0)}\in \hom(T_{\wedge}, 
\overline{T}_{\wedge})
\end{eqnarray*}
is in $\hom(H^{T_{\wedge}}, H^{T_{\wedge}})$, then we say that 
{\it $(Q_{1}, \alpha_{1})\circ_{i}(Q_{2}, \alpha_{2})$ exists}. It is
clear that when $(Q_{1}, \alpha_{1})\circ_{i}(Q_{2}, \alpha_{2})$
exists, $Q_{1}{}_{^{i}}\infty_{^{0}}Q_{2}$ exists.

We assume now that $(Q_{1}, \alpha_{1})\circ_{i}(Q_{2}, \alpha_{2})$
exists. 
Let 
$$\nu_{n}: K(n)\to \hom(T_{\wedge}^{\otimes n},
\overline{T}_{\wedge})$$
for $n\ge 0$
be the maps defining the geometric vertex operator algebra
structure of central charge $0$
on $T_{\wedge}$.
Consider the expression
\begin{eqnarray}\label{4.1}
\lefteqn{(e^{-L^{-}(A^{(0)}(Q_{1}))}u^{(0)}e^{L^{-}(A^{(0)}(Q_{1}))})
\circ \nu_{n_{1}}(Q_{1})}\nno\\
&&\quad\quad \circ (u^{(1)}\otimes 
\cdots\otimes (u^{(i)}\circ 
((e^{-L^{-}(A^{(0)}(Q_{2}))}v^{(0)}e^{L^{-}(A^{(0)}(Q_{2}))})\nno\\
&&\quad\quad
\circ \nu_{n_{2}}(Q_{2})\circ (v^{(1)}\otimes \cdots\otimes
v^{(n_{2})}))
\otimes \cdots \otimes u^{(n_{1})}))\nno\\
&&=(r_{0}^{L(0)}e^{-L^{-}(A^{(0)})}(a_{0}^{(0)})^{-L(0)}
u^{(0)})
\circ \nu_{n_{1}}(\hat{Q}_{1})
\circ ((e^{-L^{+}(A^{(1)})}(a_{0}^{(1)})^{-L(0)}
r_{1}^{L(0)}u^{(1)})\nno\\
&&\quad\quad\otimes \cdots \otimes ((e^{-L^{+}(A^{(i)})}(a_{0}^{(i)})^{-L(0)}
r_{i}^{L(0)}u^{(i)})
\circ ((s_{0}^{L(0)}e^{-L^{-}(B^{(0)})}(b_{0}^{(0)})^{-L(0)}
v^{(0)})\nno\\
&&\quad\quad \circ \nu_{n_{2}}(\hat{Q}_{2})\circ
((e^{-L^{+}(B^{(1)})}(b_{0}^{(1)})^{-L(0)}
s_{1}^{L(0)}v^{(1)})\otimes \cdots\nno\\
&&\quad\quad \otimes (e^{-L^{+}(B^{(n_{2})})}(b_{0}^{(n_{2})})^{-L(0)}
s_{n_{2}}^{L(0)}v^{(n_{2})})))\otimes
\cdots \nno\\
&&\quad\quad \otimes (e^{-L^{+}(A^{(n_{1})})}(a_{0}^{(n_{1})})^{-L(0)}
r_{n_{1}}^{-L(0)}u^{(n_{1})})))
\end{eqnarray}
where $\circ$ means the composition.

We first show that (\ref{4.1}) makes sense and can in fact be viewed as an
element of 
$$\hom(T_{\wedge}^{\otimes (n_{1}+n_{2}-1)}, 
\overline{T}_{\wedge}).$$     
Since $\hat{Q}_{1}\in K_{{\frak H}_{1}}(n_{1})$, 
$\hat{Q}_{2}\in K_{{\frak H}_{1}}(n_{2})$ and
$({\bf 0}, (a^{(i)}_{0}, A^{(i)})), ((B^{(0)}, (b_{0}^{(0)}, {\bf 0}))
\in K_{{\frak H}_{1}}(1)$,
$\nu_{n_{1}}(\hat{Q}_{1})$, $\nu_{n_{2}}(\hat{Q}_{2})$
and $e^{-L^{+}(A^{(i)}}(a_{0}^{(i)})^{-L(0)}$
are in
$\hom((H^{T_{\wedge}})^{\otimes n_{1}},
H^{T_{\wedge}})$, $\hom((H^{T_{\wedge}})^{\otimes n_{2}},
H^{T_{\wedge}})$ and $\hom(H^{T_{\wedge}},
H^{T_{\wedge}})$, respectively. By assumption, we also know that
$$r_{i}^{L(0)}u^{(i)}s_{0}^{L(0)}e^{-L^{-}(B^{(0)})}
(b_{0}^{(0)})^{-L(0)}
v^{(0)}$$ 
is in $\hom(H^{T_{\wedge}},
H^{T_{\wedge}})$.
Thus 
\begin{eqnarray*}
\lefteqn{
\nu_{n_{1}}(\hat{Q}_{1})\circ (1\otimes \cdots}\nno\\
&&\quad\quad
\otimes ((e^{-L^{+}(A^{(i)}}(a_{0}^{(i)})^{-L(0)}
r_{i}^{L(0)}u^{(i)}
s_{0}^{-L(0)}e^{-L^{-}(B^{(0)}}(b_{0}^{(0)})^{-L(0)}
v^{(0)})\nno\\
&&\quad\quad
\circ \nu_{n_{2}}(\hat{Q}_{2})))\otimes
\cdots \otimes 1)
\end{eqnarray*}
is in fact in
$\hom((H^{T_{\wedge}})^{\otimes (n_{1}+n_{2}-1)},
H^{T_{\wedge}})$. Since 
$r_{0}^{-L(0)}e^{-L^{-}(A^{(0)})}(a_{0}^{(0)})^{-L(0)}
u^{(0)}$ can be viewed as a map in $\hom(\overline{T}_{\wedge},
\overline{T}_{\wedge})$ and since 
$e^{-L^{+}(A^{(k)})}(a_{0}^{(k)})^{-L(0)}
r_{k}^{L(0)}u^{(k)}$, $k=1, \dots, n_{1}$, and
$e^{-L^{+}(B^{(l)})}(b_{0}^{(l)})^{-L(0)}
s_{l}^{L(0)}v^{(l)}$, $l=1, \dots, n_{2}$, are
maps in $\hom(T_{\wedge}, T_{\wedge})$, 
(\ref{4.1}) can be viewed as an element of 
$\hom(T_{\wedge}^{\otimes (n_{1}+n_{2}-1)}, 
\overline{T}_{\wedge}).$ 

Next we rewrite (\ref{4.1}) in a different form. Since
both 
$e^{-L^{+}(A^{(i)})}(a_{0}^{(i)})^{-L(0)}$
and
$$r_{i}^{L(0)}u^{(i)}s_{0}^{-L(0)}e^{-L^{-}(B^{(0)})}(b_{0}^{(0)})^{-L(0)}
v^{(0)}$$
are in $\hom(H^{T_{\wedge}},
H^{T_{\wedge}})$, their product is also in $\hom(H^{T_{\wedge}},
H^{T_{\wedge}})$. But this product can be rewritten as
\begin{eqnarray}\label{4.1.1}
\lefteqn{e^{-L^{+}(A^{(i)})}(a_{0}^{(i)})^{-L(0)}
e^{-L^{-}(B^{(0)})}(b_{0}^{(0)})^{-L(0)}\cdot}\nno\\
&&\quad\quad\quad \cdot(b_{0}^{(0)})^{L(0)}
e^{L^{-}(B^{(0)})}r_{i}^{L(0)}u^{(i)}s_{0}^{-L(0)}
e^{-L^{-}(B^{(0)})}(b_{0}^{(0)})^{-L(0)}
v^{(0)}.\quad\quad
\end{eqnarray}
Since as elements of $\hom(T_{\wedge},
\overline{T}_{\wedge})$,
$$e^{-L^{+}(A^{(i)})}(a_{0}^{(i)})^{-L(0)}=\nu_{1}(({\bf 0}, (a_{0}^{(i)},
A^{(i)})))$$ and
$$e^{-L^{-}(B^{(0)})}(b_{0}^{(0)})^{-L(0)}=\nu_{1}((B^{(0)}, 
(b_{0}^{(0)}, {\bf 0})))$$
and since 
$$\nu_{1}(({\bf 0}, (a_{0}^{(i)},
A^{(i)}))), \nu_{1}((B^{(0)}, 
(b_{0}^{(0)}, {\bf 0})))\in
\hom(H^{T_{\wedge}},
H^{T_{\wedge}}),$$  
we have
\begin{eqnarray}\label{4.1.2}
\lefteqn{e^{-L^{+}(A^{(i)})}(a_{0}^{(i)})^{-L(0)}
e^{-L^{-}(B^{(0)})}(b_{0}^{(0)})^{-L(0)}}\nno\\
&&=\nu_{1}(({\bf 0}, (a_{0}^{(i)},
A^{(i)}))) \nu_{1}((B^{(0)}, 
(b_{0}^{(0)}, {\bf 0}))).
\end{eqnarray}

But 
$$\nu_{1}(({\bf 0}, (a_{0}^{(i)},
A^{(i)})))\nu_{1}((B^{(0)}, 
(b_{0}^{(0)}, {\bf 0})))$$
 is in fact in $\hom(H^{T_{\wedge}},
H^{T_{\wedge}})$. By the sewing axiom for geometric vertex operator
algebras (see Chapter 5 of \cite{H}), the right-hand side of (\ref{4.1.2})
is equal to 
\begin{equation}\label{4.1.3}
\nu_{1}(({\bf 0}, (a_{0}^{(i)},
A^{(i)}))_{^{1}}\infty_{^{0}}(B^{(0)}, 
(b_{0}^{(0)}, {\bf 0}))).
\end{equation}

Let $(E^{(0)}, (e_{0}^{(0)}, {\bf 0}))$
and $({\bf 0},
(e_{0}^{(1)}, E^{(1)}))$ be elements of $K_{{\frak H}_{1}}$
such that
$$(E^{(0)}, (e_{0}^{(1)}, {\bf 0}))_{^{1}}\infty_{^{0}}({\bf 0},
(e_{0}^{(1)}, E^{(1)}))=({\bf 0}, (a_{0}^{(i)},
A^{(i)}))_{^{1}}\infty_{^{0}}(B^{(0)}, 
(b_{0}^{(0)}, {\bf 0})).$$
Then (\ref{4.1.3}) is equal to
\begin{eqnarray}\label{4.1.4}
\lefteqn{\nu_{1}((E^{(0)}, (e_{0}^{(0)}, {\bf 0}))_{^{1}}
\infty_{^{0}}({\bf 0},
(e_{0}^{(1)}, E^{(1)}))}\nno\\
&&=e^{-L^{-}(E^{(0)})}(e_{0}^{(0)})^{-L(0)}
e^{-L^{+}(E^{(1)})}(e_{0}^{(1)})^{-L(0)}.
\end{eqnarray}
Since 
$$(E^{(0)}, (e_{0}^{(0)}, {\bf 0})), ({\bf 0},
(e_{0}^{(1)}, E^{(1)}))\in K_{{\frak H}_{1}},$$
we have
$$e^{-L^{-}(E^{(0)})}(e_{0}^{(0)})^{-L(0)},
e^{-L^{+}(E^{(1)})}(e_{0}^{(1)})^{-L(0)}\in \hom(H^{T_{\wedge}},
H^{T_{\wedge}}).$$

By the calculations (\ref{4.1.2})--(\ref{4.1.4}), (\ref{4.1.1}) 
is equal to
\begin{eqnarray}\label{4.1.4.5}
\lefteqn{e^{-L^{-}(E^{(0)})}(e_{0}^{(0)})^{-L(0)}
e^{-L^{+}(E^{(1)})}(e_{0}^{(1)})^{-L(0)}\cdot}\nno\\
&&\quad\quad\cdot
e^{L^{-}(B^{(0)})}r_{i}^{L(0)}u^{(i)}s_{0}^{-L(0)}
e^{-L^{-}(B^{(0)})}(b_{0}^{(0)})^{-L(0)}
v^{(0)}\nno\\
&&=e^{-L^{-}(E^{(0)})}(e_{0}^{(0)})^{-L(0)}
\biggl(e^{-L^{+}(E^{(1)})}(e_{0}^{(1)})^{-L(0)}
e^{L^{-}(B^{(0)})}r_{i}^{L(0)}u^{(i)}s_{0}^{-L(0)}\cdot\nno\\
&&\quad\quad\cdot
e^{-L^{-}(B^{(0)})}(b_{0}^{(0)})^{-L(0)}
v^{(0)}e^{-L^{+}(E^{(1)})}(e_{0}^{(1)})^{-L(0)}\biggr)
(e_{0}^{(1)})^{L(0)}e^{L^{+}(E^{(1)})}.\nno\\
&&
\end{eqnarray}
Since 
$$e^{-L^{-}(E^{(0)})}(e_{0}^{(0)})^{-L(0)},
e^{-L^{+}(E^{(1)})}(e_{0}^{(1)})^{-L(0)}\in 
\hom(H^{T_{\wedge}},
H^{T_{\wedge}})$$ 
and (\ref{4.1.1})
are also in $\hom(H^{T_{\wedge}},
H^{T_{\wedge}})$ and since 
both $e^{-L^{-}(E^{(0)})}(e_{0}^{(0)})^{-L(0)}$,  
$e^{-L^{+}(E^{(1)})}(e_{0}^{(1)})^{-L(0)}$ are invertible, 
\begin{eqnarray}\label{4.1.5}
\lefteqn{e^{-L^{+}(E^{(1)})}(e_{0}^{(1)})^{-L(0)}
e^{L^{-}(B^{(0)})}r_{i}^{L(0)}u^{(i)}s_{0}^{-L(0)}}\nno\\
&&\quad\quad\quad \cdot
e^{-L^{-}(B^{(0)})}(b_{0}^{(0)})^{-L(0)}
v^{(0)}e^{-L^{+}(E^{(1)})}(e_{0}^{(1)})^{-L(0)}
\end{eqnarray}
is also in $\hom(H^{T_{\wedge}},
H^{T_{\wedge}})$. Moreover, we can in fact express (\ref{4.1.5}) 
as a finite linear combination of elements of 
the form $w^{-}w^{+}$ where $w^{-}$ and $w^{+}$ are 
elements of $\widetilde{T}_{\wedge}$ and
$\widetilde{T}_{\wedge}^{\ge -1}$, respectively; here we view
$\widetilde{T}_{\wedge}$ and
$\widetilde{T}_{\wedge}^{\ge -1}$ as subspaces of
$\hom(T_{\wedge},
\overline{T}_{\wedge})$. Thus by the calculations 
(\ref{4.1.2})--(\ref{4.1.4.5}), 
(\ref{4.1.1}) is a finite linear combination of 
elements of the form 
$$e^{-L^{-}(E^{(0)})}(e_{0}^{(0)})^{-L(0)}
w^{-}w^{+}
(e_{0}^{(1)})^{L(0)}e^{L^{+}(E^{(1)})}.$$

Using the calculations and discussions above, the Jacobi identity 
and the study of exponentials of infinite sums of Virasoro operators
in \cite{H}, 
we see that (\ref{4.1}) 
can  be written as 
a finite linear combination
of elements of 
$\hom(T_{\wedge}^{\otimes (n_{1}+n_{2}-1)}, 
\overline{T}_{\wedge})$
of the form
\begin{eqnarray}\label{4.2}
\lefteqn{(e^{-L^{-}(A^{(0)}(Q_{1}\;_{^{i}}\infty_{^{0}}Q_{2}))}w^{(0)}
e^{-L^{-}(A^{(0)}(Q_{1}\;_{^{i}}\infty_{^{0}}Q_{2}))})\circ}\nno\\
&&\quad\quad\quad
\nu_{n_{1}+n_{2}-1}(Q_{1}\;_{^{i}}\infty_{^{0}}Q_{2})\circ
(w^{(1)} \otimes \cdots\otimes 
w^{(n_{1}+n_{2}-1)}),
\end{eqnarray}
where 
$w^{(0)}$ is an element of $\hom(T_{\wedge}, \overline{T}_{\wedge})$ 
induced {from} an element (denoted by the same notations
$w^{(0)}$) of  $\widetilde{T}_{\wedge}$, and 
$w^{(k)}$, $k=1, \dots, n_{1}+n_{2}-1$, are 
elements of $\hom(T_{\wedge}, \overline{T}_{\wedge})$ 
induced {from} elements (denoted by the same notations
$w^{(k)}$) 
of $\widetilde{T}_{\wedge}^{\ge -1}$. 

The element (\ref{4.2}) gives an element
\begin{equation}\label{4.3}
{\frak K}_{n_{1}+n_{2}-1}^{*}(w^{(0)}\otimes w^{(1)}\otimes \cdots \otimes 
w^{(n_{1}+n_{2}-1)})
\end{equation}
of $\widetilde{\cal T}_{\wedge}(n_{1}+n_{2}-1)$.  Taking the same
linear combination as what we obtained for (\ref{4.1}) but replacing
maps of the form (\ref{4.2}) by the corresponding elements (\ref{4.3})
of $\widetilde{\cal T}_{\wedge}(n_{1}+n_{2}-1)$, we obtain an element
of $\widetilde{\cal T}_{\wedge}(n_{1}+n_{2}-1)$, that is, a section of
${\frak G}(n_{1}+n_{2}-1)$.  Let
$\alpha_{1}{}_{^{i}}*_{^{0}}\alpha_{2}$ be the value of this section
at $Q_{1}\;_{^{i}}\infty_{^{0}}Q_{2}$. It is clear {from} the
construction above that $\alpha_{1}{}_{^{i}}*_{^{0}}\alpha_{2}$ is
uniquely determined by $(Q_{1}, \alpha_{1})$ and $(Q_{2},
\alpha_{2})$.  We define
$$(Q_{1}, \alpha_{1})\circ_{i}(Q_{2}, \alpha_{2})=
(Q_{1}\;_{^{i}}\infty_{^{0}}Q_{2}, \alpha_{1}{}_{^{i}}*_{^{0}}\alpha_{2}).$$
This finishes the definition of the composition maps.

The bundle ${\frak G}(n)$ has a constant section $1$. 
When $n=1$, the value of $1$ at the identity $I$ 
of $K(1)$ gives an element $I_{\frak G}=(I, 1\mbar_{I})$.
For $n\ge 0$, there is an obvious action of $S_{n}$ on ${\frak G}(n)$. 
{From} the definition, the following result is obvious:

\begin{propo}
The sequence ${\frak G}$ equipped with the composition maps defined above,
the identity $I_{\frak G}$ and the obvious actions of $S_{n}$ is an analytic
partial operad.\epf
\end{propo}

We consider the vertex operator subalgebra $T_{G}$ of 
$T_{\wedge}$
generated by $b$, $\omega_{G}$, $f_{G}$ and $q_{G}$.
We consider the elements of 
$\widetilde{\cal T}_{\wedge}(1)$  
of the form 
${\frak K}_{n_{1}}^{*}(\sum_{j=1}^{k}
u^{(0)}_{j}\otimes u_{j}^{(1)})$, where  
$u^{(0)}_{j}\in \widetilde{T}_{\wedge}\cap \hom(T_{G}, \overline{T}_{G})$
and $u_{j}^{(1)}\in \widetilde{T}_{\wedge}^{\ge -1}
\cap \hom(T_{G}, \overline{T}_{G})$, $j=1, \dots, k$.
These elements give a holomorphic subbundle ${\cal T}_{G}(1)$
of ${\frak G}(1)$ and $T_{G}$ gives a holomorphic 
subbundle ${\cal T}_{G}(0)$
of ${\frak G}(0)$. The value at $P(1)$ 
of the section $1$ of 
${\frak G}(2)$ above gives an element $(P(1), 1\mbar_{P(1)})$
of ${\frak G}(2)$. 
Let ${\cal T}_{G}$ be the partial suboperad of 
${\frak G}$ generated by ${\cal T}_{G}(0)$, ${\cal T}_{G}(1)$ and 
$(P(1), 1\mbar_{P(1)})$. This partial operad is still not
the (strong) topological partial operad we need but it
motivates our
construction of the (strong) topological partial operad in Subsection 4.2.

\renewcommand{\theequation}{\thesection.\arabic{equation}}
\renewcommand{\therema}{\thesection.\arabic{rema}}
\setcounter{equation}{0}
\setcounter{rema}{0}
\section{(Strong) topological vertex partial operads and the
proofs of the main results}

\subsection{Properties of 
topological vertex operator algebras}

In this subsection, we prove some results on 
topological vertex operator algebras
needed in the construction of (strong) topological partial
operads and the proofs of the main theorem. 

First we need several lemmas:

\begin{lemma}\label{lm1}
Let $V$ be a ${\Bbb Z}$-graded vector space. 
Consider the space 
$\mbox{\rm End}\;V$ spanned by 
homogeneous operators on $V$. Then $\mbox{\rm End}\;V$ is a
${\Bbb Z}$-graded associative algebra. If we use $|A|$ to denote 
the grading of a homogeneous element $A$
of $\mbox{\rm End}\;V$ and define a bracket 
$[\cdot, \cdot]$ on $\mbox{\rm End}\;V$ by
$$[A, B]=AB-(-1)^{|A||B|}BA,$$
then $\mbox{\rm End}\;V$ equipped with the graded associative 
algebra structure and the bracket $[\cdot, \cdot]$ is a
graded Poisson algebra. In other words, in addition to the axioms
for a graded associative algebra, the following conditions are satisfied:
\begin{eqnarray}
{[A, B]}&=&(-1)^{|A||B|}[B, A]\label{skew-sym}\\
{[A, [B, C]]}&=&[[A, B], C] + 
(-1)^{|A| |B|}[B, [A, C]]\label{jacobi}\\
{[A, BC]}&=&[A, B]C +(-1)^{|A| |B|} B[A, C].\label{derivation}\epfd
\end{eqnarray}
\end{lemma}

\begin{lemma}\label{lm3}
Let $V$ be a topological vertex operator algebra.
Then for $i, j\in {\Bbb Z}$, 
\begin{eqnarray}\label{g-l}
{[g(i), L(j)]}&=&[L(i), g(j)]\nno\\
&=&(i-j)g(i+j).
\end{eqnarray}
\end{lemma}
\pf 
Since $L(n)g=0$ for $n>0$ and 
$L(0)g=2g$, the vertex operator $Y(g, x)$ of $g$ is a primary field 
of weight $2$. So for $i, j\in {\Bbb Z}$, we have
\begin{eqnarray}\label{lm3-1}
{[L(i), g(j)]}&=&[L(i), \res_{x}x^{j+1}Y(g, x)]\nno\\
&=&\res_{x}x^{j+1}[L(i), Y(g, x)]\nno\\
&=&\res_{x}x^{j+1}\left(x^{i+1}\frac{d}{dx}+2(i+1)x^{i}\right)Y(g, x)\nno\\
&=&(i-j)g(i+j).
\end{eqnarray}
By (\ref{skew-sym}) and (\ref{lm3-1}), we obtain
\begin{eqnarray*}
{[g(i), L(j)]}&=&-[L(j), g(i)]\nno\\
&=&-(j-i)g(j+i)\nno\\
&=&(i-j)g(i+j).\epfd
\end{eqnarray*}

\begin{propo}\label{lm2}
Let $V$ be a topological vertex operator algebra.
Then for any $n \in {\Bbb Z}$ and 
$v\in V$, 
\begin{equation}\label{qq-brkt}
[Q, [Q, v_{n}]]=0,
\end{equation}
where $v_{n}=\res_{x}x^{n}Y(v, x)$.
\end{propo}
\pf
For any $u\in V$, by the commutator formula for 
 $V$,
$$[q_{0}, Y(u, x)]=Y(q_{0}u, x).$$
Thus 
\begin{eqnarray*}
{[Q, [Q, v_{n}]]}&=&[q_{0}, [q_{0}, v_{n}]]\nno\\
&=&[q_{0}, [q_{0}, \res_{x}x^{n}Y(v, x)]]\nno\\
&=&\res_{x}x^{n}[q_{0}, [q_{0}, Y(v, x)]]\nno\\
&=&\res_{x}x^{n}[q_{0}, Y(q_{0}v, x)]\nno\\
&=&\res_{x}x^{n}Y(q^{2}_{0}v, x)\nno\\
&=&\res_{x}x^{n}Y(Q^{2}v, x)\nno\\
&=&0. \epfd
\end{eqnarray*}

\begin{propo}\label{5-1}
The central charge of a topological vertex operator algebra $V$ is $0$.
\end{propo}
\pf 
Let the central charge of $V$
be $c$. Then  for any $i\in {\Bbb Z}$,
\begin{equation}\label{6.1.5}
[L(i), L(-i)]=2iL(0)+\frac{c}{12}(i^{3}-i).
\end{equation}
On the other hand, using $L(-i)=[Q, g(-i)]$ and 
the Jacobi identity (\ref{jacobi}),
\begin{eqnarray}\label{6.2}
[L(i), L(-i)]&=&[L(i), [Q, g(-i)]]\nno\\
&=&[[L(i), Q], g(-i)]
+[Q, [L(i), g(-i)]],
\end{eqnarray}
By (\ref{g-l}), 
$[L(i), g(-i)]=2ig(0).$ Thus 
$$[Q, [L(i), g(-i)]]=[Q, 2ig(0)]=2iL(0).$$
So we see that the right-hand side of (\ref{6.2}) 
is equal to $2iL(0)$. Comparing it with the right-hand side of 
(\ref{6.1.5}), we obtain $c=0$. \epfv

\begin{propo}
Let $V$ be a topological vertex operator algebra, then 

\begin{enumerate}

\item For any $i, j \in \Bbb Z$, 
$[Q, [g(i), g(j)]]=0$. 

\item If $g(0)^2=0$, then we have $[g(i), g(j)]=0$ for any 
$i, j \in \Bbb Z$. In particular,
we have $g(i)^2=0$ for any $i \in {\Bbb Z}$.

\end{enumerate}
\end{propo}
\pf
The first conclusion follows {from} the following straightforward  
calculation using $[Q, g(i)]=L(i)$, $i\in {\Bbb Z}$, 
Lemmas \ref{lm1} and \ref{lm3}:
\begin{eqnarray}
[Q, [g(i), g(j)]]&=&[[Q, g(i)], g(j)] -[g(i), [Q, g(j)]]\nno\\
&=& [L(i), g(j)]-[g(i), L(j)]\nno\\
&=&(i-j) g(i+j)-(i-j) g(i+j)\nno\\
&=& 0.
\end{eqnarray}

To prove the second conclusion, first we claim that 
$[g(i), g(0)]=0$ for any $i \in \Bbb Z$. For $i=0$, this
follows  trivially  
{from} the condition $g(0)^2=0$. So we assume $i \neq 0$. Then 
by (\ref{derivation}) and (\ref{g-l}),
\begin{eqnarray}
0&=&[L(i), g(0)^2]\nno\\
&=&[L(i), g(0)]g(0) + g(0)[L(i), g(0)]\nno\\
&=& i (g(i)g(0)+g(0)g(i) )\nno\\
&=& i[g(i), g(0)]
\end{eqnarray}
So we also have $[g(0), g(i)]=0$ for  $i \ne 0$. 
Finally we assume $i, j \neq 0$. Using 
$[g(i), g(0)]=0$ for $i\in {\Bbb Z}$ which we have just proved, the formulas
(\ref{g-l}) and (\ref{jacobi}), we have
\begin{eqnarray}
[g(i), g(j)]&=&\frac 1{j} [g(i), [L(j), g(0)]]\nno\\
&=& \frac 1{j} [[g(i), L(j)], g(0)] +\frac 1{j} [ L(j), [g(i), g(0)]] \nno\\
&=& \frac {i-j} {j} [g(i+j), g(0)] +\frac 1{j} [ L(j), [g(i), g(0)]] \nno\\
&=& 0. \epfd
\end{eqnarray}

By this proposition, for any strong topological 
vertex operator algebra, the operators
$g(i)$, $i\in {\Bbb Z}$, are anti-commutative.

\subsection{Construction of 
(strong) topological vertex partial operads}

We construct topological vertex partial operads ${\cal T}^{k}$ of type
$k$ and strong topological partial operad $\bar{\cal T}^{k}$ of type
$k$ for $k<0$
in this subsection.

To construct (strong) topological vertex partial operads, we need 
to construct a grading-restricted (strong) 
topological vertex operator algebra $T^{k}$ ($\bar{T}^{k}$) of type $k$
for each
$k<0$ such that it is
{\it universal} in the sense that given any 
locally-grading-restricted (strong)   topological 
vertex operator algebra $V$ of type $k$, 
there is a unique homomorphism 
{from} $T^{k}$ ($\bar{T}^{k}$) to $V$.

We first construct a ${\Bbb Z}\times {\Bbb Z}$-graded
vertex algebra $U$. Consider the set
$$S=\{\omega_{G}, b, f_{G}, q_{G}\}\subset T_{G}.$$ 
Vertex operators of elements of this set have the component operators 
$L(n)$, $b(n), f_{G}(n), q_{G}(n)$, $n\in {\Bbb Z}$. Clearly, these
operators are linearly independent. Let $U_{S}$ be the 
vector space spanned by the operators $L(n), b(n)$ and 
$f_{G}(n), q_{G}(n)$ and 
let $T(U_{S})$ be the tensor algebra 
over $U_{S}$. Then as a vector space, $T(U_{S})$ has 
a basis consisting $1$ and elements of the form
$u_{1}(n_{1})\otimes \dots \otimes u_{l}(n_{l})$
for $l\ge 0$, $u_{1}, \dots, u_{l}\in S$ and 
$n_{1}, \dots, n_{l}\in {\Bbb Z}$. Let $J$ be the ideal of $T(U_{S})$
generated by $L(n)$, $b(n)$, $n\ge -1$, and $f_{G}(n)$, $q_{G}(n)$,
$n\ge 0$, and let $U=T(U_{S})/J$. We shall use 
$[u_{1}(n_{1})\otimes \dots \otimes u_{l}(n_{l})]$ to
denote the coset in $U$ containing 
$u_{1}(n_{1})\otimes \dots \otimes u_{l}(n_{l})$.
We shall denote $[\omega_{G}(-2)]$,
$[b(-2)]$, $[f_{G}(-1)]$ and
$[q_{G}(-1)]$ in $U$ by 
$\omega_{U}$, $g_{U}$, $f_{U}$ 
and $q_{U}$, respectively. 

We define the weight of the element 
$[u_{1}(n_{1})\otimes \cdots \otimes 
u_{l}(n_{l})]$ to be $-n_{1}-\cdots -n_{l}$.
We define the fermion number of the same element to be 
the number of $b$'s in $(u_{1}, \dots, u_{l})$ 
minus  the number of 
 $q_{G}$'s in $(u_{1}, \dots, u_{l})$. 
We use the notation $|v|$ to denote
the fermion number of an element $v\in U$.
Thus $U$ becomes a ${\Bbb Z}\times {\Bbb Z}$-graded vector space.

For any $u\in S$ and $n\in {\Bbb Z}$, we define an action of $u(n)$ on
$U$ as follows: We define 
$u(n)[1]=[u(n)]$ and 
$$u(n)([u_{1}(n_{1})\otimes \dots \otimes 
u_{l}(n_{l})])=[u(n)\otimes u_{1}(n_{1})\otimes \dots \otimes 
u_{l}(n_{l})]$$
for elements of $U$ of the form 
$[u_{1}(n_{1})\otimes \dots \otimes 
u_{l}(n_{l})]$, $l>0$.
We have
$$[u_{1}(n_{1})\otimes \dots \otimes 
u_{l}(n_{l})]=u_{1}(n_{1}) \dots 
u_{l}(n_{l}) [1].$$

We define a vertex operator map $Y: U
\to (\mbox{\rm End}\;U)[[x, x^{-1}]]$
as follows: We define $Y([1], x)=I_{U}$, 
the identity operator on $U$.
We also define 
\begin{eqnarray*}
Y(\omega_{U}, x)&=&\sum_{n\in {\Bbb Z}}L(n)x^{-n-2}\\
Y(g_{U}, x)&=&\sum_{n\in {\Bbb Z}}b(n)x^{-n-2}\\
Y(f_{U}, x)&=&\sum_{n\in {\Bbb Z}}f_{G}(n)x^{-n-1}\\
Y(q_{U}, x)&=&\sum_{n\in {\Bbb Z}}q_{G}(n)x^{-n-1}.
\end{eqnarray*}
Let  $v=u(m)w$.
We define
\begin{eqnarray*}
Y(v, x)&=&\res_{x_{1}}(x_{1}-x)^{m+1}
Y(u, x_{1})Y(w, x)\nno\\
&&-(-1)^{|u||w|}\res_{x_{1}}
(-x+x_{1})^{m+1}Y(w, x)Y(u, x_{1})
\end{eqnarray*}
when $u=\omega_{U}, g_{U}$, and define
\begin{eqnarray*}
Y(v, x)&=&\res_{x_{1}}(x_{1}-x)^{m}
Y(u, x_{1})Y(w, x) \nno\\
&&-(-1)^{|u||w|}\res_{x_{1}}
(-x+x_{1})^{m}Y(w, x)Y(u, x_{1})
\end{eqnarray*}
when $u=f_{U}, q_{U}$. 
By recursion, we obtain a vertex operator map 
$Y: U\to (\mbox{\rm End}\;U)[[x, x^{-1}]]$.

\begin{propo}\label{tilde-t}
The ${\Bbb Z}\times {\Bbb Z}$-graded vector space $U$ equipped
with the vertex operator map $Y$ defined above, the vacuum $[1]$
and the Virasoro element $\omega_{U}$ 
is a ${\Bbb Z}\times {\Bbb Z}$-graded
vertex operator algebra without grading restrictions. 
Given any locally-grading-restricted topological  vertex operator algebra 
$V$, there exists a unique homomorphism of 
${\Bbb Z}\times {\Bbb Z}$-graded vertex operator algebra without grading 
restrictions
{from} $U$ to $V$ such that 
the images of $g_{U}$, $f_{U}$ and $q_{U}$ are
$g$, $f$ and $q$, respectively.
\end{propo}
\pf
The proof is the same as
the proof that $M(c)$ is a vertex operator algebra in \cite{H1}.
\epfv

The ${\Bbb Z}\times {\Bbb Z}$-graded vertex  algebra $U$ is not a topological
vertex operator algebra because  the additional axioms for $g_{U}$, $f_{U}$ and 
$q_{U}$ are not
satisfied. For $k<0$, the grading-restricted 
topological vertex operator algebra $T^{k}$ 
of type $k$ is
constructed {from} $U$ as follows: Let $R^{k}$ be the ideal
of $U$ generated by elements of the forms
$[u_{1}(n_{1})\otimes \cdots\otimes u_{l}(n_{l})]$ for 
$l>0$, $u_{1}, \dots, u_{l}\in S$, $n_{1}, \dots, n_{l}\in {\Bbb Z}$
satisfying $-n_{1}-\cdots -n_{l}<k$,
$f_{U}(0)v-|v|v$, $L(n)q_{U}$, 
$q_{U}(0)^{2}v$, $L(n)g_{U}$,
and $q_{U}(0)g_{U}-\omega_{U}$ for $v\in U$,
$n>0$.  
Let
$T^{k}=U/R^{k}$.  Let
$g_{T^{k}}$, $f_{T^{k}}$ and $q_{T^{k}}$ be the cosets in $T^{k}$ containing
$g_{U}$, $f_{U}$ and $q_{U}$,
respectively. The following result is obvious {from} Proposition \ref{tilde-t}
and the construction of
$T^{k}$:

\begin{propo}\label{univ-t}
For any $k<0$, 
$T^{k}$ together with $g_{T^{k}}$, $f_{T^{k}}$ and $q_{T^{k}}$
is a grading-restricted 
topological vertex operator algebra of type $k$. It satisfies the 
following universal property:
Given any locally-grading-restricted topological  vertex operator algebra 
$V$ of type $k$, there exists a unique homomorphism of 
locally-grading-restricted topological  vertex operator
algebras {from} $T^{k}$ to $V$.\epf
\end{propo}

For $k<0$, using $T^{k}$, we construct a 
grading-restricted  strong topological vertex
operator algebra $\bar{T}^{k}$ of type $k$
as follows: Let $O^{k}$ be the ideal
of $T^{k}$ generated by elements of the form $g_{T^{k}}^{2}(0)v$ for
$v\in T^{k}$. We define
$\bar{T}=T^{k}/O^{k}$. Let
$g_{\bar{T}^{k}}$, $f_{\bar{T}^{k}}$ and $q_{\bar{T}^{k}}$ 
be the cosets in $T^{k}$ containing
$g_{T^{k}}$, $f_{T^{k}}$ and $q_{T^{k}}$,
respectively. The following result is obvious:

\begin{propo}
For any $k<0$,  $\bar{T}^{k}$ together with $g_{\bar{T}^{k}}$, 
$f_{\bar{T}^{k}}$ and $q_{\bar{T}^{k}}$ 
 is a grading-restricted strong topological 
vertex operator algebra of type $k$. It satisfies the 
following universal property:
Given any  locally-grading-restricted  strong
topological  vertex operator
algebra 
$V$ satisfying of type $k$, there exists a unique homomorphism of 
locally-grading-restricted topological vertex operator
algebras {from} $U$ to $V$.\epf
\end{propo}

Motivated by the partial operad ${\cal T}_{G}$, we now construct the 
topological vertex 
partial operad ${\cal T}^{k}$ of type $k$ for $k<0$
using the topological vertex operator algebra
$T^{k}$ of type $k$ and the sphere partial operad $K$ as follows: 
Since components of vertex operators of elements of $U_{S}\subset T_{G}$ can
be identified with holomorphic sections of the holomorphic bundle ${\cal
T}_{G}(1)$ over $K(1)$, there exists a holomorphic bundle  
over $K(1)$ such that 
 components of vertex operators of elements
of $U$  span a space of holomorphic sections of
this holomorphic bundle over $K(1)$ 
and fibers over $Q\in K(1)$ of this holomorphic bundle are 
spanned by the values of the sections in this space at $Q$. Note
that the space $Y(U)$ of components of vertex operators 
of elements of $U$ is a ${\Bbb Z}$-graded vector space.
Its algebraic completion $\overline{Y(U)}$ 
is a subspace of $\hom(U, \overline{U})$. 
Let $H^{U}$
be the locally convex topological completion of $U$ constructed in  
\cite{H6} and \cite{H7} (see Subsection 3.4) and let 
 $\widetilde{Y(U)}$  
be the subspace of $\overline{Y(U)}$ consisting of elements $u\in 
\overline{Y(U)}$ satisfying the following property: There exists 
a positive number $a$ such that $a^{L(0)}ua^{-L(0)}\in
\hom(U, \overline{U})$ is in fact in $\hom(H^{U}, H^{U})\subset
\hom(U, \overline{U})$. Then there exists a holomorphic vector bundle
$\tilde{\cal U}(1)$ over $K(1)$ such that elements of $\widetilde{Y(U)}$ 
span
a space of holomorphic sections of $\tilde{\cal U}(1)$ and 
fibers of $\tilde{\cal U}(1)$ are 
spanned by the values of the sections in this space. Using 
${\frak E}: K(0)\to K(1)$, we pull $\tilde{\cal U}(1)$ back to a holomorphic 
bundle $\tilde{\cal U}(0)$ over $K(0)$. Also we have 
 a holomorphic bundle
$\tilde{\cal U}(1)\mbar_{K_{\ge -1}(1)}$ over $K_{\ge -1}(1)$. 
It is clear that 
 $\tilde{\cal U}(1)$ is canonically isomorphic to 
the pullback of the exterior product bundle $\tilde{\cal
U}(0)\boxtimes \tilde{\cal U}(1)\mbar_{K_{\ge -1}(1)}$ by
the map ${\frak K}_{1}: K(1)\to K(0)\times
K_{\ge -1}(1)$. For $n\ge
2$, using the canonical injective map ${\frak K}_{n}: K(n)\to K(0)\times
(K_{\ge -1}(1))^{n}$, we pull the exterior product bundle $\tilde{\cal
U}(0)\boxtimes (\tilde{\cal U}(1)\mbar_{K_{\ge -1}(1)})^{\boxtimes n}$
over $K(0)\times (K_{\ge -1}(1))^{n}$ to obtain a holomorphic bundle
$\tilde{\cal U}(n)$ over $K(n)$. 
We obtain a sequence $\tilde{\cal
U}=\{\tilde{\cal U}(n)\}_{n\ge 0}$ of holomorphic vector bundles. 

Let $\widetilde{Y(R)}^{k}$ be the intersection
of $\widetilde{Y(U)}$ and 
the algebraic completion $\overline{R}^{k}$ of $R^{k}$. 
Since $R^{k}$ is a ${\Bbb Z}\times {\Bbb Z}$-graded 
vertex operator subalgebra without grading-restrictions 
of $U$ and since fibers of $\tilde{\cal U}(1)$ are 
spanned by the values of the sections in $\widetilde{Y(U)}$ 
which is viewed as a space of holomorphic sections
of $\tilde{\cal U}(1)$, 
$\widetilde{Y(R)}^{k}$ is isomorphic to a space of holomorphic sections of 
a holomorphic subbundle  $\tilde{\cal R}(1)$ of $\tilde{\cal U}(1)$
such that the fibers of  $\tilde{\cal R}(1)$ 
are spanned by the values of the sections in $\widetilde{Y(R)}^{k}$. 
The pullback of 
$\tilde{\cal R}(1)$ by ${\frak E}$ is a holomorphic subbundle 
$\tilde{\cal R}^{k}(0)$  of $\tilde{\cal U}(0)$.
Using the same construction as the one for $\tilde{\cal U}(n)$, $n\ge 2$,
we construct holomorphic subbundles $\tilde{\cal R}^{k}(n)$, $n\ge 2$, 
of $\tilde{\cal U}(n)$. Let ${\cal T}^{k}(n)
=\tilde{\cal U}(n)/\tilde{\cal R}^{k}(n)$ for $n\ge 0$. 

Let $Y(T^{k})$ be the space of components of vertex operators
of elements of the topological vertex operator algebra
$T^{k}$ and $\overline{Y(T^{k})}$ its 
algebraic completion. Let $H^{T^{k}}$ be the locally convex completion of 
$T^{k}$ and $\widetilde{Y(T^{k})}$  the subspace of $\overline{Y(T^{k})}$ 
consisting of elements $u\in 
\overline{Y(T^{k})}$ satisfying the following property: There exists 
a positive number $a$ such that $a^{L(0)}ua^{-L(0)}\in
\hom(T^{k}, \overline{T}^{k})$ is in fact in $\hom(H^{T^{k}}, H^{T^{k}})\subset
\hom(T^{k}, \overline{T}^{k})$.
Then by definition, 
we see that  $\overline{Y(T^{k})}$
is isomorphic to a generating subspace of holomorphic sections of 
${\cal T}^{k}(1)$. 
The holomorphic bundles ${\cal T}^{k}(n)$ can also be constructed {from}
${\cal T}^{k}(1)$ using the 
same construction as the one for $\tilde{\cal U}(n)$, $n\ge 2$.
That is, we have holomorphic bundles ${\cal T}^{k}(0)$ and 
${\cal T}^{k}(1)|_{K_{\ge -1}(1)}$ over $K(0)$ and $K_{\ge -1}(1)$, 
respectively, such that 
$${\cal T}^{k}(n)={\frak K}_{n}^{*}({\cal T}^{k}(0)\boxtimes
({\cal T}^{k}(1)|_{K_{\ge -1}(1)})^{\boxtimes n}.$$
Thus we see that the sequence ${\cal T}^{k}=\{{\cal T}^{k}(n)\}_{n\ge 0}$
of holomorphic bundles 
can be constructed {from} the topological vertex operator algebra $T^{k}$
of type $k$. Since $T^{k}$ is a grading-restricted 
topological vertex operator algebra,
the same method as the one in 
the construction of the composition maps for the 
partial operad ${\frak G}$ (see Subsection 3.4)
gives composition maps for ${\cal T}^{k}$.
The
symmetric group $S_{n}$ acts on ${\cal T}^{k}(n)$  in the obvious way
for $n\ge 0$. The constant term of the vertex operator of the vacuum of
$T^{k}$ is
a holomorphic section of ${\cal T}^{k}(1)$. We denote the value of
this section at $I\in K(1)$ by $I_{{\cal T}^{k}}$.

Similarly, using the strong topological vertex operator algebra
$\bar{T}^{k}$ of type $k$, we can construct
a sequence $\bar{\cal T}^{k}=\{\bar{\cal T}^{k}(n)\}_{n\ge 0}$ 
of holomorphic bundles together with composition maps, an  identity
$I_{\bar{\cal T}^{k}}$ and actions of $S_{n}$
{from} ${\cal T}^{k}$. 
We have:

\begin{propo}
The sequence ${\cal T}^{k}$ equipped with the composition maps,
the identity $I_{{\cal T}^{k}}$ and the  actions of 
$S_{n}$ induced {from} ${\cal T}^{k}$
is an analytic
partial operad. The sequence $\bar{\cal T}^{k}$ 
equipped with the composition maps,
the identity $I_{\bar{\cal T}^{k}}$ and the  actions of 
$S_{n}$ 
is also an analytic
partial operad.\epf
\end{propo}

For $k<0$, the partial operad ${\cal T}^{k}$
is called the {\it topological vertex partial operad of type $k$} and
the partial operad $\bar{\cal T}^{k}$ is called the 
{\it strong topological vertex partial operad of type $k$}. 
Clearly these partial 
operads are  ${\Bbb C}^{\times}$-rescalable (see \cite{HL1}, \cite{HL2}
and Appendix C of \cite{H} for the definition). By definition, for any 
$k<0$,
there is a  canonical morphism of ${\Bbb C}^{\times}$-rescalable analytic
partial operads  {from} 
${\cal T}^{k}$ to $\bar{\cal T}^{k}$.

We now define the notion of (weakly) meromorphic 
${\Bbb Z}\times {\Bbb Z}$-graded algebra over 
${\cal T}^{k}$ and
$\bar{\cal T}^{k}$, $k<0$. 
First, recall {from} \cite{H2} that 
for a ${\Bbb Z}\times {\Bbb Z}$-graded 
space $V$ (graded by fermion numbers and by weights) and 
a ${\Bbb Z}\times {\Bbb Z}$-graded subspace $W$, we can define the 
graded endomorphism partial pseudo-operad ${\cal H}_{V}$ to be
the sequence $\{\hom(V^{\otimes n}, \overline{V})\}_{n\in \ge 0}$
as in \cite{HL1}, \cite{HL2} and Appendix C of \cite{H},
except that the left actions of the symmetric groups are defined such that 
for any $f\in \hom(V^{\otimes n}, \overline{V})$
and $\sigma_{i,i+1}\in S_{n}$ which is the transposition permuting $i$
and $i+1$,
\begin{equation}
(\sigma_{i, i+1}(f))(v_{1}, \dots, v_{i}, v_{i+1}, \dots, v_{n})
=(-1)^{|v_{i}||v_{i+1}|}f(v_{1}, \dots, v_{i+1}, v_{i}, \dots, v_{n})
\end{equation}
for any $v_{1}, \dots, v_{n}\in V$ with homogeneous fermion numbers.
Using this graded endomorphism partial pseudo-operad, we can define 
the notions of {\it ${\Bbb Z}\times {\Bbb Z}$-graded algebra} 
over any ${\Bbb C}^{\times}$-rescalable
partial operads in the obvious way. 

Let ${\bf 1}_{T^{k}}$ be the vacuum of $T^{k}$ 
and ${\bf 1}_{T^{k}}(0)=\res_{x} x^{-1}Y({\bf 1}_{T^{k}}, x)\in 
\widetilde{Y(T^{k})}$.
Then the holomorphic sections ${\frak
K}_{n}^{*}({\frak E}^{*}({\bf 1}_{T^{k}}(0))\boxtimes 
({\bf 1}_{T^{k}}(0)|_{K_{\ge -1}(1)})^{\boxtimes n})$  
of ${\cal T}^{k}(n)$, 
$n\ge 0$, give
a morphism of partial operads {from} $K$ to ${\cal T}^{k}$. 
Composing with the morphisms of partial operads
{from} 
${\cal T}^{k}$ to $\bar{\cal T}^{k}$, we obtain a morphism of partial operads
{from}  $K$ to $\bar{\cal T}^{k}$.
In particular, 
an ${\Bbb Z}\times {\Bbb Z}$-graded algebra over 
${\cal T}^{k}$ or $\bar{\cal T}^{k}$ 
is also an ${\Bbb Z}\times {\Bbb Z}$-graded algebra over $K$, or equivalently,
an algebra over $\tilde{K}^{0}$. 

Recall the following definition of 
meromorphic algebra over $\tilde{K}^{c}$ 
in \cite{H}: 

\begin{defi}
{\rm A  {\it meromorphic algebra over  
$\tilde{K}^{c}$} is a ${\Bbb Z}$-graded vector space 
$V=\coprod_{n\in {\Bbb Z}}V_{(n)}$ equipped with a 
morphism $\nu$ of partial operads {from} $\tilde{K}^{c}$ to 
${\cal H}_{V}=\{\hom(V^{\otimes n}, \overline{V})$
satisfying the following axioms:

\begin{enumerate}

\item  The grading-restriction axiom:
$\dim V_{(n)}<\infty$ for $n\in {\Bbb Z}$ and 
$V_{(n)}=0$ for $n$ sufficiently small.

\item For any $n\in {\Bbb N}$, $\nu_{n}$ is linear on any fiber
of $\tilde{K}^{c}(n)$.

\item For any positive integer $n$, $v'\in V'$,
$v_{1}, \dots, v_{n}\in V$, the function
$$Q\to \langle v', \nu_{n}(\psi_{n}(Q))) (v_{1}\otimes \cdots
\otimes v_{n})\rangle$$ on $K(n)$ is meromorphic (in the sense of
Section 3.1) and if $z_{i}$ and $z_{j}$ are the $i$-th and $j$-th
punctures of $Q\in K(n)$ respectively, then for any $v_{i}$ and
$v_{j}$ in $V$ there exists a positive integer $N(v_{i}, v_{j})$ such
that for any $v'\in V'$, $v_{k}\in V$, $k\ne i, j$, the order of the
pole $z_{i}=z_{j}$ (we use the convention $z_{n}=0$) of $$\langle v',
\nu_{n}(\psi_{n}(Q)))(v_{1}\otimes \cdots \otimes v_{n})\rangle$$
 is  less than $N(v_{i}, v_{j})$.
\end{enumerate}}
\end{defi}

We also need:

\begin{defi}
{\rm A {\it weakly meromorphic algebra over  
$\tilde{K}^{c}$} is a ${\Bbb Z}$-graded vector space 
$V=\coprod_{n\in {\Bbb Z}}V_{(n)}$ equipped with a 
morphism $\nu$ of partial operads {from} $\tilde{K}^{c}$ to 
${\cal H}_{V}=\{\hom(V^{\otimes n}, \overline{V})$
satisfying all the axioms for meromorphic algebras over  
$\tilde{K}^{c}$ except the grading-restriction axioms.
The notions of {\it (weakly) meromorphic (${\Bbb Z}\times {\Bbb Z}$-graded) 
algebra
over $\tilde{K}$} and {\it (weakly) meromorphic ${\Bbb Z}\times 
{\Bbb Z}$-graded algebra
over $\tilde{K}^{c}$}  are defined similarly.}
\end{defi}

 We have the following notions:

\begin{defi}
{\rm For $k<0$, a ${\Bbb Z}\times {\Bbb Z}$-graded algebra $V$ 
over the topological vertex partial operad 
${\cal T}^{k}$ of type $k$ 
(or over the strong topological vertex partial operad $\bar{\cal T}^{k}$
of type $k$) 
is called 
{\it (weakly) meromorphic} if it is  (weakly) meromorphic  as 
a ${\Bbb Z}\times {\Bbb Z}$-graded algebra over $\tilde{K}^{0}$ and 
if the map 
{from} $V$ to $\overline{V}$
corresponding to 
$(({\bf 0}, (1, {\bf 0})), f_{T^{k}}(0))\in {\cal T}^{k}(1)$ 
(or $(({\bf 0}, (1, {\bf 0})), f_{\bar{T}^{k}}(0))\in \bar{\cal T}^{k}(1)$) 
defines the fermion grading on $V$. }
\end{defi}

\subsection{Proof of Theorem \ref{main}}

We prove the main theorem, Theorem \ref{main}, of the present paper in this 
subsection using the main theorem of \cite{H} and its generalization 
to locally-grading-restricted
conformal vertex algebra
(that is,  vertex operator algebras without grading-restrictions but satisfying 
the local grading-restriction conditions, see Remark \ref{term})
given in \cite{H7}. 

The case of grading-restricted 
(strong) topological vertex operator algebras of type $k<0$
and the case of  locally-grading-restricted strong topological 
vertex operator algebras of 
type $k$ is an 
easy consequence of the case of locally-grading-restricted  topological 
vertex operator algebras
of type $k$.
We only prove the case of locally-grading-restricted  topological 
vertex operator  
algebras of type $k$.  

Let $V$ be a locally-grading-restricted  topological 
vertex operator algebra of type $k$. Then by 
Proposition \ref{5-1}, the central charge of $V$ is $0$.
In \cite{H7}, 
the main theorem of \cite{H} is generalized
to locally-grading-restricted 
conformal vertex algebras, that is, 
vertex operator algebras without grading-restrictions but satisfying 
the local grading-restriction conditions (see Remark \ref{term}). 
This generalization states 
that the category of locally-grading-restricted 
conformal vertex algebras of central charge $c$ is isomorphic 
to the category of weakly meromorphic algebras over $\tilde{K}^{c}$.
See \cite{H7} for details. 
By this generalization,  $V$ gives a weakly meromorphic  
${\Bbb Z}\times {\Bbb Z}$-graded algebra
over $\tilde{K}^{0}$. {From} the definitions of the partial operad
structure on $\tilde{K}^{0}$ and the definition of a weakly meromorphic 
${\Bbb Z}\times {\Bbb Z}$-graded algebra
over $\tilde{K}^{0}$, we see that $V$ in fact gives 
a ${\Bbb Z}\times {\Bbb Z}$-graded algebra over $K$.

Let $H^{V}$ be the locally convex topological completion of 
$V$ constructed in \cite{H6} and \cite{H7}.
By Proposition \ref{univ-t}, there is a homomorphism of 
locally-grading-restricted
topological 
 vertex operator
algebras {from} $T^{k}$ to $V$.
In particular, elements of $\widetilde{Y(T^{k})}$ can be viewed as
elements of $\hom(V, \overline{V})$. 
Recall that   ${\cal T}^{k}(1)$ is canonically isomorphic to 
$${\frak K}_{1}^{*}({\cal T}^{k}(0)\boxtimes 
{\cal T}^{k}(1)|_{K_{\ge -1}(1)}).$$
So we need only consider elements of the subbundles
${\frak K}_{1}^{*}({\cal T}^{k}(0)\boxtimes 
{\bf 1}_{T^{k}}(0)|_{K_{\ge -1}(1)})$ and 
${\frak K}_{1}^{*}({\frak E}^{*}({\bf 1}_{T^{k}}(0))\boxtimes 
{\cal T}^{k}(1)|_{K_{\ge -1}(1)})$ of ${\cal T}^{k}(1)$.
Any element of ${\cal T}^{k}(1)$ is of the form 
$\hat{Q}=(Q, u\mbar_{Q}),$
 where 
$u\in \widetilde{Y(T^{k})}$. 
(Recall that the fibers of ${\cal T}^{k}(1)$ are 
spanned by the values of the elements of $\widetilde{Y(T^{k})}$.) 
Let $\nu_{n}: K(n)\to \hom(V^{\otimes n}, 
\overline{V})$, $n\ge 0$, be the maps defining the structure of 
a ${\Bbb Z}\times {\Bbb Z}$-graded algebra over $K$. 
We define 
$$\nu_{1}^{{\cal T}^{k}}(\hat{Q})
=e^{-L^{-}(A^{(0)}(Q))}u e^{L^{-}(A^{(0)}(Q))}\nu_{1}(Q)$$
when 
$\hat{Q}\in {\frak K}_{1}^{*}({\cal T}^{k}(0)\boxtimes 
{\bf 1}_{T^{k}}(0)|_{K_{\ge -1}(1)})$,  and we define 
$$\nu_{1}^{{\cal T}^{k}}(\hat{Q})
=\nu_{1}(Q)u$$
when $\hat{Q}\in 
{\frak K}_{1}^{*}({\frak E}^{*}({\bf 1}_{T^{k}}(0))\boxtimes 
{\cal T}^{k}(1)|_{K_{\ge -1}(1)})$.
It is easy to see that when 
$\hat{Q}\in {\frak K}_{1}^{*}({\cal T}^{k}(0)\boxtimes 
{\bf 1}_{T^{k}}(0)|_{K_{\ge -1}(1)})$, $u$ as an element of 
$\hom(H^{V}, H^{V})\subset \hom(V, \overline{V})$
is in fact an
infinite sum of homogeneous elements of $\hom(V, V)$ 
of weights larger than an integer depending on $u$, and when 
$$\hat{Q}\in 
{\frak K}_{1}^{*}({\frak E}^{*}({\bf 1}_{T^{k}}(0))\boxtimes 
{\cal T}^{k}(1)|_{K_{\ge -1}(1)}),$$
$u$ as an element of 
$\hom(H^{V}, H^{V})\subset \hom(V, \overline{V})$ is in fact in 
$\hom(V, V)$. Thus in both cases,
$\nu_{1}^{{\cal T}^{k}}(\hat{Q})$ is indeed in 
$\hom(V, \overline{V})$.

We also define 
$$\nu_{2}^{{\cal T}^{k}}({\frak
K}_{2}^{*}({\frak E}^{*}({\bf 1}_{T^{k}}(0))\boxtimes 
{\bf 1}_{T^{k}}(0)\boxtimes {\bf 1}_{T^{k}}(0)\mbar_{P(z)})
=\nu_{2}(P(z)).$$
Note that ${\cal T}^{k}(0)$ is constructed {from} ${\cal T}^{k}(1)$ and
 ${\cal T}^{k}$ is generated by elements of ${\cal T}^{k}(0)$,
$\hat{Q}$ and 
$${\frak
K}_{2}^{*}({\frak E}^{*}({\bf 1}_{T^{k}}(0))\boxtimes 
{\bf 1}_{T^{k}}(0)\boxtimes {\bf 1}_{T^{k}}(0)\mbar_{P(z)}.$$ 
Also note that 
the definition of the composition maps for
${\cal T}^{k}$ depends only on the vertex operator algebra structure on 
$T^{k}$. Thus
we can extend  $\nu_{1}^{{\cal T}^{k}}$ and 
$\nu_{2}^{{\cal T}^{k}}$ defined 
on these particular elements of ${\cal T}^{k}$ to a morphism 
$\nu^{{\cal T}^{k}}$ of 
partial pseudo-operads {from} ${\cal T}^{k}$ to the endomorphism 
partial pseudo-operad ${\cal H}_{V}=\{\hom(V^{\otimes n}, 
\overline{V})\}_{n\ge 0}$ of $V$. So $V$ gives a ${\Bbb Z}\times 
{\Bbb Z}$-graded 
algebra over ${\cal T}^{k}$.
Since as a ${\Bbb Z}\times {\Bbb Z}$-graded algebra over 
$\tilde{K}^{0}$, it is meromorphic and since
the map {from} $V$ to $\overline{V}$
corresponding to
$({\bf 0}, (1, {\bf 0}), f_{T^{k}}(0))\in {\cal T}^{k}(1)$ 
is the fermion grading operator 
$f_{0}$ 
defining the fermion grading on $V$, the ${\Bbb Z}\times 
{\Bbb Z}$-graded algebra over ${\cal T}^{k}$
constructed above
is weakly meromorphic. We define a functor
{from} the category of locally-grading-restricted topological 
vertex operator algebras of type $k$
to the category of weakly meromorphic 
${\Bbb Z}\times {\Bbb Z}$-graded algebras over ${\cal T}^{K}$ 
by assigning to 
a locally-grading-restricted topological 
vertex operator 
algebra $V$ of type $k$ the weakly meromorphic 
${\Bbb Z}\times {\Bbb Z}$-graded algebra over ${\cal T}^{k}$ 
constructed above {from} $V$.

Conversely, given a weakly meromorphic 
${\Bbb Z}\times {\Bbb Z}$-graded algebra $V$ over ${\cal T}^{k}$, 
we obtain a locally-grading-restricted topological  
vertex operator algebra as follows:
Since by definition 
$V$ is a weakly meromorphic ${\Bbb Z}\times {\Bbb Z}$-graded 
algebra over $\tilde{K}^{0}$,
$V$ has a structure of locally-grading-restricted 
${\Bbb Z}\times {\Bbb Z}$-graded 
conformal  vertex algebra of central charge $0$ by the generalization of 
the main theorem of \cite{H} discussed above and in \cite{H7}. 
Let $\nu_{n}^{{\cal T}^{k}}$, $n\ge 0$, be the maps {from} 
${\cal T}^{k}(n)$ to 
$\hom(V^{\otimes n}, 
\overline{V})$ defining the structure of a weakly 
meromorphic ${\Bbb Z}\times {\Bbb Z}$-graded algebra over 
${\cal T}^{k}$. Let 
\begin{eqnarray*}
g=\nu_{0}^{{\cal T}^{k}}(({\bf 0}, {\frak E}^{*}(g_{T^{k}}(-2))
\mbar_{\bf 0})),\\
f=\nu_{0}^{{\cal T}^{k}}(({\bf 0}, {\frak E}^{*}(f_{T^{k}}(-1))
\mbar_{\bf 0})),\\
q=\nu_{0}^{{\cal T}^{k}}(({\bf 0}, {\frak E}^{*}(q_{T^{k}}(-1))
\mbar_{\bf 0})).
\end{eqnarray*}
Then it is clear that $V$ together 
with the elements $g$, $f$ and $q$ is a locally-grading-restricted 
topological vertex operator
algebra. {From} the construction of $V$, we also see that there is 
a canonical homomorphism of locally-grading-restricted 
topological vertex operator
algebras {from} $T_{k}$ to $V$. So $V$ is of type $k$.
Thus we obtain a functor {from} the category of weakly meromorphic 
${\Bbb Z}\times {\Bbb Z}$-graded algebras over ${\cal T}^{k}$ to the 
category of 
locally-grading-restricted topological 
 vertex operator algebras of type $k$. 

The two functors constructed above are clearly inverse to each other. 
Thus we see that the two categories are isomorphic.

\subsection{Proof of Theorem {\ref{theo2}}}

In this subsection, we prove Theorem \ref{theo2} using Theorem \ref{main}.

Let $V$ be a 
locally-grading-restricted strong topological vertex operator algebra.
First we know that there exists $k<0$ such that 
$V$ is of type $k$.
By Theorem \ref{main}, $V$ has a structure of a 
weakly meromorphic ${\Bbb Z}\times {\Bbb Z}$-graded
algebra over $\bar{\cal T}^{k}$. In particular, $V$ has a structure of a
weakly meromorphic ${\Bbb Z}\times {\Bbb Z}$-graded algebra over $K$. 
Let $H^{V}$ be the locally convex topological completion
of $V$ constructed in  \cite{H6} and
\cite{H7}. Then it has been shown in \cite{H6} and
\cite{H7} that
the maps in $\hom(V^{\otimes n}, \overline{V})$
corresponding to elements of $K_{{\frak H}_{1}}(n)$ 
is in fact in $\hom(H^{\otimes n}, H)$. 
Moreover, $H^{V}$ has
a structure of an algebra over the operad (not partial) $K_{{\frak H}_{1}}$.

For $n\ge 0$, let $\bar{\cal T}^{k}_{{\frak H}_{1}}(n)$ be the
holomorphic subbundle over $K_{{\frak H}_{1}}(n)$
of the holomorphic bundle $\bar{\cal T}^{k}(n)$
consisting elements of the form $(Q, \alpha)$ where 
$Q\in K_{{\frak H}_{1}}(n)$ and $\alpha$ is the restriction to $Q$ of 
a linear combination of sections
$${\frak K}_{n_{1}}^{*}(u^{(0)}\otimes
u^{(1)}\otimes \cdots\otimes u^{(n_{1})})$$
for $u^{(0)}\in {\frak E}^{*}(\widetilde{Y(\bar{T}^{k})}\cap 
\hom(H^{\bar{T}^{k}}, H^{\bar{T}^{k}}))$
satisfying 
$$e^{-L^{-}(A^{(0)}(Q))}u^{(0)}e^{L^{-}(A^{(0)}(Q))}
\in \hom(H^{\bar{T}^{k}}, H^{\bar{T}^{k}})$$
and
$u^{(k)}\in (\widetilde{Y(\bar{T}^{k})}\cap \hom(H^{\bar{T}^{k}}, 
H^{\bar{T}^{k}}))|_{K_{\ge -1}(1)}$, 
$k=1, \dots, n$. Then it is clear
that $\bar{\cal
T}^{k}_{{\frak H}_{1}}=\{\bar{\cal
T}^{k}_{{\frak H}_{1}}(n)\}_{n\ge 0}$ is a suboperad (not partial)
of $\bar{\cal T}^{k}$.
{From} the definition of $\bar{\cal T}^{k}$ in Subsection 4.2
and the construction of the structure of a weakly meromorphic ${\Bbb
Z}\times {\Bbb Z}$-graded algebra over $\bar{\cal T}^{k}$ on $V$
in the preceding subsection, we
see that the maps in $\hom(V^{\otimes n}, \overline{V})$ corresponding
to elements of $\bar{\cal T}^{k}_{{\frak H}_{1}}(n)$ 
are actually in $\hom((H^{V})^{\otimes n}, H^{V})$. Thus we see 
that $H^{V}$ has
a structure of an algebra over the operad (not partial) $\bar{\cal
T}^{k}_{{\frak H}_{1}}$. 

On the other hand, 
{from} the construction of
$\bar{\cal T}^{k}$ in Subsection 4.2, there is a morphism of operads
{from} $\wedge TK_{{\frak H}_{1}}$ to $\bar{\cal T}^{k}_{{\frak
H}_{1}}$ constructed as follows: Consider the  vertex operator
subalgebra $T_{b}$ of $T_{G}$ generated by $b$ and $\omega_{G}$. Using
$T_{b}$ and the
same construction as the one used to construct ${\cal T}_{G}$ at the
end of Subsection 3.3, we obtain a partial suboperad  ${\cal T}_{b}$
of ${\cal T}_{G}$. Let $\wedge TK(n)$, $n\ge 0$,
be  the direct sums of all wedge powers of the tangent bundles
of the moduli spaces $K(n)$. Then it is clear that 
$\wedge TK=\{\wedge TK(n)\}_{n\ge 0}$
is a partial operad. For any $n\ge 0$, the direct sums of all wedge 
powers of the holomorphic tangent bundles of the moduli space $K(n)$
is a subbundle of $\wedge TK(n)$. In fact it is clear that 
this subbundle is isomorphic to the bundle ${\cal T}_{b}(n)$ and
there is a bundle map {from} $\wedge TK(n)$ to ${\cal T}_{b}(n)$ such 
on each fiber, this bundle map is the projection {from} the fiber 
of $\wedge TK(n)$ to the fiber of ${\cal T}_{b}(n)$. 
Moreover, by the definition of the composition maps for ${\cal T}_{b}$,
we see that ${\cal T}_{b}$ is in fact a partial suboperad of 
$\wedge TK$ and the bundle maps {from} $\wedge TK(n)$ to ${\cal T}_{b}(n)$
give a morphism of partial operads. 

Clearly the restriction of $\wedge TK$ to
$K_{{\frak H}_{1}}(n)$ is the operad $\wedge TK_{{\frak H}_{1}}$.
Since $\bar{T}$ is a strong topological vertex operator algebra, 
there is a homomorphism 
of vertex operator algebras {from} $T_{b}$ to $\bar{T}$. 
This homomorphism induces a morphism of partial operad {from}
${\cal T}_{b}$ to $\bar{\cal T}^{k}$. Composing this 
morphism with the morphism {from} $\wedge TK$ to ${\cal T}_{b}$,
we obtain a morphism {from} $\wedge TK$ to $\bar{\cal T}^{k}$.
This morphism gives a morphism
{from} $\wedge TK_{{\frak H}_{1}}$ to $\bar{\cal T}^{k}_{{\frak
H}_{1}}$. 
Composing the morphism {from} $\wedge TK_{{\frak H}_{1}}$ to 
$\bar{\cal T}^{k}_{{\frak
H}_{1}}$  and the morphism {from} $\bar{\cal T}^{k}_{{\frak
H}_{1}}$ to $\{\hom((H^{V})^{\otimes n},
H^{V})\}_{n\ge 0}$, we obtain a morphism {{from}} $\wedge
TK_{{\frak H}_{1}}$ to ${\cal E}_{H^{V}}=\{\hom((H^{V})^{\otimes n},
H^{V})\}_{n\ge 0}$.  

The differential $Q$ on $H^{V}$ is the natural extension to $H^{V}$ of
the operator $Q=q_{0}=\res_{x}Y(q, x)$ on $V$.  The proof of
(\ref{6.1}) is the same as the proof of the same formula for
topological vertex algebras in \cite{H2} (cf. also \cite{KSV}). In
fact, since in our case $b(j)$, $j\in {\Bbb Z}$, are anti-commutative
with each other, the proof of (\ref{6.1}) is simpler than the
corresponding proof in  \cite{H2}.

\renewcommand{\theequation}{\thesection.\arabic{equation}}
\renewcommand{\therema}{\thesection.\arabic{rema}}
\setcounter{equation}{0}
\setcounter{rema}{0}
\section{Appendix: Examples of 
 locally-grading-re\-stricted strong  topological 
vertex operator algebras}

\subsection{Locally-grading-restricted conformal vertex algebra 
of central charge $26$ tensored with the ghost vertex operator algebra}

These examples  are substructures of string backgrounds in 
string theory. They were discussed extensively by Lian and
Zuckerman. See for example \cite{LZ2} and the references there. 

Let $V$ be a locally-grading-restricted conformal vertex algebra 
of central charge $26$, that is, a vertex operator algebra of without 
grading restrictions of central charge $26$
satisfying local grading-restriction conditions (see Remark \ref{term}). Then 
$V\otimes \wedge_{\infty} \psi_{0}^{*}(\tilde{\frak D}_{*}({\Bbb L}))$
is a locally-grading-restricted 
${\Bbb Z}\times {\Bbb Z}$-graded conformal vertex algebra  of central charge $0$.
Let
\begin{eqnarray*}
g&=&{\bf 1}_{V}\otimes b,\\
q&=&L_{V}(-2){\bf 1}\otimes c(1)1
+1\otimes (b(-2)c(1)c(0)1)\\
&=&L_{V}(-2){\bf 1}\otimes c+{\bf 1}\otimes 
(\psi_{0}^{*}({\Bbb L}'(0))\wedge \psi_{0}^{*}({\Bbb L}'(-1))\wedge 
\psi_{0}^{*}(\tilde{\frak D}_{*}({\Bbb L}(-2)))),\\
f_{G}&=&{\bf 1}\otimes (c(1)b(-2)1)\\
&=&{\bf 1}\otimes (\psi_{0}^{*}({\Bbb L}'(-1))\wedge 
\psi_{0}^{*}(\tilde{\frak D}_{*}({\Bbb L}(-2))))
\end{eqnarray*}
be three elements of 
$V\otimes \wedge_{\infty} \psi_{0}^{*}(\tilde{\frak D}_{*}({\Bbb L}))$.
Then $V\otimes \wedge_{\infty} \psi_{0}^{*}(\tilde{\frak D}_{*}({\Bbb
L}))$ together with these three elements is a 
locally-grading-restricted strong  topological
vertex operator algebra. 

The above construction begin with a locally-grading-restricted 
conformal vertex algebra of central charge $26$.
Such an algebra can be constructed {from} a 
locally-grading-restricted 
conformal vertex algebra of arbitrary central charge $c$
this algebra with a familiar locally-grading-restricted 
conformal vertex algebra
of central charge $26-c$. One example of such familiar 
locally-grading-restricted 
conformal vertex  algebras 
is the one coming {from} the so-called Liouville models. We can also use
other familiar examples associated to the Virasoro algebra, affine Lie
algebras and lattices. 

One important example of this class of 
 locally-grading-restricted strong topological
vertex operator  algebra is the 
tensor product of the moonshine module vertex operator algebra constructed
by Frenkel, Lepowsky and Meurman \cite{FLM1} \cite{FLM}, 
the locally-grading-restricted 
conformal vertex algebra constructed {from} a rank $2$ Lorentz lattice
and the ghost vertex operator algebra $\wedge_{\infty} 
\psi_{0}^{*}(\tilde{\frak D}_{*}({\Bbb L}))$. This example was studied
in detail by Lian and Zuckerman in \cite{LZ3}. Since the lattice 
used is Lorentzian, this example does not satisfy the 
grading-restriction conditions. So it is not a topological vertex 
operator algebra in the sense of \cite{H2} and Definition \ref{tvoa}. 
The tensor product of the first two algebras
was used by Borcherds in his famous 
proof of the Monstrous Moonshine Conjecture in
\cite{Bor2}. The Monster Lie algebra constructed in \cite{Bor2} 
can also be constructed using the Gerstenhaber algebra structure 
on the cohomology of this example.

\subsection{Twisted $N=2$ superconformal vertex operator superalgebras}

In \cite{EY}, Eguchi and Yang 
constructed explicitly  topological conformal field 
theories in a physical sense
by  twisting $N=2$ superconformal field theories. 
A good exposition can be found in, for example,  \cite{Wa} by
Warner. Here we formulate the 
theory in
terms of the mathematical language of vertex operator algebras. These
topological vertex operator algebras are closely related to the mirror
symmetry. 

\begin{defi}
{\rm An {\it $N=2$ superconformal vertex operator superalgebra} is a
vertex operator superalgebra $(V, Y, {\bf 1}, \omega)$ together with
odd elements $\tau^{+}$, $\tau^{-}$ and even $h$ satisfying the
following axioms: Let
\begin{eqnarray*}
Y(\tau^{+}, x)&=&\sum_{n\in {\Bbb Z}}G^{+}(n+\frac{1}{2})x^{-n-2},\\
Y(\tau^{-}, x)&=&\sum_{n\in {\Bbb Z}}G^{-}(n-\frac{1}{2})x^{-n-1},\\
Y(h, x)&=&\sum_{n\in {\Bbb Z}}J(n)x^{-n-1}.
\end{eqnarray*}
Then $V$ is a direct sum of eigenspaces of $J(0)$ with integral
eigenvalues which modulo $2{\Bbb Z}$ give the ${\Bbb Z}_{2}$ grading
for the vertex operator superalgebra structure, and  the following 
$N=2$ Neveu-Schwarz relations hold: For $m, n\in {\Bbb Z}$,
\begin{eqnarray*}
{[L(m), L(n)]}&=&(m-n)L(m+n)+\frac{c}{12}(m^{3}-m)\delta_{m+n, 0},\\
{[J(m), J(n)]}&=&\frac{c}{3}m\delta_{m+n, 0},\\
{[L(m), J(n)]}&=&-nJ_{m+n},\\
{[L(m), G^{\pm}_{n\pm \frac{1}{2}}]}&=&\left(\frac{m}{2}-(n\pm
\frac{1}{2})\right)G^{\pm}_{m+n\pm \frac{1}{2}},\\
{[J(m), G^{\pm}_{n\pm\frac{1}{2}}]}&=&\pm
G^{\pm}_{m+n\pm\frac{1}{2}},\\
{[G^{+}_{m+\frac{1}{2}}, G^{-}_{n-\frac{1}{2}}]}&=& 
2L(m+n)+(m-n+1)J(m+n)+\frac{c}{3}(m^{2}+m)\delta_{m+n, 0},\\
{[G^{\pm}_{m\pm \frac{1}{2}}, G^{\pm}_{n\pm \frac{1}{2}}]}&=&0
\end{eqnarray*}
where $L(m)$, $m\in {\Bbb Z}$, are the Virasoro operators on $V$ and
$c$ is the central charge of $V$.}
\end{defi}

We shall denote the $N=2$ superconformal vertex operator superalgebra
defined above by 
$(V, Y, {\bf 1}, \omega, \tau^{+}, \tau^{-}, h)$ or simply by $V$.

Let $(V, Y, {\bf 1}, \omega, \tau^{+}, \tau^{-}, h)$
be an $N=2$ superconformal vertex operator superalgebra.  
We define
\begin{eqnarray*}
\omega_{T}&=&\omega+\frac{1}{2}J(-2){\bf 1}, \\
f&=&h\\
q&=&\tau^{+},\\
g&=&\tau^{-}.
\end{eqnarray*}

An easy calculation gives the following:

\begin{propo}
Let $(V, Y, {\bf 1}, \omega, \tau^{+}, \tau^{-}, h)$
be an $N=2$ superconformal vertex operator superalgebra.  Then
$(V, Y, {\bf 1}, \omega_{T}, f, q, g)$ is a grading-restricted
strong topological vertex
operator algebra. \epf
\end{propo}

{\small \sc Department of Mathematics, Rutgers University,
110 Frelinghuysen Rd., Piscataway, NJ 08854-8019}

{\em E-mail address}: yzhuang@math.rutgers.edu

\vspace{2em}

{\small \sc Department of Mathematics, University of Chicago, 
Chicago, IL 60637}

{\em E-mail address}: zhao@math.uchicago.edu

\end{document}